\documentclass[a4paper, reqno]{amsart}

\title[{Twisted Zhu algebras}]{Twisted Zhu algebras}

\author{Naoki Genra}
\address{Graduate School of Mathematical Sciences, The University of Tokyo 3-8-1 Komaba, Meguro, 153-8914 Tokyo, Japan}
\email{genra@ms.u-tokyo.ac.jp}

\usepackage{latexsym}
\usepackage{amsmath}
\usepackage{amsfonts}
\usepackage{amssymb}
\usepackage{amsxtra}
\usepackage{mathrsfs}
\usepackage[all]{xy}
\usepackage{mathtools}
\usepackage{enumerate}
\usepackage{enumitem}
\usepackage[dvipdfmx]{graphicx,color}
\usepackage[dvipdfmx]{hyperref}

\hypersetup{% hyperref_options
setpagesize=false,
 bookmarksnumbered=true,%
 bookmarksopen=true,%
 colorlinks=flase,%
 linkcolor=red,
 citecolor=green,
 urlcolor=cyan,
}

\newtheorem{definition}{Definition}[section]
\newtheorem{proposition}[definition]{Proposition}
\newtheorem{theorem}[definition]{Theorem}
\newtheorem{corollary}[definition]{Corollary}
\newtheorem{lemma}[definition]{Lemma}

\newtheorem*{thmA}{Theorem A}
\newtheorem*{thmB}{Theorem B}
\newtheorem*{thmC}{Theorem C}
\newtheorem*{thmD}{Theorem D}
\newtheorem*{thmE}{Theorem E}

\numberwithin{equation}{section}

\newcommand{\N}{\mathbb{N}}
\newcommand{\Z}{\mathbb{Z}}
\newcommand{\Q}{\mathbb{Q}}
\newcommand{\R}{\mathbb{R}}
\newcommand{\C}{\mathbb{C}}
\newcommand{\End}{\operatorname{End}}

\newcommand{\Ker}{\operatorname{Ker}}
\newcommand{\Img}{\operatorname{Im}}
\newcommand{\Span}{\operatorname{Span}}
\newcommand{\id}{\operatorname{id}}
\newcommand{\gr}{\operatorname{gr}}
\newcommand{\ad}{\operatorname{ad}}

\newcommand{\NO}[1]{:\!#1\!:}
\newcommand{\vac}{|0\rangle}

\newcommand{\Zhu}{\operatorname{Zhu}}
\newcommand{\base}{\mathfrak{v}}
\newcommand{\conf}{R}
\newcommand{\monoid}{\Gamma}
\newcommand{\gH}[1]{#1_{g,H}}
\newcommand{\Charge}{\Omega}
\newcommand{\Coh}[1]{\overline{#1}}
\newcommand{\Cohindex}{\overline{I}}
\newcommand{\wt}{\operatorname{wt}}
\newcommand{\projZhu}{\tau}
\newcommand{\twist}{\theta_H}

\begin{document}
\maketitle

\markboth{Twisted Zhu algebras}{Twisted Zhu algebras}

\begin{abstract}
Let $V$ be a freely generated pregraded vertex superalgebra, $H$ a Hamiltonian operator of $V$, and $g$ a diagonalizable automorphism of V commuting with $H$ with modulus $1$ eigenvalues. We prove that the $(g, H)$-twisted Zhu algebra of $V$ has a PBW basis, is isomorphic to the universal enveloping algebra of some non-linear Lie superalgebra, and satisfies the commutativity of BRST cohomology functors, which generalizes results of De Sole and Kac. As applications, we compute the twisted Zhu algebras of affine vertex superalgebras and affine $W$-algebras.
\end{abstract}

\section{Introduction}

Let $V$ be a vertex superalgebra, $H$ a Hamiltonian operator whose eigenvalues are real numbers, and $g$ an automorphism of $V$ commuting with $H$ whose eigenvalues are modulus $1$, i.e., in $S^1=\{z \in \C \mid |z|=1\}$. The $(g, H)$-twisted Zhu algebra of $V$ is defined by
\begin{align*}
\gH{\Zhu}V = \gH{V}/\gH{J}(V),
\end{align*}
where $\gH{V}= V^{ g \circ \theta_H^{-1}}$, $\theta_H = \mathrm{e}^{2\pi i H}$ is the twist automorphism \cite{HKL}, and $\gH{J}(V)$ is a subspace defined by \eqref{eq:def-VandJ}. Then $\gH{\Zhu}V$ has a structure of a unital associative superalgebra by the product induced by the $*$-product
\begin{align*}
a * b = \sum_{j=0}^\infty\binom{\Delta_a}{j}a_{(j-1)}b,\quad
a, b \in \gH{V},
\end{align*}
where $\Delta_a$ is the eigenvalue of $a$ for $H$, called the \textit{conformal weight} of $a$. The twisted Zhu algebras enjoy the Zhu theorem: there exists a one-to-one correspondence between simple lower-bounded $g$-twisted $V$-modules and simple $\gH{\Zhu}V$-modules, originally proved by Zhu \cite{Zhu} when $V$ is a vertex operator algebra and $g=1$, and generalized by \cite{KWan, DLM, DZ, DK}.
\smallskip

Suppose that $V$ is freely generated and pregraded in the sense of \cite{DK}. See also Section \ref{sec:def-va} for the definitions. In the present paper, we prove that the $(g, H)$-twisted Zhu algebra of $V$ has a PBW basis, is isomorphic to the universal enveloping algebra of some non-linear Lie superalgebra, and satisfies the commutativity of BRST cohomology functors, which generalizes results of De Sole and Kac in \cite{DK}. As applications, we prove that the twisted Zhu algebras of affine vertex superalgebras $V^k(\mathfrak{g})$ are isomorphic to the universal enveloping algebras of certain Lie subalgebras of $\mathfrak{g}$, and that the twisted Zhu algebras of affine $W$-algebras $W^k(\mathfrak{g}, f)$ are isomorphic to the finite $W$-algebras of certain Lie subalgebras of $\mathfrak{g}$ with the same nilpotent element $f$.

To clarify our results, let $\{e_\alpha\}_{\alpha \in I}$ be free generators of $V$ homogeneous for $g, H$ and parity with an ordered index set $I$, and $\tau \colon \gH{V} \twoheadrightarrow \gH{\Zhu}V$ be the canonical projection. Set
\begin{align*}
&\gH{\Zhu} \conf := \projZhu(\gH{\conf}),\quad
\gH{\conf} = \C[\partial]\otimes\gH{\base} \subset \gH{V},\\
&\gH{\base} = \Span_\C\{e_\alpha \mid \alpha \in I,\ (g \circ \theta_H^{-1})(e_\alpha) = e_\alpha\},
\end{align*}
where $\partial$ is the translation operator. 

\begin{thmA}[Theorem \ref{thm:main PBW thm}]
For any $g$ and $H$, $\gH{\Zhu} \conf$ has a structure of non-linear Lie superalgebra such that
\begin{align*}
\gH{\Zhu} V \simeq U(\gH{\Zhu} \conf).
\end{align*}
In particular, $\gH{\Zhu} V$ is an associative superalgebra generated by $\projZhu(\gH{\base}) \simeq \gH{\base}$ with the defining relations
\begin{align*}
[\projZhu(a), \projZhu(b)] = \sum_{j=0}^{\infty}\binom{\Delta_a-1}{j}\projZhu(a_{(j)}b),\quad
a, b \in \gH{\base}.
\end{align*}
\end{thmA}

See \cite[Section 3]{DK}, or Section \ref{sec:non-lin Liesalg} for the definitions of non-linear Lie superalgebras $\base$ and the universal enveloping algebra $U(\base)$. De Sole and Kac proved Theorem A in \cite[Theorem 3.25]{DK} in the case of $g = \theta_H$, called the \textit{$H$-twisted} case, and also raised a question about the existence of a non-linear Lie superalgebra $\gH{\Zhu} \conf$ that realizes the twisted Zhu algebra $\gH{\Zhu} V$ as the universal enveloping algebra in \cite[Remark 3.31]{DK}. Theorem A gives a positive answer to the problem.
\smallskip

Suppose that $V$ satisfies the same assumptions in Theorem A with a $\Z$-grading $V= \bigoplus_{n \in \Z}V^n$, and that $(V^\bullet, d)$ forms a cochain complex by an odd differential $d \colon V^n \rightarrow V^{n+1}$. Then the cohomology $H(V, d)$ is a vertex superalgebra. If $d$ commutes with $g$ and $H$, $H(V, d)$ has the induced Hamiltonian operator $H$ and the induced automorphism $g$, and $d$ induces an odd differential $\overline{d}$ on $\gH{\Zhu}V$.

\begin{thmB}[Theorem \ref{thm:main_cohomology_thm}]
Suppose (A.1)--(A.7) in Section \ref{sec:commutativity}. Then
\begin{align*}
\gH{\Zhu}H(V, d) \simeq H(\gH{\Zhu}V, \overline{d})
\end{align*}
as associative superalgebras.
\end{thmB}

Theorem B is a slight generalization of \cite[Theorem 4.20]{DK} and is essentially proved in the same way as in that paper using Theorem A.
\smallskip

Let $\mathfrak{g}$ be a simple basic classical Lie superalgebra with the normalized even supersymmetric invariant bilinear form $(\cdot|\cdot)$, $V^k(\mathfrak{g})$ the (universal) affine vertex superalgebra of $\mathfrak{g}$ at level $k \in \C$, $L^\mathfrak{g}$ the Sugawara conformal vector of $V^k(\mathfrak{g})$, and $x$ a semisimple element of $\mathfrak{g}$ whose eigenvalues are real numbers. Set a Hamiltonian operator $H$ on $V^k(\mathfrak{g})$ by
\begin{align}\label{intro:H}
H = (L^\mathfrak{g} + \partial x)_0 = L^\mathfrak{g}_0 - x_{(0)}.
\end{align}
The twist $\theta_H$ induces an automorphism of $\mathfrak{g}$ by
\begin{align}
\theta_H(a) = \mathrm{e}^{2\pi i \Delta_a}a = \mathrm{e}^{-2 \pi i j}a,\quad
a \in \mathfrak{g}_j = \{ u \in \mathfrak{g} \mid [x, u] = j u\}.
\end{align}
Let $g$ be a diagonalizable Lie superalgebra automorphism of $\mathfrak{g}$ with modulus $1$ eigenvalues. Assume that $g(x) = x$ and $(g(a)|g(b)) = (a|b)$ for $a, b \in \mathfrak{g}$. The first condition guarantees $g \circ \ad x = \ad x \circ g$. The second condition is automatically true if the dual Coxeter number is non-zero since the normalized form is proportional to the Killing one. Then $g$ induces an automorphism of $V^k(\mathfrak{g})$ commuting with $H$.

\begin{thmC}[Theorem \ref{thm:twisted Zhu of affine}]
$\Zhu_{g, H}V^k(\mathfrak{g}) \simeq U(\mathfrak{g}^{g\circ \twist^{-1}})$ for all $g$ and $H$.
\end{thmC}

Theorem C is an application of Theorem A. Theorem C was proved by \cite{FZ} when $\mathfrak{g}$ is a Lie algebra,  $x=0$ and $g=1$, by \cite{Yang} when $\mathfrak{g}$ is a Lie algebra and $x=0$ but general $g$, and by \cite{SRW} when $\mathfrak{g} = \mathfrak{osp}(1|2)$, $x=0$ and $g=1$ or the parity automorphism.

Let $f$ be a nilpotent element with the even parity in $\mathfrak{g}$. Choose a semisimple element $x$ in $\mathfrak{g}$ whose eigenvalues are half integers and satisfying that $(f, x)$ is a good pair in the sense of \cite{KRW}. See, e.g., Section \ref{sec:good} for the definitions of goodness, and \cite{EK, Hoyt} for the classifications. Then the (universal) affine $W$-algebra $W^k(\mathfrak{g}, f)$ of $\mathfrak{g}, f$ at level $k$ is defined as the generalized Drinfeld-Sokolov reduction of $V^k(\mathfrak{g})$ by \cite{FF90a, KRW}. The Hamiltonian operator \eqref{intro:H} on $V^k(\mathfrak{g})$ may be extended to the vertex superalgebra $C^k(\mathfrak{g}, f, x)$ of the cochain complex and induces the Hamiltonian operator of the affine $W$-algebra
\begin{align*}
W^k(\mathfrak{g}, f) = H(C^k(\mathfrak{g}, f, x), d).
\end{align*}
Let $g$ be a diagonalizable Lie superalgebra automorphism of $\mathfrak{g}$ with modulus $1$ eigenvalues such that $g(x) = x$, $g(f) = f$ and $(g(a)|g(b)) = (a|b)$ for $a, b \in \mathfrak{g}$. Then $g$ induces an automorphism of $W^k(\mathfrak{g}, f)$ commuting with $H$. Denote by $U(\mathfrak{g}, f)$ the \textit{finite $W$-algebra} of $\mathfrak{g}$, $f$, introduced by \cite{Premet02, Premet07, Kostant, Lynch, BT, RaSo, GG, Losev}. As an application of Theorem B, we have the following results, proved by \cite{KWan, Milas, Arakawa07, DK} in special cases.

\begin{thmD}[Theorem \ref{thm:ZhuofW}, Corollary \ref{cor: g=1-finite-W}]
$\Zhu_{g, H}W^k(\mathfrak{g}, f) \simeq U(\mathfrak{g}^{g\circ \twist^{-1}}, f)$ for all $g$. In particular,
\begin{align*}
\Zhu_{1, H}W^k(\mathfrak{g}, f) \simeq U(\mathfrak{g}_\Z, f),
\end{align*}
where $\mathfrak{g}_\Z = \{ a \in \mathfrak{g} \mid [x, a]=j a\ \mathrm{for}\ \mathrm{some}\ j \in \Z\}$.
\end{thmD}

Let $V$ be a vertex superalgebra with $g$ and $H$ as before. Then we have a Poisson superalgebra analog of the $(g, H)$-twisted Zhu algebra of $V$, denoted by $\gH{R}(V)$, introduced by Zhu \cite{Zhu} in the untwisted cases, and generalized by \cite{DK}. Then, there exists the canonical projection
\begin{align}\label{intro:cano-proj}
\gH{R}(V) \twoheadrightarrow \operatorname{gr}\gH{\Zhu}V,
\end{align}
where $\operatorname{gr}\gH{\Zhu}V$ is the associated graded Poisson superalgebra of $\gH{\Zhu}V$ for the filtration induced by the Li filtration on $V$ \cite{Li}. Using Theorem A, we have the following results, proved in \cite[Proposition 3.30]{DK} in the $H$-twisted cases.

\begin{thmE}[Theorem \ref{thm:R_surj}]
If $V$ satisfies the assumptions in Theorem A, \eqref{intro:cano-proj} is an isomorphism of Poisson superalgebras.
\end{thmE}

Let $U$ be a quotient vertex superalgebra of $V$ by a $(g, H)$-invariant ideal, and $\gH{X}(U) = \operatorname{Specm}\gH{R}(U)$ the associated Poisson variety. If $g = \theta_H$, $X_U := {X}_{\theta_H, H}(U)$ is called the \textit{associated variety of $U$}, introduced by \cite{Arakawa12}. The $(g \circ \theta_H^{-1})$-action on $V$ induces a $\C^*$-action on $X_V$ and the $\C^*$-invariant space is isomorphic to $\gH{X}(V)$ under the assumptions in Theorem A. Then
\begin{align*}
X_V \simeq \gH{X}(V) \sqcup (X_V)_{\neq0}.
\end{align*}
We prove that $\dim\gH{X}(U) = 0$ if $X_U$ is \textit{lisse along $(X_V)_{\neq0}$} in the sense of \cite{BFM} (Proposition \ref{prop:lisse along}). In our future work, we hope to study the relationship between the $g$-twisted $U$-modules and the associated Poisson variety $\gH{X}(U)$.
\smallskip

The paper is organized as follows. In Section \ref{sec:def-va}, we recall the definitions of freely generated pregraded vertex superalgebras. In Section \ref{sec:non-lin Liesalg} and \ref{sec:non-linear lie conf salg}, we review the definitions of non-linear Lie superalgebras and non-linear Lie conformal superalgebras and results of \cite{DK05, DK} about the universalities. In Section {sec:twisted modules} and {sec:Gamma/Z}, we recall the definitions of the $g$-twisted lower-bounded $V$-modules, the $*_n$-products, the $(g, H)$-twisted Zhu algebras, and the Zhu theorem. In Section \ref{sec:properties_*}, we collect properties of the $*_n$-products stated in \cite{DK} and prove them for completeness. In Section \ref{sec:main}, we prove Theorem A. The idea of the proof is almost the same as in \cite{DK}: we introduce a subset $\Sigma$ of $\gH{J}(V)$ consisting of linearly independent elements and show that $\Sigma$ is a basis of $\gH{J}(V)$. In Section \ref{sec:asso var}, we introduce the Poisson superalgebras $\gH{R}(V)$, prove Theorem E by using Theorem A, and also study the lisse conditions of the associated variety $\gH{X}(V) = \operatorname{Specm}\gH{R}(V)$. In Section \ref{sec:coh of non-linear}, we review (and give a slight generalization of) results of \cite{DK} about sufficient conditions for the vanishing of cohomologies of cochain complexes of non-linear Lie superalgebras and non-linear Lie conformal superalgebras. In Section \ref{sec:commutativity}, we prove Theorem B by using Theorem A and the results of Section \ref{sec:coh of non-linear}. In Section \ref{sec:applications}, we prove Theorem C and D by using Theorem A and B, respectively.
\smallskip

{\it Acknowledgments}\quad The author is deeply grateful to Hiroshi Yamauchi for valuable comments and crucial suggestions that generalize his results about the twisted Zhu algebras of the principal affine $W$-algebras of $\mathfrak{osp}(1|2n)$ to Theorem A. He wishes to thank Shigenori Nakatsuka, Tomoyuki Arakawa, Anne Moreau, Thomas Creutzig, Masahiko Miyamoto, Robert McRae, Xingjun Lin, and Qing Wang for valuable comments and discussions. He also thanks Drazen Adamovic, David Ridout, Ching Hung Lam, Jinwei Yang, Vyacheslav Futorny, and Ana Kontrec for organizing great conferences and providing him with meaningful encounters with various researchers and problems. Some part of this work was done while he was visiting Inter-University Centre, Dubrovnik, Croatia in October 2022, Centre de Recherches Math\'{e}matiques, Universit\'{e} de Montr\'{e}al, Qu\'{e}bec, Canada in October and November 2022, Institute of Mathematics, Academia Sinica, Taipei, Taiwan in December 2022, Research Center for Theoretical Physics, Jagna, Bohol and Belmont Hotel, Mactan, Lapu-Lapu city, Cebu, Philippines in July and August 2023, Shanghai Jiao Tong University, China in June 2024, and Research Institute for Mathematical Sciences, Kyoto University, Japan in September 2024. He is grateful to those institutes, universities, schools, and hotels for their hospitality. He is supported by the World Premier International Research Center Initiative (WPI), MEXT, Japan, and JSPS KAKENHI Grant Number JP21K20317 and JP24K16888.

\section{Universality of Non-linear superalgebras}\label{sec:non-linear salg}

In this section, we recall results of De Sole and Kac in \cite{DK05, DK} about the universality of non-linear Lie superalgebras and non-linear Lie conformal superalgebras.
\smallskip

Let $\monoid$ be a discrete additive subgroup of $\R$. Set
\begin{align*}
\monoid_{\geq0} := \monoid \cap \R_{\geq0},\quad
\monoid_{>0} := \monoid \cap \R_{>0}.
\end{align*}
Denote by $\zeta_-$ the largest number in $\monoid_{\geq0}$ that is strictly smaller than $\zeta$ for each $\zeta \in \monoid_{>0}$.

\subsection{Pregraded vertex superalgebras}\label{sec:def-va}
Let $V$ be a vertex superalgebra. Denote by $\vac$ the vacuum vector, by $\partial$ the translation operator, by $p(a)$ the parity of $a \in V$ and by
\begin{align*}
a(z) = \sum_{n \in \Z}a_{(n)}z^{-n-1}
\end{align*}
the field on $V$ corresponding to $a \in V$. Let
\begin{align*}
[a_\lambda b] = \sum_{n=0}^\infty \frac{\lambda^n}{n!}a_{(n)}b \in \C[\lambda]\otimes V
\end{align*}
be the $\lambda$-bracket of $a$ and $b$ for $a, b \in V$, which satisfies the sesquilinearity, Jacobi identity, and the non-commutative Wick formulae. See, e.g., \cite{Kac}.
\smallskip

Let $\{e_\alpha\}_{\alpha \in I}$ be a subset of $V$ with an index set $I$.
\smallskip

Suppose that $I$ is ordered. Then $\mathbb{I}:= \Z_{\geq0}\times I$ is lexicographically ordered. That is, for $(k, \alpha), (l, \beta) \in \mathbb{I}$,
\begin{align}\label{eq:I_order}
(k, \alpha) < (l, \beta) \iff \alpha < \beta,\ \mathrm{or}\ \alpha=\beta\ \mathrm{and}\ k<l.
\end{align}
We shall equip $\mathbb{I}$ with this lexicographic order. For example, if $I$ is finite, we can fix such an order on $\mathbb{I}$. Set
\begin{align*}
e_i := \partial^{(k)}e_\alpha,\quad
i=(k, \alpha) \in \mathbb{I},
\end{align*}
where $\partial^{(k)} = \partial^k/k!$. We call that $V$ is \textit{freely generated by free  generators $\{e_\alpha\}_{\alpha \in I}$} if
\begin{align*}
\left\{ \NO{e_{i_1}e_{i_2}\cdots e_{i_s}}\ \mid s \geq0, i_j \in \mathbb{I}, i_j \leq i_{j+1}, \mathrm{and}\ i_{j}<i_{j+1}\ \mathrm{if}\ p(e_{i_j}) = \bar{1}\right\}
\end{align*}
forms a basis of $V$. Here each element $\NO{e_{i_1}e_{i_2}\cdots e_{i_s}}$ denotes
\begin{align*}
e_{i_1 (-1)}e_{i_2 (-1)}\cdots e_{i_s (-1)}\vac
= e_{\alpha_1 (-k_1-1)}e_{\alpha_2 (-k_2-1)}\cdots e_{\alpha_s (-k_s-1)}\vac
\end{align*}
if $i_j = (\alpha_j, k_j)$ for $j=1,2, \ldots, s$.
\smallskip

In addition, suppose that $\{e_\alpha\}_{\alpha \in I}$ is $\monoid_{>0}$-graded. That is, the following map is given:
\begin{align*}
\zeta_\bullet \colon I\rightarrow \monoid_{>0},\quad
\alpha \mapsto \zeta_\alpha,
\end{align*}
and the degree of $e_\alpha$ is defined by $\wt(e_\alpha) = \zeta_\alpha$. Set
\begin{align*}
\base = \Span_\C\{e_\alpha\}_{\alpha \in I}.
\end{align*}
Then $\base$ is $\monoid_{>0}$-graded:
\begin{align*}
\base = \bigoplus_{\zeta \in \monoid_{>0}}\base_\zeta,\quad
\base_\zeta = \Span_\C\{a \in \base \mid \wt(a) = \zeta\}.
\end{align*}
We extend the degree of $\base$ to
\begin{align*}
\C[\partial]\otimes\base = \sum_{k=0}^\infty\partial^k\base \subset V
\end{align*}
by $\wt(\partial a) = \wt(a)$. Define an increasing $\monoid_{\geq0}$-filtration of $V$ by
\begin{align}\label{eq:K-filt}
V_\zeta = \sum_{\begin{subarray}{c}k_1, \ldots, k_s \in \Z_{\geq0} \\ \zeta_1 + \cdots + \zeta_s \leq \zeta \end{subarray}}\NO{(\partial^{k_1}\base_{\zeta_1}) \cdots (\partial^{k_s}\base_{\zeta_s})},\quad
\zeta \in \monoid_{\geq0},
\end{align}
where the index $s$ in the summation of the RHS runs over all non-negative integers, and in the case of $s=0$, the summand is defined to be the one-dimensional subspace spanned by the vacuum vector $\vac$. In particular, $V_\zeta$ is non-zero for any $\zeta$.

\begin{definition}[\cite{DK05}]
Let $V$ be a vertex superalgebra freely generated by $\{e_\alpha\}_{\alpha \in I}$. $V$ is called pregraded by $\monoid_{\geq0}$ if there exists a $\monoid_{>0}$-grading on $\{e_\alpha\}_{\alpha \in I}$ such that
\begin{align*}
[{\base_{\zeta_1}}_\lambda\base_{\zeta_2}] \subset \C[\lambda]\otimes V_{(\zeta_1+\zeta_2)_-}
\end{align*}
for all $\zeta_1, \zeta_2 \in \monoid_{>0}$.
\end{definition}

\subsection{Non-linear Lie superalgebras}\label{sec:non-lin Liesalg}
Let $\base$ be a superspace and $T(\base)$ the tensor superalgebra of $\base$. Suppose that $\base$ is $\monoid_{>0}$-graded $\displaystyle \base = \bigoplus_{\zeta \in \monoid_{>0}}\base_\zeta$. Define an increasing $\monoid_{\geq0}$-filtration on $T(\base)$ by
\begin{align*}
T(\base)_\zeta = \sum_{\zeta_1+\cdots+\zeta_s \leq \zeta}\base_{\zeta_1}\otimes\cdots\otimes\base_{\zeta_s},\quad
\zeta \in \monoid_{\geq0}.
\end{align*}
A $\monoid_{>0}$-graded superspace $\base$ is called a non-linear Lie superalgebra if $\base$ is equipped with a parity-preserving linear map $[\ ,\ ] \colon \base\otimes\base \rightarrow T(\base)$ such that
\begin{enumerate}
\item $[\base_{\zeta_1}, \base_{\zeta_2}] \subset T(\base)_{(\zeta_1+\zeta_2)_-}$,
\item $[a, b] = -(-1)^{p(a)p(b)}[b, a]$ for all $a, b \in \base$,
\item $[a,[b, c]]-[[a, b], c]-(-1)^{p(a)p(b)}[b, [a,c]] \in M(\base)_{(\zeta_1+\zeta_2+\zeta_3)_-}$ for all $a\in\base_{\zeta_1}, b\in\base_{\zeta_2}, c\in\base_{\zeta_3}$, where
\begin{align*}
M(\base)_\zeta = \{& A \otimes (a\otimes b - (-1)^{p(a)p(b)}b\otimes a - [a, b])\otimes B\ |\ a\in \base_{\zeta_1}, b \in \base_{\zeta_2},\\
&A \in T(\base)_{\zeta_3}, B \in T(\base)_{\zeta_4}, \zeta_1+\zeta_2+\zeta_3+\zeta_4 \leq \zeta\}.
\end{align*}
\end{enumerate}
Then $M(\base) = \sum_{\zeta \in \monoid_{\geq0}}M(\base)_\zeta$ is a two-sided ideal in $T(\base)$. Define an associative superalgebra
\begin{align*}
U(\base) = T(\base)/M(\base),
\end{align*}
called the universal enveloping algebra of $\base$.

\begin{definition}
Let $U$ be an associative $\C$-superalgebra and $\{e_\alpha\}_{\alpha \in I}$ be a set of parity-homogeneous generators of $U$.\\
(1) $U$ is called \textit{PBW generated by $\{e_\alpha\}_{\alpha \in I}$} if $I$ is an ordered set and
\begin{align*}
\{ e_{\alpha_1}\cdots e_{\alpha_s} \mid s \geq0, \alpha_j \in I, \alpha_j \leq \alpha_{j+1}, \mathrm{and}\ \alpha_{j}<\alpha_{j+1}\ \mathrm{if}\ p(e_{\alpha_j}) = \bar{1}\}
\end{align*}
forms a basis of $U$. Then we call $\{e_\alpha\}_{\alpha \in I}$ PBW generators of $U$.\\
(2) Suppose that $\{e_\alpha\}_{\alpha \in I}$ is $\monoid_{>0}$-graded by $\wt(e_\alpha) = \zeta_\alpha \in \monoid_{>0}$. Let
\begin{align*}
\base = \bigoplus_{\zeta \in \monoid_{>0}}\base_\zeta,\quad
\base_\zeta = \Span_\C\{e_\alpha \mid \alpha \in I, \wt(e_\alpha) = \zeta\} \subset U.
\end{align*}
Define an increasing $\monoid_{\geq0}$-filtration on $U$ by
\begin{align*}
U_\zeta = \sum_{\zeta_1+\cdots+\zeta_s \leq \zeta}\base_{\zeta_1}\cdots\base_{\zeta_s}.
\end{align*}
Then $U$ is called \textit{pregraded by $\monoid_{\geq0}$} if
\begin{align*}
[\base_{\zeta_1}, \base_{\zeta_2}] \subset U_{(\zeta_1+\zeta_2)_-},
\end{align*}
where $[a, b] = ab - (-1)^{p(a)p(b)}b a$.
\end{definition}

\begin{proposition}[Universality of non-linear Lie superalgebras by {\cite[Theorem 3.6]{DK}}]\label{prop:DK-U(r)}
Let $U$ be an associative superalgebra and $\base$ be a $\monoid_{>0}$-graded subspace of $U$ with an ordered homogenous basis $\{e_\alpha\}_{\alpha \in I}$. Suppose that $\{e_\alpha\}_{\alpha \in I}$ are PBW generators of $U$ and $U$ is pregraded by $\monoid_{\geq0}$. Let $\rho \colon U \rightarrow T(\base)$ be a linear map defined by
\begin{align*}
e_{\alpha_1}\cdots e_{\alpha_s} \mapsto e_{\alpha_1}\otimes\cdots \otimes e_{\alpha_s}
\end{align*}
for ordered monomials $e_{\alpha_1}\cdots e_{\alpha_s}$ in $U$. Then $\base$ is a non-linear Lie superalgebra by
\begin{align*}
[a, b] = \rho(ab-(-1)^{p(a)p(b)}ba),\quad
a, b \in \base
\end{align*}
such that $U \simeq U(\base)$.
\end{proposition}

\subsection{Non-linear Lie conformal superalgebras}\label{sec:non-linear lie conf salg}

Let $\conf$ be a $\C[\partial]$-module equipped with a $\monoid_{>0}$-grading
\begin{align*}
\conf = \bigoplus_{\zeta \in \monoid_{>0}}\conf_\zeta
\end{align*}
such that each $\conf_\zeta$ is a $\C[\partial]$-submodule. Then, the tensor superalgebra $T(\conf)$ is naturally a $\C[\partial]$-module by extending the action of $\partial$ as an even derivation. Define an increasing $\monoid_{\geq0}$-filtration $T(\conf)_\zeta$ as in Section \ref{sec:non-lin Liesalg}. Suppose that $\conf$ is equipped with a parity-homogeneous linear map $[{\ }_\lambda\ ] \colon \conf\otimes \conf \rightarrow \C[\lambda] \otimes T(\conf)$ called the $\lambda$-bracket such that
\begin{enumerate}
\item $[{\conf_{\zeta_1}}_\lambda \conf_{\zeta_2}] \subset \C[\lambda] \otimes T(\conf)_{(\zeta_1+\zeta_2)_-}$,
\item $[\partial a_\lambda b] = -\lambda[a_\lambda b]$ and $[a_\lambda \partial b] = (\lambda + \partial)[a_\lambda b]$.
\item $[a_\lambda b] = -(-1)^{p(a)p(b)}[b_{-\lambda-\partial}a]$,
\end{enumerate}
\begin{lemma}[{\cite[Lemma 3.2]{DK05}}]The $\lambda$-bracket can be uniquely extended to $T(\conf)$ as a parity-homogeneous linear map
\begin{align*}
T(\conf) \otimes T(\conf) \ni a \otimes b \mapsto [a_\lambda b] \in \C[\lambda] \otimes T(\conf),
\end{align*}
and the normally ordered product on $T(\conf)$ is uniquely defined as a parity-homogeneous linear map
\begin{align*}
T(\conf) \otimes T(\conf) \ni a \otimes b \mapsto\ \NO{ab}\ \in T(\conf)
\end{align*}
by the following equations for $a, b \in \conf$ and $A, B \in T(\conf)$:
\begin{align}\label{eq:lambda-NO-def}
\begin{split}&
[1_\lambda A] = 0 = [A_\lambda 1],\quad
\NO{1 A} = A = \NO{A 1},\quad
\NO{a \otimes A} = a \otimes A,
\\&
\NO{(\NO{a A})B}\ =\ \NO{a (\NO{AB})} + \NO{\left(\int_0^\partial a\ d\lambda\right)[A_\lambda B]}
\\&\qquad
+ (-1)^{p(a)p(A)}\NO{\left(\int_0^\partial A d\lambda\right)[a_\lambda B]},
\\&
[a_\lambda\NO{b A}]\ =\ \NO{[a_\lambda b]A} + (-1)^{p(a)p(b)}\NO{b[a_\lambda A]}+\int_0^\lambda[[a_\lambda b]_\mu A]d\mu,
\\&
[{\NO{a A}}_\lambda B]\ =\ \NO{(\mathrm{e}^{\partial\frac{d}{d\lambda}}a)[A_\lambda B]}+(-1)^{p(a)p(A)}\NO{(\mathrm{e}^{\partial\frac{d}{d\lambda}}A)[a_\lambda B]}
\\&\qquad
+ (-1)^{p(a)p(A)} \int_0^\lambda [A_\mu[a_{\lambda-\mu}B]]d\mu.
\end{split}
\end{align}
Then we have
\begin{align*}
\NO{T(\conf)_{\zeta_1}T(\conf)_{\zeta_2}}\ \subset T(\conf)_{\zeta_1+\zeta_2},\quad
[{T(\conf)_{\zeta_1}}_\lambda T(\conf)_{\zeta_2}] \subset T(\conf)_{(\zeta_1+\zeta_2)_-}.
\end{align*}
\end{lemma}

A $\monoid_{>0}$-graded $\C[\partial]$-module $\conf$ equipped with the $\lambda$-bracket above is called a non-linear Lie conformal superalgebra if
\begin{align*}
[a_\lambda [b_\mu c]]-(-1)^{p(a)p(b)}[b_\mu[a_\lambda c]]-[[a_\lambda b]_{\lambda+\mu}c] \in \C[\lambda, \mu] \otimes M(\conf)_{(\zeta_1+\zeta_2+\zeta_3)_-}
\end{align*}
for all $a \in \conf_{\zeta_1}$, $b \in \conf_{\zeta_2}$, $c \in \conf_{\zeta_3}$, where
\begin{align}\label{eq:M(conf)_zeta}
\begin{split}
&M(\conf)_\zeta = \Big\{A \otimes (a \otimes b \otimes B - (-1)^{p(a)p(b)}b\otimes a \otimes B - \NO{\left(\int_{-\partial}^0 [a_\lambda b]d\lambda\right)B}\ )\ \Big|\\
&\qquad a \in \conf_{\zeta_1}, b \in \conf_{\zeta_2} \in \conf, A \in T(\conf)_{\zeta_3}, B \in T(\conf)_{\zeta_4}, \zeta_1+\zeta_2+\zeta_3+\zeta_4 \leq \zeta\Big\}.
\end{split}
\end{align}
Let $M(\conf) = \sum_{\zeta \in \monoid_{\geq0}}M(\conf)_{\zeta}$. By \cite[Theorem 3.9]{DK05},
\begin{align*}
V(\conf) = T(\conf)/M(\conf)
\end{align*}
has a unique vertex superalgebra structure such that
\begin{enumerate}
\item the vacuum vector is $1$,
\item the translation operator is $\partial$,
\item the normally ordered operator and $\lambda$-bracket are induced from $T(\conf)$.
\end{enumerate}
Then, $V(\conf)$ is called the universal enveloping vertex superalgebra of $\conf$.

\begin{proposition}[Universality of non-linear Lie conformal superalgebras by {\cite[Theorem 6.4]{DK05}}]\label{prop:DK-V(R)}
Let $V$ be a vertex superalgebra and $\conf$ be a $\monoid_{>0}$-graded free $\C[\partial]$-submodule of $V$ with an ordered homogeneous basis $\{e_i\}_{i\in\mathbb{I}}$. Suppose that $\{e_\alpha\}_{\alpha \in I}$ are free generators of $V$ and $V$ is pregraded by $\monoid_{>0}$. Let $\rho \colon V \hookrightarrow T(\conf)$ be a linear embedding by
\begin{align*}
\NO{e_{i_1} \cdots e_{i_s}}\ \mapsto e_{i_1} \otimes \cdots \otimes e_{i_s}
\end{align*}
for ordered monomials $\NO{e_{i_1} \cdots e_{i_s}}$ in $V$. Then $\conf$ has a structure of a non-linear Lie conformal superalgebra by
\begin{align*}
\conf \otimes \conf \ni a \otimes b \mapsto \rho([a_\lambda b]) \in T(\conf)
\end{align*}
such that $V \simeq V(\conf)$.
\end{proposition}

\section{Twisted Zhu algebras}\label{sec:twisted Zhu}

\subsection{Twisted modules}\label{sec:twisted modules}

A Hamiltonian operator $H$ on $V$ is a semisimple operator on $V$ satisfying that
\begin{align*}
[H, Y(a, z)] = z\partial_zY(a, z) + Y(H(a), z)
\end{align*}
for all $a \in V$. The eigenvalue of $H$ is called the conformal weight. If $V$ is conformal and $\displaystyle L(z) = \sum_{n \in \Z}L_nz^{-n-2}$ is the field corresponding to the conformal vector of $V$, we may choose $H = L_0$ as a Hamiltonian operator. We will use the formal expansions
\begin{align*}
&\left.\frac{1}{(z-w)^{n+1}}\right|_{|z|>|w|} = \frac{(-\partial_z)^{n}}{n!}\sum_{m=0}^\infty z^{-m-1}w^m,\\
&\left.\frac{1}{(z-w)^{n+1}}\right|_{|z|<|w|} = -\frac{(-\partial_z)^{n}}{n!}\sum_{m=0}^\infty z^{m}w^{-m-1}
\end{align*}
for $n \in \Z_{\geq0}$ and the formal delta function
\begin{align*}
\delta(z-w) = \sum_{m\in\Z}z^{-m-1}w^{m} = \left.\frac{1}{z-w}\right|_{|z|>|w|} - \left.\frac{1}{z-w}\right|_{|z|<|w|}.
\end{align*}

Let $V$ be a vertex superalgebra equipped with a Hamiltonian operator $H$ whose eigenvalues belong to $\R$, $\operatorname{Aut}_H V$ be the group of automorphisms on $V$ commuting with $H$, and $g \in \operatorname{Aut}_H V$ be a diagonalizable automorphism with modulus $1$ eigenvalues. Then $V$ is simultaneously diagonalized by $g$ and $H$:
\begin{align*}
V = \bigoplus_{\begin{subarray}{c} \bar{\gamma} \in \R/\Z, \\ \Delta \in \R \end{subarray}}V(\bar{\gamma}, \Delta),\quad
V(\bar{\gamma}, \Delta) = \{ a \in V \mid g(a) = \mathrm{e}^{2\pi i \gamma}a,\ H(a) = \Delta a\},
\end{align*}
where $\gamma$ is a representative of $\bar{\gamma}$ in $\R$. If $a \in V(\bar{\gamma}, \Delta)$, denote by $\Delta_a := \Delta$, and set
\begin{align*}
\gamma_a := \epsilon_a + \Delta_a,
\end{align*}
where $\epsilon_a$ is a unique real number such that $-1 < \epsilon_a \leq 0$ and $\epsilon_a + \Delta_a \in \bar{\gamma}$. Then $\gamma_a$ is a representative of $\bar{\gamma}$. Since
\begin{equation}\label{eq:Delta+epsilon}
\begin{split}
&\Delta_{\vac} = 0,\quad
\Delta_{\partial a} = \Delta_a + 1,\quad
\Delta_{a_{(n)}b} = \Delta_a + \Delta_b - n - 1,\\
&\epsilon_{\vac} = 0,\quad
\epsilon_{\partial a} = \epsilon_a,\quad
\epsilon_{a_{(n)}b} = \epsilon_a + \epsilon_b + \lfloor -\epsilon_a-\epsilon_b \rfloor,
\end{split}
\end{equation}
we have
\begin{align}\label{eq:gamma}
\gamma_{\vac} = 0,\quad
\gamma_{\partial a} = \gamma_a + 1,\quad
\gamma_{a_{(n)}b} = \gamma_a + \gamma_b - n - 1 + \lfloor -\epsilon_a-\epsilon_b \rfloor,
\end{align}
where $\lfloor x \rfloor = \max\{ n \in \Z \mid n \leq x\}$ is the floor function. Set
\begin{align}\label{eq:def-chi}
\chi_{a, b} := \lfloor -\epsilon_a-\epsilon_b \rfloor.
\end{align}
Using the fact that $-1 < \epsilon_a, \epsilon_b \leq 0$, it follows that
\begin{align*}
\chi_{a, b} = \lfloor -\epsilon_a-\epsilon_b \rfloor =
\begin{cases}
1 & \mathrm{if}\ \epsilon_a + \epsilon_b \leq -1,\\
0 & \mathrm{if}\ \epsilon_a + \epsilon_b > -1.
\end{cases}
\end{align*}
A $g$-twisted $V$-module $M$ is a superspace equipped with a parity-preserving linear map
\begin{align*}
Y_M \colon V \ni a \rightarrow Y_M(a, z) = a^M(z) = \sum_{n \in \Z + \gamma_a}a^M_{(n)} z^{-n-1} \in \sum_{n \in \R}(\End M)z^n
\end{align*}
such that $Y_M(\vac, z) = \id_M$, $a^M_{(n)}m = 0$ for $n \gg 0$ and $Y_M$ satisfies the $g$-twisted Borcherds identity: for $a, b \in V$, $\ell \in \Z + \gamma_a$, $m \in \Z + \gamma_b$, $n \in \Z$,
\begin{equation}\label{eq:Borcherds}
\begin{split}
&\sum_{j=0}^\infty(-1)^j\binom{n}{j}\left(
a^M_{(\ell+n-j)}b^M_{(m+j)} - (-1)^{p(a)p(b)+n}b^M_{(m+n-j)}a^M_{(\ell+j)}
\right)\\
&= \sum_{j=0}^\infty \binom{\ell}{j}\left(a_{(n+j)}b
\right)^M_{(\ell+m-j)}.
\end{split}
\end{equation}
Set $a^M_n = a^M_{(n+\Delta_a-1)}$ for $a \in V$ so that
\begin{align*}
a^M(z) = \sum_{n \in \Z + \epsilon_a} a^M_n z^{-n-\Delta_a}.
\end{align*}
Then the $g$-twisted Borcherds identity \eqref{eq:Borcherds} is equivalent to
\begin{equation}\label{eq:Borcherds-shift}
\begin{split}
&\sum_{j=0}^\infty(-1)^j\binom{n}{j}\left(
a^M_{\ell+n-j}b^M_{m+j-n} - (-1)^{p(a)p(b)+n}b^M_{m-j}a^M_{\ell+j}
\right)\\
&= \sum_{j=0}^\infty \binom{\ell+\Delta_a-1}{j}\left(a_{(n+j)}b
\right)^M_{\ell+m}
\end{split}
\end{equation}
for $a, b \in V$, $\ell \in \Z+ \epsilon_a$, $m \in \Z+ \epsilon_b$, $n \in \Z$.

\subsection{Zhu's Theorem}\label{sec:Gamma/Z}

For $a, b \in V$ and $n \in \Z$, the $*_n$-product $a*_n b$ is defined by
\begin{align*}
(1+z)^{\gamma_a}a(z)b = \sum_{n \in \Z}(a*_n b)z^{-n-1},
\end{align*}
where $\displaystyle (1+z)^{\gamma_a} = \sum_{j=0}^\infty \binom{\gamma_a}{j}z^j$. Hence
\begin{align*}
a*_n b = \sum_{j=0}^\infty\binom{\gamma_a}{j}a_{(j+n)}b.
\end{align*}
For $a, b \in V$, the $*$-bracket $[a_*b]$ is defined by
\begin{align}\label{eq:*-br}
[a_* b] = \underset{z=0}{\operatorname{Res}}(1+z)^{\gamma_a-1}a(z)b = \sum_{j =0}^\infty(-1)^j a *_{j}b = \sum_{j =0}^\infty\binom{\gamma_a-1}{j} a_{(j)}b.
\end{align}
For $a, b \in V$, the $*$-product $a * b$ and $\circ$-product $a \circ b$ are defined by
\begin{align*}
a * b = a *_{-1} b,\quad
a \circ b = a *_{-2} b.
\end{align*}
Define a subspace $\gH{V}$ of $V$ and a subspace $\gH{J}(V)$ of $\gH{V}$ by
\begin{align}\label{eq:def-VandJ}
\begin{array}{l}
\gH{V} = \Span_\C\{a \in V \mid \epsilon_a = 0\},\\
\gH{J}(V) = \Span_\C\{a * _{-2 + \chi_{a, b}}b \in V \mid \epsilon_a + \epsilon_b \in \Z\}.
\end{array}
\end{align}
We often write $\gH{J}$ for $\gH{J}(V)$ for simplicity. Since
\begin{align*}
\gH{V} = V^{g\circ \twist^{-1}} = \{a \in V \mid (g\circ \twist^{-1})(a) = a\},\quad
\twist := \mathrm{e}^{2\pi i H}\in \operatorname{Aut}_H V,
\end{align*}
$\gH{V}$ is a vertex subalgebra of $V$.
\begin{proposition}[{\cite[Theorem 2.13]{DK}}]\label{prop:J-ideal}\quad
\begin{enumerate}
\item $\vac$ is the identity element for the $*$-product.
\item $\gH{J}$ contains all elements $a *_{-j+\chi_{a, b}}b$ for $j \geq 2$ and $a, b \in V$ such that $\epsilon_a + \epsilon_b \in \Z$.
\item $(\partial + H) \gH{V} \subset \gH{J}$.
\item $\gH{J}$ is a two-sided ideal of $\gH{V}$ with respect to the $*$-product and $*$-bracket.
\item $a * b -(-1)^{p(a)p(b)}b*a \equiv (1-\chi_{a, b})[a_* b]$ mod $\gH{J}$ for $a, b \in V$ such that $\epsilon_a + \epsilon_b \in \Z$.
\item The $*$-product is associative in $\gH{V}$ mod $\gH{J}$.
\end{enumerate}
\end{proposition}
See Section \ref{sec:properties_*} for the proof of Proposition \ref{prop:J-ideal}.
\smallskip

By Proposition \ref{prop:J-ideal},
\begin{align*}
\gH{\Zhu}V = \gH{V}/\gH{J}
\end{align*}
has a structure of an associative $\C$-superalgebra by the induced $*$-product. $\gH{\Zhu}V$ is called the $(g, H)$-twisted Zhu algebra of $V$, introduced by Zhu \cite{Zhu} for $g=1$ on a $\Z$-graded (non-super) $V$, by Kac--Wang \cite{KWan} for $g=1$ on a $\frac{1}{2}\Z$-graded $V$ such that $V_{\bar{0}} = \bigoplus_{\Delta \in \Z}V(\Delta)$ and $V_{\bar{1}} = \bigoplus_{\Delta \in \Z + \frac{1}{2}}V(\Delta)$ (correct statistic), by Dong--Li--Mason \cite{DLM} for $g$ of finite order on a $\Z$-graded (non-super) $V$, by Dong--Zhao for $g$ of finite order on a correct statistic $V$, and by De Sole--Kac \cite{DK} in general cases.

A $g$-twisted $V$-module $M$ is called lower-bounded (or positive-energy) if $M$ has an $\R_{\geq 0}$-grading
\begin{align*}
M = \bigoplus_{j \geq 0}M_j,\quad
M_0 \neq 0
\end{align*}
such that
\begin{align*}
a^M_n M_j \subset M_{j-n},\quad
a \in V,\quad
n \in \Z + \epsilon_a,\quad
j \in \R,
\end{align*}
where $M_j = 0$ for $j < 0$. Then $M_0$ is called the top space, denoted by $M_\mathrm{top}$. By \cite[Lemma 2.22]{DK}, a linear map
\begin{align*}
\gH{V} \ni a \mapsto a^M|_{M_\mathrm{top}} \in \End M_\mathrm{top}
\end{align*}
induces a homomorphism $\gH{\Zhu}V \rightarrow \End M_\mathrm{top}$. Thus, we have a functor
\begin{align*}
\psi \colon M \mapsto M_\mathrm{top}
\end{align*}
from the category of lower-bounded $g$-twisted $V$-modules to the category of $\gH{\Zhu}V$-modules. The following theorem is essentially proved by Zhu (see, e.g., a proof of \cite[Theorem 2.30]{DK}):

\begin{theorem}[{\cite{Zhu}}]
The functor $\psi$ establishes a bijection (up to isomorphisms) between simple lower-bounded $g$-twisted $V$-modules and simple $\gH{\Zhu}V$-modules.
\end{theorem}

\subsection{Properties of products}\label{sec:properties_*}

We collect properties of the $*_n$-products and $*$-bracket given in \cite{DK}, which are proved only in the case of $g = \mathrm{e}^{2\pi i H}$, but essentially the same proofs apply. We also give a proof of Proposition \ref{prop:J-ideal}.

Let $H_g \colon V \rightarrow V$ be a linear map defined by $H_g(a) = \gamma_a a$, and
\begin{align*}
Y_*(a, z) = (1+z)^{\gamma_a}Y(a, z),\quad
a \in V.
\end{align*}

\begin{proposition}[{\cite[(2.7)(2.8)]{DK}}]\label{prop:trans-1}Let $a, b \in V$ and $n \in \Z$.
\begin{enumerate}
\item $\left(\partial a + (n+1+\gamma_a)a\right)*_n b
= -n a*_{n-1}b$.
\item $\displaystyle \partial(a*_{n}b)
= \sum_{j=0}^\infty (-1)^j(\partial a)*_{n+j}b + a*_{n}\partial b$.
\end{enumerate}
\begin{proof}
Since
$
(1+z)^{\gamma_a +1}Y(\partial a, z)
= (1+z)^{\gamma_a +1}[\partial, Y(a, z)]
= (1+z)^{\gamma_a +1}\partial_z Y(a, z),
$
we have
$
Y_{*}(\partial a, z)
= (1+z)[\partial, Y_*(a, z)]
= \left((1+z)\partial_z - \gamma_a\right)Y_*(a, z),
[\partial, Y_*(a, z)]
= (1 + z)^{-1}Y_*(\partial a, z)
= \sum_{j=0}^\infty (-z)^j Y_*(\partial a, z).
$
These prove (1) and (2), respectively.
\end{proof}
\end{proposition}
\begin{corollary}[{\cite[(2.12)]{DK}}]\label{cor:a(-k-1)b-ind}
For $k \in \Z_{\geq0}$,
\begin{align*}
a*_{-k-1}b = \sum_{j =0}^k \binom{\gamma_a}{k-j}(\partial^{(j)}a)*b.
\end{align*}
\begin{proof}
The case $k=0$ is clear. For $k \geq 1$, by Proposition \ref{prop:trans-1}(1),
\begin{align}\label{eq:*_k-ind}
a*_{-k-1}b = \frac{1}{k}(\partial a)*_{-k}b + \frac{\gamma_a -k +1}{k} a*_{-k}b.
\end{align}
Then, the assertion follows inductively.
\end{proof}
\end{corollary}
\begin{proposition}[{\cite[(2.9)(2.10)(2.11)]{DK}}]\label{prop:trans-id}
Let $a, b \in V$.
\begin{enumerate}
\item $a\circ b = \left((\partial + H_g)a\right)*b$.
\item $ [(\partial + H_g)a_*b] = 0$.
\item $\left(\partial +  H_g- \chi_{a, b} \right)[a_* b] = [a_* (\partial +  H_g)b]$.
\end{enumerate}
\begin{proof}
(1) follows from Proposition \ref{prop:trans-1}(1) for $n=-1$. (2) is easy to check by definition. To show (3), we have
\begin{align*}
&\left(\partial +  H_g\right)[a_*b] - [a_* (\partial + H_g)b]\\
=& \sum_{j=0}^\infty \binom{\gamma_a-1}{j}\left\{ \partial(a_{(j)}b) + (\gamma_a+\gamma_b+\chi_{a, b}-j-1)a_{(j)}b - a_{(j)}\partial b - \gamma_b a_{(j)} b\right\}\\
=& \sum_{j=0}^\infty \binom{\gamma_a-1}{j} \left( (\gamma_a+\chi_{a, b})-j-1)a_{(j)}b + (\partial a)_{(j)}b\right)\\
=& \sum_{j=0}^\infty \left( (\gamma_a+\chi_{a, b}-j-1)\binom{\gamma_a-1}{j} - (j+1)\binom{\gamma_a-1}{j+1}\right)a_{(j)}b\\
=& \chi_{a, b}\sum_{j=0}^\infty\binom{\gamma_a -1}{j}  a_{(j)}b
= \chi_{a, b}[a_* b].
\end{align*}
This proves (3).
\end{proof}
\end{proposition}
\begin{proposition}[{\cite[Proposition 2.2]{DK}}]\label{prop:n-th prod}For $a, b \in V$ and $n \in \Z$,
\begin{align*}
(1+w)^{n+1-\chi_{a, b}}Y_*(a*_{n}b, w) = Y_*(a, w)_{(n)}Y_*(b, w).
\end{align*}
\begin{proof}
First of all,
$%\displaystyle
(1+ w)^{n+1-\chi_{a, b}}Y_*(a*_{n}b, w)
= \sum_{j=0}^\infty \binom{\gamma_a}{j} (1+ w)^{n+1-\chi_{a, b}}Y_*(a_{(n+j)}b, w)
= \sum_{j=0}^\infty \binom{\gamma_a}{j} (1+ w)^{\gamma_a + \gamma_b -j}Y(a_{(n+j)}b, w).
$
Since
\begin{align*}
&Y(a_{(n+j)}b, w)\\
=& \left.\underset{z=0}{\operatorname{Res}}\ (z-w)^{n+j}Y(a, z)Y(b, w)\right|_{|z|>|w|} - (-1)^{p(a)p(b)} \left.\underset{z=0}{\operatorname{Res}}\ (z-w)^{n+j}Y(b, w)Y(a, z)\right|_{|z|<|w|},
\end{align*}
we have
\begin{align*}
&(1+ w)^{n+1-\chi_{a, b}}Y_*(a*_{n}b, w)\\
=&\underset{z=0}{\operatorname{Res}}\ \sum_{j=0}^\infty\binom{\gamma_a}{j}(z-w)^j(1+ w)^{\gamma_a + \gamma_b -j}\\
\cdot &\left\{\left.(z-w)^{n}Y(a, z)Y(b, w)\right|_{|z|>|w|} - (-1)^{p(a)p(b)} \left.(z-w)^{n}Y(b, w)Y(a, z)\right|_{|z|<|w|}\right\}\\
=&\underset{z=0}{\operatorname{Res}}\ (1+ z)^{\gamma_a}(1+ w)^{\gamma_b}\\
\cdot &\left\{\left.(z-w)^{n}Y(a, z)Y(b, w)\right|_{|z|>|w|} - (-1)^{p(a)p(b)} \left.(z-w)^{n}Y(b, w)Y(a, z)\right|_{|z|<|w|}\right\}\\
=&\left.\underset{z=0}{\operatorname{Res}}\ (z-w)^{n}Y_*(a, z)Y_*(b, w)\right|_{|z|>|w|} - (-1)^{p(a)p(b)}\left.\underset{z=0}{\operatorname{Res}}\ (z-w)^{n}Y_*(b, w)Y_*(a, z)\right|_{|z|<|w|}\\
=&Y_*(a, w)_{(n)}Y_*(b, w).
\end{align*}
This completes the proof.
\end{proof}
\end{proposition}
\begin{corollary}[{\cite[(2.14)]{DK}}]Let $a, b \in V$.
\begin{enumerate}
\item $\displaystyle [Y_*(a, z), Y_*(b, w)] = \sum_{j=0}^\infty(1+ w)^{j+1-\chi_{a, b}}Y_*(a*_{j}b, w)\partial^{(j)}_w\delta(z-w)$.
\item $Y_*(a*b, w) = (1+ w)^{\chi_{a, b}}\NO{Y_*(a, w)Y_*(b, w)}$\ .
\end{enumerate}
\begin{proof}
Apply Proposition \ref{prop:n-th prod} for $n=0$ and $n=-1$.
\end{proof}
\end{corollary}
\begin{theorem}[{\cite[(2.13)]{DK}}]\label{thm:Borcherds-id}
Let $a, b, c \in V$ and $m, n \in \Z$. Then
\begin{align*}
&\sum_{j=0}^\infty (-1)^j \binom{n}{j} (a*_{m+n-j}Y_*(b, w)c\ w^j - (-1)^{p(a)p(b)+n}Y_*(b, w)a*_{m+j}c\ w^{n-j})\\
&= \sum_{j=0}^\infty(1+ w)^{n+j+1-\chi_{a, b}}\binom{m}{j}Y_*(a*_{n+j}b, w)c\ w^{m-j}.
\end{align*}
\begin{proof}
Using Proposition \ref{prop:n-th prod}, it follows that
\begin{align*}
&\left.(z-w)^{n}Y_*(a, z)Y_*(b, w)\right|_{|z|>|w|} - (-1)^{p(a)p(b)}\left.(z-w)^{n}Y_*(b, w)Y_*(a, z)\right|_{|z|<|w|}\\
=& \sum_{j=0}^\infty \underset{z=0}{\operatorname{Res}}
\Bigl\{\left.(z-w)^{n+j}Y_*(a, z)Y_*(b, w)\right|_{|z|>|w|}\\
&\quad- (-1)^{p(a)p(b)}\left.(z-w)^{n+j}Y_*(b, w)Y_*(a, z)\right|_{|z|<|w|}\Bigr\}\partial^{(j)}_w\delta(z-w)\\
=& (1+ w)^{n+j+1-\chi_{a, b}}Y_*(a*_{n+j}b, w)\partial^{(j)}_w\delta(z-w).
\end{align*}
Now we compare the coefficients of $z^{-m-1}$ of both sides above applied to $c \in V$ and then get the assertion. 
\end{proof}
\end{theorem}
\begin{corollary}[{\cite[(2.15)]{DK}}]\label{cor:Borcherds-id}
Let $a, b, c \in V$ and $m, n, k \in \Z$. Then
\begin{multline*}
\sum_{i, j=0}^\infty \binom{n}{j} \binom{-n-1+\chi_{a, b}}{i}
\left(
a*_{m+n-j} b*_{k+j+i} c
-(-1)^{p(a)p(b)+n} b*_{k+n-j+i} a*_{m+j}c
\right)\\
=\sum_{i, j=0}^\infty \binom{m}{j} \binom{j}{i} 
\left(
a*_{n+j}b
\right)*_{k+m-j+i}c.
\end{multline*}
\begin{proof}
Multiply both sides of Theorem \ref{thm:Borcherds-id} by $(1+ w)^{-n-1+\chi_{a, b}}$ and then compare the coefficients of $w^{-k-1}$.
\end{proof}
\end{corollary}
\begin{proposition}[{\cite[(2.16)]{DK}}]\label{prop:Borcherds-id}
For $a, b, c \in V$ and $n, k \in\Z$,
\begin{multline*}
(a*_{n}b)*_{k}c\\
= \sum_{i, j=0}^\infty(-1)^j\binom{n}{j}\binom{-n-1+\chi_{a, b}}{i}\left(
a*_{n-j}b*_{k+j+i}c-(-1)^{p(a)p(b)+n}b*_{k+n-j+i}a*_{j}c
\right).
\end{multline*}
\begin{proof}
Apply Corollary \ref{cor:Borcherds-id} for $m=0$.
\end{proof}
\end{proposition}
\begin{corollary}[{\cite[(2.17)]{DK}}]\label{cor:quasi-asso}
For $a, b, c \in V$,
\begin{multline*}
(a*b)*c - a*(b*c)\\
=\sum_{j=0}^\infty
\left\{
(a*_{-j-2}+\chi_{a, b} a*_{-j-1})b*_{j}c
+(-1)^{p(a)p(b)}(b*_{-j-2}+\chi_{a, b} b*_{-j-1})a*_{j}c
\right\}.
\end{multline*}
\begin{proof}
Apply Proposition \ref{prop:Borcherds-id} for $n=k=-1$.
\end{proof}
\end{corollary}
\begin{proposition}[{\cite[(2.18)]{DK}}]\label{prop:AcircB*C}
For $a, b, c\in V$,
\begin{align*}
&(a\circ b)*c = a\circ (b*c)\\
&+\sum_{j=0}^\infty
\Bigl\{
\left(
(\partial a)*_{-j-2}+\chi_{a, b} (\partial a)*_{-j-1}
\right)b*_{j}c
+(-1)^{p(a)p(b)}
\left(
b*_{-j-2}+\chi_{a, b} b*_{-j-1}
\right)
(\partial a)*_{j}c
\Bigr\}\\
&+\gamma_a
\sum_{j=0}^\infty
\Bigl\{
\left(
a*_{-j-2}+\chi_{a, b} a*_{-j-1}
\right)b*_{j}c
+(-1)^{p(a)p(b)}\gamma_a
\left(
b*_{-j-2}+\chi_{a, b} b*_{-j-1}
\right)
a*_{j}c
\Bigr\}.
\end{align*}
\begin{proof}
By Proposition \ref{prop:trans-id} (1) and Corollary \ref{cor:quasi-asso},
\begin{align*}
&(a\circ b)*c-a\circ (b*c)
=\left(
(\partial a+ \gamma_a a)*b
\right)*c
-(\partial a + \gamma_a a)*(b*c)\\
=&\sum_{j=0}^\infty
\Bigl\{
\left(
(\partial a+ \gamma_a a)*_{-j-2}+\chi_{a, b} (\partial a+ \gamma_a a)*_{-j-1}
\right)b*_{j}c\\
&\hspace{15mm}+(-1)^{p(a)p(b)}(b*_{-j-2}+\chi_{a, b} b*_{-j-1})(\partial a+ \gamma_a a)*_{j}c
\Bigr\}\\
=&\sum_{j=0}^\infty
\Bigl\{
\left(
(\partial a)*_{-j-2}+\chi_{a, b} (\partial a)*_{-j-1}
\right)b*_{j}c
+(-1)^{p(a)p(b)}(b*_{-j-2}+\chi_{a, b} b*_{-j-1})(\partial a)*_{j}c\\
&\quad+\gamma_a\left(
a*_{-j-2}+\chi_{a, b} a*_{-j-1}
\right)b*_{j}c
+(-1)^{p(a)p(b)}\gamma_a(b*_{-j-2}+\chi_{a, b} b*_{-j-1})a*_{j}c
\Bigr\}.
\end{align*}
Here we use $\chi_{\partial a, b} = \chi_{a, b}$.
\end{proof}
\end{proposition}
\begin{proposition}[{\cite[(2.19)]{DK}}]\label{prop:A*Bcomm}
For $a, b \in V$,
\begin{align*}
&a*b-(-1)^{p(a)p(b)}b*a-(1-\chi_{a, b})[a_*b]\\
&=\sum_{n=1}^\infty\sum_{l=1}^n\binom{\gamma_b}{n-l}(-1)^{n+1}
\left\{
\partial^{(l)}(a_{(n-1)}b)-\binom{-\gamma_{a_{(n-1)}b}}{l}a_{(n-1)}b
\right\}.
\end{align*}
\begin{proof}
Recall the skew-symmetry
$%\displaystyle
(-1)^{p(a)p(b)}b_{(n)}a = \sum_{i=0}^\infty(-1)^{n+i+1}\partial^{(i)}(a_{(n+i)}b).
$
Then we have
\begin{align*}&
-(-1)^{p(a)p(b)}b*a
=-(-1)^{p(a)p(b)}\sum_{j=0}^\infty\binom{\gamma_b}{j} b_{(j-1)}a
=\sum_{i, j=0}^\infty\binom{\gamma_b}{j}(-1)^{i+j+1}\partial^{(i)}(a_{(i+j-1)}b)
\\&
=\sum_{n=0}^\infty\sum_{l=0}^{n}\binom{\gamma_b}{n-l}(-1)^{n+1}\partial^{(l)}(a_{(n-1)}b).
\end{align*}
Thus
\begin{align*}&
-(-1)^{p(a)p(b)}b*a - \sum_{n=1}^\infty\sum_{l=1}^n\binom{\gamma_b}{n-l}(-1)^{n+1}
\left\{
\partial^{(l)}(a_{(n-1)}b)-\binom{-\gamma_{a_{(n-1)}b}}{l}a_{(n-1)}b
\right\}\\&
=\sum_{n=1}^\infty\sum_{l=1}^n(-1)^{n+1}\binom{\gamma_b}{n-l}\binom{-\gamma_{a_{(n-1)}b}}{l}a_{(n-1)}b+\sum_{n=0}^\infty(-1)^{n+1}\binom{\gamma_b}{n}a_{(n-1)}b
\\&
=\sum_{n=0}^\infty\sum_{l=0}^n(-1)^{n+1}\binom{\gamma_b}{n-l}\binom{-\gamma_{a_{(n-1)}b}}{l}a_{(n-1)}b
\\&
=\sum_{n=0}^\infty(-1)^{n+1}\binom{-\gamma_a+n-\chi_{a, b}}{n}a_{(n-1)}b
%\\&
=-\sum_{n=0}^\infty\binom{\gamma_a-1+\chi_{a, b}}{n}a_{(n-1)}b
\\&
=-\sum_{n=0}^\infty\binom{\gamma_a}{n}a_{(n-1)}b+(1-\chi_{a,b})\sum_{n=1}^\infty\binom{\gamma_a-1}{n-1}a_{(n-1)}b
\\&
=-a*b+(1-\chi_{a,b})[a_*b]
\end{align*}
Here, we use the trivial identity
\begin{align*}
&\sum_{l=0}^n\binom{\gamma_b}{n-l}\binom{-\gamma_{a_{(n-1)}b}}{l} = \binom{\gamma_b -\gamma_{a_{(n-1)}b}}{n}= \binom{-\gamma_a+n+\chi_{a, b}}{n}.
\end{align*}
\end{proof}
\end{proposition}
\begin{proposition}[{\cite[(2.23)]{DK}}]\label{prop:brA*B*nC}
For $a, b, c \in V$,
\begin{align*}
[a_* \left(Y_*(b, w)c\right)] = (1+ w)^{-\chi_{a, b}}Y_*\left([a_* b], w\right)c + (-1)^{p(a)p(b)}Y_*(b, w)[a_*c],
\end{align*}
or equivalently,
\begin{align*}
[a_*\left( b*_n c\right)]= ([a_*b])*_nc + (-1)^{p(a)p(b)}b*_n[a_*c] + \chi_{a, b}\sum_{j=1}^\infty (-1)^j([a_* b])*_{n-j}c
\end{align*}
for all $n \in \Z$.
\begin{proof}
By Theorem \ref{thm:Borcherds-id} for $n=0$,
$%\displaystyle
a*_{m}Y_*(b, w)c - (-1)^{p(a)p(b)}Y_*(b, w)a*_{m}c
=\sum_{j=0}^\infty(1+w)^{j+1-\chi_{a, b}}\binom{m}{j}Y_*\left(a*_{j}b, w\right)c\ w^{m-j}.
$
Thus, we have
\begin{align*}
&[a_*\left(Y_*(b, w)c\right)] - (-1)^{p(a)p(b)}Y_*(b, w)[a_*c]
\\&
=\sum_{m=0}^\infty(-1)^m
\left(
a*_{m}Y_*(b, w)c-(-1)^{p(a)p(b)}Y_*(b, w)a*_{m}c
\right)
\\&
=\sum_{m, j=0}^\infty (-1)^m(1+ w)^{j+1-\chi_{a, b}}\binom{m}{j}Y_*\left(a*_{j}b, w\right)c\ w^{m-j}
\\&
=(1+w)^{-\chi_{a, b}}\sum_{j=0}^\infty(-1)^j Y_*\left(a*_{j}b, w\right)c
=(1+w)^{-\chi_{a, b}}Y_*\left([a_*b], w\right)c.
\end{align*}
Here, we use the following identity:
$%\displaystyle
\sum_{m=0}^\infty\binom{m}{j}(- w)^m = \frac{(- w)^j}{(1+ w)^{j+1}}.
$
\end{proof}
\end{proposition}

\begin{proposition}[{\cite[(2.38)]{DK}}]\label{prop:*-property}
Let $a, b, c \in V$. For $k \in \Z_{\geq1}$,
\begin{align*}
a*_{-k-1}b = \sum_{j=1}^k \binom{\gamma_a}{k-j}
\left\{
(\partial^{(j)}a)*b - \binom{-\gamma_a}{j}a*b
\right\}.
\end{align*}
\begin{proof}
Apply Corollary \ref{cor:a(-k-1)b-ind}. Then
$%\displaystyle
a*_{-k-1}b
-\sum_{j=1}^k \binom{\gamma_a}{k-j}
\left\{
(\partial^{(j)}a)*b - \binom{-\gamma_a}{j}a*b
\right\}\\
=\sum_{j=0}^k\binom{\gamma_a}{k-j}\binom{-\gamma_a}{j}a*b
=\binom{0}{k}a*b = 0.
$
\end{proof}
\end{proposition}

\hspace{-3.5mm}\textit{Proof of Proposition \ref{prop:J-ideal}.} Although the statements are known in \cite{DK}, we include proof for completeness. The assertion (1) is clear by definition. The assertion (2) follows from \eqref{eq:*_k-ind} and the fact that $\chi_{\partial a, b} = \chi_{a, b}$. Proposition \ref{prop:trans-id} with $b=\vac$ gives $(\partial + H_g)a = a \circ \vac$ for all $a \in V$. If $a \in \gH{V}$, $\epsilon_a = 0$ and $H_g(a) = H(a)$. Thus $(\partial + H)a \in \gH{J}$ for $a \in \gH{V}$. This proves (3). Notice that by (3),
\begin{align}\label{eq:J-ideal1}
\partial^{(n)}a \equiv \binom{-\Delta_a}{n}a \mod \gH{J},\quad
a \in \gH{V}.
\end{align}
The assertion (5) follows from Proposition \ref{prop:A*Bcomm} and \eqref{eq:J-ideal1}. The assertion (6) follows from (1) and Corollary \ref{cor:quasi-asso}. Let $a, b \in V$ such that $\epsilon_a + \epsilon_b \in \Z$ and $c \in \gH{V}$. If $\epsilon_a + \epsilon_b = -1$, we have $\chi_{a, b} = 1$ and $(a*b)*c \equiv a*(b*c)$ mod $\gH{J}$ by (6). Since $\epsilon_c = 0$, we have $\chi_{b, c} = 0$ so that $\epsilon_a + \epsilon_{b*c} = \epsilon_a + \epsilon_b = -1$. Thus $(a*b)*c \in \gH{J}$ and also $c*(a*b) \in \gH{J}$ by (5). If $a, b \in \gH{V}$, $(a \circ b)*c \in \gH{J}$ by (2) and Proposition \ref{prop:AcircB*C}. Moreover, by Proposition \ref{prop:brA*B*nC} with $n=-2$,
\begin{align*}
[c_*(a\circ b)] = ([c_*a])\circ b + (-1)^{p(c)p(a)}a \circ [c_*b] \in \gH{J}.
\end{align*}
Thus $c*(a \circ b) \in \gH{J}$ by (5). This proves (4).\qed

\section{Structure theorem of twisted Zhu algebras}\label{sec:main}

Let $V$ be a vertex superalgebra with a Hamiltonian operator $H$ whose eigenvalues belong to $\R$ and $g \in \operatorname{Aut}_H V$ be a diagonalizable automorphism with modulus $1$ eigenvalues. Recall that $V(\bar{\gamma}, \Delta) = \{ a \in V \mid g(a) = \mathrm{e}^{2\pi i \gamma}a,\ H(a) = \Delta a\}$ for $\bar{\gamma} \in \R/\Z$ and $\Delta \in \R$, and that $\twist = \mathrm{e}^{2\pi i H} \in\operatorname{Aut}_H V$. For $\bar{\mu} \in \R/\Z$, we set
\begin{align*}
V(\bar{\gamma}, \bar{\mu}) = \bigoplus_{\Delta \in \bar{\mu}}V(\bar{\gamma}, \Delta) = \{a \in V \mid g(a) = \mathrm{e}^{2\pi i \gamma}a,\ \twist(a) = \mathrm{e}^{2\pi i \mu}a\}.
\end{align*}
It follows from \eqref{eq:Delta+epsilon}\eqref{eq:gamma} that
\begin{align}\label{eq:Vgamma-mu}
\partial\,V(\bar{\gamma}, \bar{\mu}) \subset V(\bar{\gamma}, \bar{\mu}),\quad
V(\bar{\gamma}_1, \bar{\mu}_1)_{(n)}V(\bar{\gamma}_2, \bar{\mu}_2) \subset V(\bar{\gamma}_1+\bar{\gamma}_2, \bar{\mu}_1+\bar{\mu}_2).
\end{align}
Hence
\begin{align*}
V(\bar{\gamma}_1, \bar{\mu}_1) *_{n}V(\bar{\gamma}_2, \bar{\mu}_2) \subset V(\bar{\gamma}_1+\bar{\gamma}_2, \bar{\mu}_1+\bar{\mu}_2).
\end{align*}
By definition, $\epsilon_a = \epsilon_b$ for all $a, b \in V(\bar{\gamma}, \bar{\mu})$, which we denote by $\epsilon(\bar{\gamma},\bar{\mu})$. Notice that $\epsilon(\bar{\gamma},\bar{\mu}) = 0$ if and only if $\bar{\gamma} = \bar{\mu}$. 
The $(g, H)$-twisted Zhu algebra of $V$ is defined by
\begin{align*}
&\gH{\Zhu}V = \gH{V}/\gH{J},\quad
\gH{V} = \bigoplus_{\bar{\gamma} \in \R/\Z}V(\bar{\gamma}, \bar{\gamma}),\\
&\gH{J} = \gH{V} \circ \gH{V}
+\sum_{
\begin{subarray}{c}
\bar{\gamma}_1+\bar{\gamma}_2 = \bar{\mu}_1+\bar{\mu}_2,\\
\bar{\gamma}_1 \neq \bar{\mu}_1, \bar{\gamma}_2 \neq \bar{\mu}_2
\end{subarray}
}
V(\bar{\gamma}_1, \bar{\mu}_1) * V(\bar{\gamma}_2, \bar{\mu}_2).
\end{align*}

Let $\monoid$ be a discrete additive subgroup of $\R$. Denote by $\monoid_{\geq0} = \monoid \cap \R_{\geq0}$, by $\monoid_{>0} = \monoid \cap \R_{>0}$ and by $\zeta_-$ the largest number in $\monoid_{\geq0}$ that is strictly smaller than $\zeta$ for each $\zeta \in \monoid_{>0}$ as in Section \ref{sec:non-linear salg}. Let $\{e_\alpha\}_{\alpha \in I}$ be a set consisting of homogenous elements in $V$ with respect to $g$, $H$ and the parity. Set $\mathbb{I} = \Z_{\geq0}\times I$ and $e_i = \partial^{(k)}e_\alpha$ for $i=(k, \alpha) \in \mathbb{I}$. %If $I$ is an ordered set, $\mathbb{I}$ is also ordered by \eqref{eq:I_order}.

Suppose that $V$ is freely generated by free generators $\{e_\alpha\}_{\alpha \in I}$ with an ordered index set $I$ and is pregraded by $\monoid_{\geq0}$ in the sense of Section \ref{sec:def-va}. Denote by $\wt(e_\alpha) = \zeta_\alpha \in \monoid_{>0}$ the given $\monoid_{>0}$-grading on $\{e_\alpha\}_{\alpha \in I}$ and set $\base_\zeta = \Span_\C\{e_\alpha \mid \alpha \in I, \zeta_\alpha = \zeta\}$ for $\zeta \in \monoid_{>0}$. We have
\begin{align*}
[{\base_{\zeta_1}}_\lambda\base_{\zeta_2}] \subset \C[\lambda]\otimes V_{(\zeta_1+\zeta_2)_-},\quad
\zeta_1, \zeta_2 \in \monoid_{>0},
\end{align*}
where $V_\zeta$ is an increasing $\monoid_{\geq0}$-filtration on $V$ defined by \eqref{eq:K-filt}. Let
\begin{align*}
\conf = \C[\partial]\otimes\base \subset V,\quad
\base = \bigoplus_{\zeta \in \monoid_{>0}}\base_\zeta = \Span_\C\{e_\alpha\}_{\alpha \in I}.
\end{align*}
By Proposition \ref{prop:DK-V(R)}, $\conf$ is naturally equipped with a structure of a non-linear Lie conformal superalgebra such that $V \simeq V(\conf)$. It follows from \cite[Lemma 6.5]{DK05} that
\begin{align*}
\NO{V_{\zeta_1} V_{\zeta_2}}\ \subset V_{\zeta_1 + \zeta_2},\quad
\partial V_\zeta \subset V_\zeta,\quad
[{V_{\zeta_1}}_\lambda V_{\zeta_2}] \subset V_{(\zeta_1 + \zeta_2)_{-}} \otimes \C[\lambda].
\end{align*}
Therefore
\begin{align}\label{eq:Vzeta1}
\left(V_{\zeta_1}\right)_{(-n-1)}\left(V_{\zeta_2}\right) \subset V_{\zeta_1+\zeta_2},\quad
\left(V_{\zeta_1}\right)_{(n)}\left(V_{\zeta_2}\right) \subset V_{(\zeta_1+\zeta_2)_{-}},\quad
n \geq 0,
\end{align}
which implies that
\begin{align}
\label{eq:Vzeta2}
&\left(V_{\zeta_1}\right)*_{-n-1}\left(V_{\zeta_2}\right) \subset V_{\zeta_1+\zeta_2},\quad
\left(V_{\zeta_1}\right)*_{n}\left(V_{\zeta_2}\right) \subset V_{(\zeta_1+\zeta_2)_{-}},\quad
n \geq 0,\\
\label{eq:Vzeta3}
&[\left(V_{\zeta_1}\right)_*\left(V_{\zeta_2}\right)] \subset V_{(\zeta_1+\zeta_2)_{-}}.
\end{align}
For $\alpha \in I$, there exist $\bar{\gamma}_\alpha \in \R/\Z$ and $\Delta_\alpha \in \R$ such that
\begin{align*}
e_\alpha \in V(\bar{\gamma}_\alpha, \Delta_\alpha) \subset V(\bar{\gamma}_\alpha, \bar{\mu}_\alpha),\quad
\bar{\mu}_\alpha = \Delta_\alpha + \Z \in \R/\Z.
\end{align*}
Set $\epsilon_\alpha = \epsilon(\bar{\gamma}_\alpha, \bar{\mu}_\alpha)$. For $\overrightarrow{i} = (i_1, \ldots, i_s) \in \mathbb{I}^s$ with $i_j = (k_j, \alpha_j) \in \mathbb{I}$, define
\begin{align*}
&\zeta(\overrightarrow{i}) = \zeta_{\alpha_1} + \cdots + \zeta_{\alpha_s},\quad
\bar{\gamma}(\overrightarrow{i}) = \bar{\gamma}_{\alpha_1} + \cdots + \bar{\gamma}_{\alpha_s},\\
&\Delta(\overrightarrow{i}) = \Delta_{\alpha_1} + k_1 + \cdots + \Delta_{\alpha_s} + k_s,\quad
\bar{\mu}(\overrightarrow{i}) = \Z + \Delta(\overrightarrow{i}) = \bar{\mu}_{\alpha_1} + \cdots + \bar{\mu}_{\alpha_s},\\
&\epsilon(\overrightarrow{i}) = \epsilon\left(\bar{\gamma}(\overrightarrow{i}), \bar{\mu}(\overrightarrow{i})\right),\quad
\gamma(\overrightarrow{i}) = \Delta(\overrightarrow{i}) + \epsilon(\overrightarrow{i}),\\
&\#(\overrightarrow{i}) = \# \left\{ (p, q) \mid 1 \leq p < q \leq s\ \mathrm{and}\ \mathrm{either}\ i_p > i_q,\ \mathrm{or}\ i_p = i_q\ \mathrm{and}\ p(e_{i_p}) = \bar{1}\right\}.
\end{align*}
Since $e_\alpha \in V_{\zeta_\alpha}$, we have $\NO{e_{i_1} \cdots e_{i_s}}\ \in V_{\zeta(\overrightarrow{i})}$. By \eqref{eq:Vgamma-mu}, we also have $\bar{\gamma}(\overrightarrow{i}) = \gamma(\overrightarrow{i}) + \Z$ and $\NO{e_{i_1} \cdots e_{i_s}}\ \in V(\bar{\gamma}(\overrightarrow{i}), \Delta(\overrightarrow{i}))$. Furthermore,
\begin{align*}
\{ \NO{e_{i_1} \cdots e_{i_s}}\, \mid \#(\overrightarrow{i}) = 0 \}
\end{align*}
is a basis of $V$. Define subsets $B_\zeta$, $B_\zeta^*$ in $V_\zeta$ for $\zeta \in \monoid$ by
\begin{align*}
&B_\zeta = \{\NO{e_{i_1} \cdots e_{i_s}}\, \mid \zeta(\overrightarrow{i}) \leq \zeta,\  \#(\overrightarrow{i})  = 0\},\\
&B^*_\zeta = \{ e_{i_1} * \cdots * e_{i_s} \mid \zeta(\overrightarrow{i}) \leq \zeta,\  \#(\overrightarrow{i})  = 0\}.
\end{align*}
Then $\{ \NO{e_{i_1} \cdots e_{i_s}}\, \mid \#(\overrightarrow{i}) = 0 \} = \bigcup_{\zeta \in \monoid_{\geq0}} B_\zeta$. Similarly we define $B^* = \bigcup_{\zeta \in \monoid_{\geq0}} B_\zeta^*$, where $e_{i_1} * \cdots * e_{i_s}$ is defined by taking products from right to left: $e_{i_1} * \cdots * e_{i_s}  = e_{i_1} * ( e_{i_2} * (\cdots * e_{i_s} ))$. By definition, we have
\begin{align*}
V_\zeta = \Span_\C\{ \NO{e_{i_1} \cdots e_{i_s}}\, \mid \zeta(\overrightarrow{i}) \leq \zeta \}.
\end{align*}
Similarly, we define $V_\zeta^*$ by $*$-products instead of normally ordered products in $V_\zeta$
\begin{align*}
V^*_\zeta = \Span_\C \{ e_{i_1} * \cdots * e_{i_s} \mid \zeta(\overrightarrow{i}) \leq \zeta \}.
\end{align*}
De Sole and Kac proved the following lemma for $g = \twist$, but the same proof applies.
\begin{lemma}[{\cite[Lemma 3.18]{DK}}]\label{lem:PBW}\ \\
\begin{enumerate}
\item For every $\overrightarrow{i} = (i_1, \ldots, i_s)$, we have
\begin{align*}
e_{i_1} * \cdots * e_{i_s} - \NO{e_{i_1} \cdots e_{i_s}}\ \in V_{\zeta(\overrightarrow{i})_{-}}.
\end{align*}
\item $V_\zeta = \Span_\C B_\zeta = V^*_\zeta = \Span_\C B^*_\zeta$.\\
\item $B^*$ is a basis of $V$.
\end{enumerate}
\end{lemma}
For any subspace $A$ of $V$, we use the following notations
\begin{align*}
A(\bar{\gamma}, \bar{\mu}) = A \cap V(\bar{\gamma}, \bar{\mu}),\quad
A_\zeta = A \cap V_\zeta,\quad
A(\bar{\gamma}, \bar{\mu})_\zeta = A \cap V(\bar{\gamma}, \bar{\mu}) \cap V_\zeta
\end{align*}
for $\bar{\gamma}, \bar{\mu} \in \R/\Z$ and $\zeta \in \monoid_{\geq0}$. The following are analogs of \cite[Lemma 3.20]{DK} and can be proved similarly by using results in Section \ref{sec:properties_*}.
\begin{lemma}\label{lem:Vzeta}
Let $a \in V(\bar{\gamma}_1, \bar{\mu}_1)_{\zeta_1}$, $b \in V(\bar{\gamma}_2, \bar{\mu}_2)_{\zeta_2}$ and $c \in V(\bar{\gamma}_3, \bar{\mu}_3)_{\zeta_3}$ for $\bar{\gamma}_j, \bar{\mu}_j \in \R/\Z$, $\zeta_j \in \monoid_{\geq0}$ for $j=1, 2, 3$. Then we have the following properties.
\begin{enumerate}
\item For $k, h \in \Z_{\geq0}$,
\begin{align}
\label{eq:der-*-eq1}&(\partial^{(k)}a)*b-\binom{-\gamma_a}{k} a*b \in V(\bar{\gamma}_1, \bar{\mu}_1)_{\zeta_1} \circ V(\bar{\gamma}_2, \bar{\mu}_2)_{\zeta_2},\\
\label{eq:der-*-eq2}&[\partial^{(h)}a_*\partial^{(k)}b]-\binom{-\gamma_a}{h}\binom{-\gamma_b-\chi_{a,b}}{k}[a_*b]\in V(\bar{\gamma}_1+\bar{\gamma}_2, \bar{\mu}_1+\bar{\mu}_2)_{(\zeta_1+\zeta_2)_{-}} \circ \vac.
\end{align}
\item For $k \in \Z_{\geq1}$, $a *_{-k-1} b \in V(\bar{\gamma}_1, \bar{\mu}_1)_{\zeta_1} \circ V(\bar{\gamma}_2, \bar{\mu}_2)_{\zeta_2}$.\\
\item $(a*b)*c-a*(b*c) \in \Bigl(V(\bar{\gamma}_1, \bar{\mu}_1) *_{-2+\chi_{a,b}}V(\bar{\gamma}_2+\bar{\gamma}_3, \bar{\mu}_2+\bar{\mu}_3)+V(\bar{\gamma}_2, \bar{\mu}_2)*_{-2+\chi_{a,b}}V(\bar{\gamma}_1+\bar{\gamma}_3, \bar{\mu}_1+\bar{\mu}_3)\Bigr)_{(\zeta_1+\zeta_2+\zeta_3)_{-}}$.\\
\item $a * b - (-1)^{p(a)p(b)}b * a-(1-\chi_{a,b})[a_*b] \in V(\bar{\gamma}_1+\bar{\gamma}_2, \bar{\mu}_1+\bar{\mu}_2)_{(\zeta_1+\zeta_2)_{-}} \circ \vac$.
\item $(a \circ b) * c \in \Bigl(V(\bar{\gamma}_1, \bar{\mu}_1) *_{-2+\chi_{a, b}}V(\bar{\gamma}_2+\bar{\gamma}_3, \bar{\mu}_2+\bar{\mu}_3) + V(\bar{\gamma}_2, \bar{\mu}_2) *_{-2+\chi_{a, b}}V(\bar{\gamma}_1+\bar{\gamma}_3, \bar{\mu}_1+\bar{\mu}_3)\Bigr)_{\zeta_1+\zeta_2+\zeta_3}$.\\
\item $[a_*(b\circ c)] \in \Bigl(V(\bar{\gamma}_1+\bar{\gamma}_2, \bar{\mu}_1+\bar{\mu}_2) \circ V(\bar{\gamma}_3, \bar{\mu}_3) + V(\bar{\gamma}_2, \bar{\mu}_2) \circ V(\bar{\gamma}_1+\bar{\gamma}_3, \bar{\mu}_1+\bar{\mu}_3)\Bigr)_{(\zeta_1+\zeta_2+\zeta_3)_{-}}$.\\
\item If two of $a, b, c$ belong to $\gH{V}$,
\begin{align*}
&a * (b \circ c) \in
\Bigl(V(\bar{\gamma}_2, \bar{\mu}_2)\circ V(\bar{\gamma}_1+\bar{\gamma}_3, \bar{\mu}_1+\bar{\mu}_3)
+ V(\bar{\gamma}_3, \bar{\mu}_3) \circ V(\bar{\gamma}_1+\bar{\gamma}_2, \bar{\mu}_1+\bar{\mu}_2)\\
&+ V(\bar{\gamma}_1+\bar{\gamma}_2, \bar{\mu}_1+\bar{\mu}_2) \circ V(\bar{\gamma}_3, \bar{\mu}_3)
+ V(\bar{\gamma}_1+\bar{\gamma}_2+\bar{\gamma}_3, \bar{\mu}_1+\bar{\mu}_2+\bar{\mu}_3) \circ \vac\Bigr)_{\zeta_1+\zeta_2+\zeta_3}.
\end{align*}
\item If either $a \circ b$ or $c$ belongs to $\gH{V}$,
\begin{align*}
[(a \circ b)_* c] \in
&\Bigl(V(\bar{\gamma}_1, \bar{\mu}_1) \circ V(\bar{\gamma}_2+\bar{\gamma}_3, \bar{\mu}_2+\bar{\mu}_3) + V(\bar{\gamma}_1+\bar{\gamma}_3, \bar{\mu}_1+\bar{\mu}_3) \circ V(\bar{\gamma}_2, \bar{\mu}_2)\\
&+ V(\bar{\gamma}_1+\bar{\gamma}_2+\bar{\gamma}_3, \bar{\mu}_1+\bar{\mu}_2+\bar{\mu}_3) \circ \vac\Bigr)_{(\zeta_1+\zeta_2+\zeta_3)_{-}}.
\end{align*}
\end{enumerate}
\begin{proof}
Proposition \ref{prop:trans-id} (1) shows \eqref{eq:der-*-eq1} for $k=1$. For $k \geq 2$,
\begin{align*}
&(\partial^{(k)}a)*b-\binom{-\gamma_a}{k} a*b\\
&= \frac{1}{k}\left(
(\partial^{(k-1)}\partial a) * b - \binom{-\gamma_a-1}{k-1}(\partial a)*b
\right) + \frac{1}{k}\binom{-\gamma_a-1}{k-1}a\circ b,
\end{align*}
so belongs to $V(\bar{\gamma}_1, \bar{\mu}_1)_{\zeta_1} \circ V(\bar{\gamma}_2, \bar{\mu}_2)_{\zeta_2}$ by inductive assumptions on $k$ and \eqref{eq:Vgamma-mu}. To prove \eqref{eq:der-*-eq2}, by Proposition \ref{prop:trans-id} (2), we have $[\partial^{(h)}a_*b] = \binom{-\gamma_a}{h}[a_*b]$. Thus, it suffices to show \eqref{eq:der-*-eq2} for $h = 0$. For $k=1$, $[a_*\partial b] + (\gamma_b+\chi_{a,b})[a_*b]=[a_*b]\circ\vac \in V(\bar{\gamma}_1+\bar{\gamma}_2, \bar{\mu}_1+\bar{\mu}_2)_{(\zeta_1+\zeta_2)_{-}} \circ \vac$ by Proposition \ref{prop:trans-id} (3) and \eqref{eq:Vgamma-mu}\eqref{eq:Vzeta3}. For $k \geq 2$,
\begin{align*}
&[a_*\partial^{(k)}b]-\binom{-\gamma_b-\chi_{a,b}}{k}[a_*b]\\
&=\frac{1}{k}\left(
[a_*\partial^{(k-1)}\partial b]-\binom{-\gamma_b-\chi_{a,b}-1}{k-1}[a_*\partial b]
\right)
+\frac{1}{k}\binom{-\gamma_b-\chi_{a,b}-1}{k-1}[a_*b]\circ\vac,
\end{align*}
so belongs to $V(\bar{\gamma}_1+\bar{\gamma}_2, \bar{\mu}_1+\bar{\mu}_2)_{(\zeta_1+\zeta_2)_{-}} \circ \vac$ by inductive assumption on $k$ and \eqref{eq:Vgamma-mu}. The assertion (2) follows from Proposition \ref{prop:*-property} and \eqref{eq:der-*-eq1}. The assertion (3) follows from Corollary \ref{cor:quasi-asso} and \eqref{eq:Vgamma-mu}\eqref{eq:Vzeta2}. The assertion (4) follows from Proposition \ref{prop:A*Bcomm} and \eqref{eq:Vgamma-mu}\eqref{eq:Vzeta2}\eqref{eq:der-*-eq1}. The assertion (5) follows from Proposition \ref{prop:AcircB*C} and \eqref{eq:Vgamma-mu}\eqref{eq:Vzeta2}. The assertion (6) follows from Proposition \ref{prop:brA*B*nC} and  \eqref{eq:Vgamma-mu}\eqref{eq:Vzeta2}\eqref{eq:Vzeta3}. To prove (7), if two of $a, b, c$ belong to $\gH{V}$, two of $\epsilon_a, \epsilon_b, \epsilon_c$ are equal to zero so that $\chi_{a, b\circ c} = \chi_{b, c} =0$. Then
\begin{align*}
a*(b\circ c) &=  (-1)^{p(a)p(b\circ c)}(b \circ c) * a + [a_*(b \circ c)]\\
&+\bigl(a*(b\circ c) - (-1)^{p(a)p(b\circ c)}(b \circ c)*a - [a_*(b\circ c)]\bigr).
\end{align*}
Thus, assertion (7) follows from assertion (4), as well as (5) and (6). To prove (8), if either $a \circ b$ or $c$ belongs to $\gH{V}$, either $\epsilon_{a\circ b}$ or $\epsilon_c$ are equal to zero so that $\chi_{a \circ b, c} = 0$. Since
\begin{align*}
&[(a \circ b)_* c] = -\bigl((a \circ b)*c-(-1)^{p(a \circ b)p(c)}c*(a \circ b)-[(a\circ b)_*c]\bigr)\\
&-(-1)^{p(a \circ b)p(c)}\bigl(c*(a\circ b)-(-1)^{p(c)p(a \circ b)}(a\circ b)*c - [c_*(a\circ b)]\bigr)\\
&-(-1)^{p(a \circ b)p(c)}[c_*(a\circ b)],
\end{align*}
the assertion (8) follows from (4)(6).
\end{proof}
\end{lemma}

Let
\begin{align*}
\Sigma_\zeta = \Sigma_{1, \zeta} \sqcup \Sigma_{2, \zeta},\quad
W_\zeta = W_{1, \zeta} + W_{2, \zeta},\quad
\zeta \in \monoid_{\geq0},
\end{align*}
where
\begin{align*}
&\Sigma_{1, \zeta} =
\ \left\{\ 
\begin{lgathered}
\partial^{(k_1)}e_{\alpha_1}*\cdots*\partial^{(k_s)}e_{\alpha_s}\\
- \binom{-\Delta_1}{k_1}\cdots\binom{-\Delta_s}{k_s}e_{\alpha_1}*\cdots*e_{\alpha_s}
\end{lgathered}
\ \middle|\ 
\begin{lgathered}
s \geq1,\ 
\sum_{j=1}^s k_j \geq1,\ 
i_j = (k_j, \alpha_j) \in \mathbb{I},\\
\zeta(\overrightarrow{i}) \leq \zeta,\ \#\overrightarrow{i}=0\ 
\mathrm{for}\ 
\overrightarrow{i} = (i_1, \ldots, i_s),\\
\Delta_j = \Delta_{\alpha_j},\ 
e_{\alpha_j} \in \gH{V}\ 
\mathrm{for}\ \mathrm{all}\ j
\end{lgathered}
\ \right\},\\
&\Sigma_{2, \zeta} =
\ \left\{\ 
e_{i_1}*\cdots*e_{i_s}
\ \middle|\ 
\begin{lgathered}
s >1,\ 
i_j \in \mathbb{I},\ 
\zeta(\overrightarrow{i}) \leq \zeta,\ 
\epsilon(\overrightarrow{i}) =0,\ 
\#\overrightarrow{i}=0\\
\mathrm{for}\ 
\overrightarrow{i} = (i_1, \ldots, i_s),\ 
{}^\exists j\ \mathrm{such}\ \mathrm{that}\ e_{i_j} \notin \gH{V}
\end{lgathered}
\ \right\}
\end{align*}
and
\begin{align*}
&W_{1, \zeta} =
\Span_\C\left\{\ 
\begin{lgathered}
\partial^{(k_1)}e_{\alpha_1}*\cdots*\partial^{(k_s)}e_{\alpha_s}\\
- \binom{-\Delta_1}{k_1}\cdots\binom{-\Delta_s}{k_s}e_{\alpha_1}*\cdots*e_{\alpha_s}
\end{lgathered}
\ \middle|\ 
\begin{lgathered}
s \geq1,\ 
i_j = (k_j, \alpha_j) \in \mathbb{I},\\
\zeta(\overrightarrow{i}) \leq \zeta\ 
\mathrm{for}\ 
\overrightarrow{i} = (i_1, \ldots, i_s),\\
\Delta_j = \Delta_{\alpha_j},\ 
e_{\alpha_j} \in \gH{V} \ 
\mathrm{for}\ \mathrm{all}\ j
\end{lgathered}
\ \right\},\\
&W_{2, \zeta} =
\Span_\C\left\{\ 
e_{i_1}*\cdots*e_{i_s}
\ \middle|\ 
\begin{lgathered}
s >1,\ 
i_j \in \mathbb{I},\ 
\zeta(\overrightarrow{i}) \leq \zeta,\ 
\epsilon(\overrightarrow{i}) =0\\
\mathrm{for}\ 
\overrightarrow{i} = (i_1, \ldots, i_s),\ 
{}^\exists j\ \mathrm{such}\ \mathrm{that}\ e_{i_j} \notin \gH{V}
\end{lgathered}
\ \right\}.
\end{align*}
Set $\Sigma = \bigcup_{\zeta \in \monoid_{\geq 0}}\Sigma_{\zeta}$ and $W = \bigcup_{\zeta \in \monoid_{\geq0}}W_{\zeta}$. Obviously, $\Span_\C\Sigma_\zeta \subset W_\zeta$. Recall that $\gamma_a = \Delta_a$ for $a \in \gH{V}$.

\begin{lemma}[{\cite[Lemma 3.19]{DK}}]\label{lem:e*W_1}
Let $i \in \mathbb{I}$ such that $e_i \in \gH{V}$ and $\zeta(i) = \zeta_0$. Then $e_i * W_{1, \zeta} \subset W_{1, \zeta_0+\zeta}$.
\end{lemma}

\begin{lemma}\label{lem:W-ideal}
Let $i \in \mathbb{I}$ such that $\zeta(i) = \zeta_0$.\\ 
(1) If $e_i \in \gH{V}$, then $e_i * W_{\zeta} \subset W_{\zeta_0 + \zeta}$.\\
(2) If $e_i \in \gH{V}$, then $e_i \circ (\gH{V})_{\zeta} \subset W_{\zeta_0 + \zeta}$.\\
(3) If $e_i \in V(\bar{\gamma}_1, \bar{\mu}_1)$ and $\bar{\gamma}_1 \neq \bar{\mu}_1$, then $e_i *_{-k-1} V(\bar{\gamma}_2, \bar{\mu}_2)_{\zeta} \subset W_{\zeta_0 + \zeta}$ for $k \geq0$, $\bar{\gamma}_2, \bar{\mu}_2 \in \R/\Z$ such that $\bar{\gamma}_1+\bar{\gamma}_2 = \bar{\mu}_1+\bar{\mu}_2$.
\begin{proof}
(1) We have $e_i * W_{1, \zeta} \subset W_{1, \zeta_0+\zeta}$ by Lemma \ref{lem:e*W_1}. It is also clear that $e_i * W_{2, \zeta} \subset W_{2, \zeta_0 + \zeta}$ by definition. Thus, the assertion follows.\\
(2) By Lemma \ref{lem:PBW},
\begin{align}\label{eq:VZzeta}
(\gH{V})_{\zeta} =
\Span_\C\left\{e_{i_1} * \cdots * e_{i_s} 
\middle|\ 
\begin{lgathered}
s \geq0,\ 
i_j \in \mathbb{I},\ 
\zeta(\overrightarrow{i}) \leq \zeta,\ 
\epsilon(\overrightarrow{i}) =0,\\
\#\overrightarrow{i}=0\ 
\mathrm{for}\ 
\overrightarrow{i} = (i_1, \ldots, i_s)
\end{lgathered}
\ \right\}.
\end{align}
Thus it is enough to show that $e_i \circ (e_{i_1} * \cdots * e_{i_s}) \in W_{\zeta_0 + \zeta}$ for $e_{i_1} * \cdots * e_{i_s}$ such that $\zeta(\overrightarrow{i}) \leq \zeta$, $\epsilon(\overrightarrow{i}) =0$ and $\#(\overrightarrow{i})=0$. If $e_{i_j} \in \gH{V}$ for all $j$, $e_i \circ (e_{i_1} * \cdots * e_{i_s}) \in W_{1,\zeta_1 + \zeta}$ by the proof of \cite[Lemma 3.21]{DK}. Thus, we only need to consider the case that there exists $j$ such that $e_{i_j} \notin \gH{V}$. Then, by Proposition \ref{prop:trans-id} (1),
\begin{align}\label{eq:W-ideal}
e_i \circ (e_{i_1} * \cdots * e_{i_s}) = (\partial e_i + \gamma_{e_i} e_i)*e_{i_1} * \cdots * e_{i_s} \in W_{2, \zeta_0 + \zeta}.
\end{align}
This completes the proof.\\
(3) Similarly, by Lemma \ref{lem:PBW},
\begin{align}\label{eq:Vgammazeta}
V(\bar{\gamma}, \bar{\mu})_{\zeta} =
\Span_\C\left\{e_{i_1} * \cdots * e_{i_s} 
\middle|\ 
\begin{lgathered}
s \geq0,\ 
i_j \in \mathbb{I},\ 
\zeta(\overrightarrow{i}) \leq \zeta,\ 
\bar{\gamma}(\overrightarrow{i}) = \bar{\gamma},\\
\bar{\mu}(\overrightarrow{i}) = \bar{\mu},\ 
\#\overrightarrow{i}=0\ 
\mathrm{for}\ 
\overrightarrow{i} = (i_1, \ldots, i_s)
\end{lgathered}
\ \right\}.
\end{align}
Let $e_{i_1} * \cdots * e_{i_s}$ such that $\zeta(\overrightarrow{i}) \leq \zeta$, $\bar{\gamma}(\overrightarrow{i}) = \bar{\gamma}_2$, $\bar{\mu}(\overrightarrow{i}) = \bar{\mu}_2$ and  $\#(\overrightarrow{i})=0$. Then it is enough to show that $e_i *_{-k-1}(e_{i_1} * \cdots * e_{i_s}) \in W_{2, \zeta_0 + \zeta}$ for all $k \geq 0$ by induction on $k$. For $k =0$, the statement is clear by definition. For $k \geq1$, the statement is also immediate from \eqref{eq:*_k-ind} and the inductive assumption on $k$. This completes the proof.
\end{proof}
\end{lemma}

Let
\begin{align*}
(\gH{J})_\zeta = \sum_{
\begin{subarray}{c}
\zeta_1, \zeta_2 \in \monoid_{\geq 0}\\
\zeta_1 + \zeta_2 \leq \zeta
\end{subarray}}
\Bigl(
(\gH{V})_{\zeta_1} \circ (\gH{V})_{\zeta_2}
+\sum_{\begin{subarray}{c} \bar{\gamma}_1 + \bar{\gamma}_2 = \bar{\mu}_1 + \bar{\mu}_2 \\ \bar{\gamma}_1 \neq \bar{\mu}_1, \bar{\gamma}_2 \neq \bar{\mu}_2 \end{subarray}} V(\bar{\gamma}_1, \bar{\mu}_1)_{\zeta_1} * V(\bar{\gamma}_2, \bar{\mu}_2)_{\zeta_2}
\Bigr).
\end{align*}
Then $\gH{J} = \sum_{\zeta \in \monoid_{\geq 0}}(\gH{J})_\zeta$.

\begin{lemma}\label{lem:J-trivial}
Let $a \in V_{\zeta_1}$ and $b \in V_{\zeta_2}$ such that $a *_{-k-1} b \in \gH{V}$ for some $k \geq 1$. Then $a *_{-k-1} b \in (\gH{J})_{\zeta_1+\zeta_2}$.
\begin{proof}
We will show the assertion by induction on $k$. We may suppose that $a \in V(\bar{\gamma}_1, \bar{\mu}_1)_{\zeta_1}$ and $b \in V(\bar{\gamma}_2, \bar{\mu}_2)_{\zeta_2}$ such that $\bar{\gamma}_1 + \bar{\gamma}_2 = \bar{\mu}_1 + \bar{\mu}_2$. Consider the case $k=1$. If $\bar{\gamma}_1=\bar{\mu}_1$, the assertion follows from Lemma \ref{lem:Vzeta} (2). If $\bar{\gamma}_1 \neq \bar{\mu}_1$, we have $a\circ b = (\partial a + \gamma_a a) * b \in V(\bar{\gamma}_1, \bar{\mu}_1)_{\zeta_1} * V(\bar{\gamma}_2, \bar{\mu}_2)_{\zeta_2} \subset (\gH{J})_{\zeta_1+\zeta_2}$ by Proposition \ref{prop:trans-id} (1) and 
\eqref{eq:Vzeta2}. For $k\geq2$, the assertion is immediate from \eqref{eq:*_k-ind}\eqref{eq:Vgamma-mu}\eqref{eq:Vzeta2} and the inductive assumption on $k$. This completes the proof.
\end{proof}
\end{lemma}

\begin{lemma}\label{lem:J-ideal}
$a * (\gH{J})_\zeta \subset (\gH{J})_{\zeta_0 + \zeta}$ and $(\gH{J})_\zeta * a \subset (\gH{J})_{\zeta_0 + \zeta}$ for all $a \in (\gH{V})_{\zeta_0}$.
\begin{proof}
The assertion is immediate from Proposition \ref{prop:J-ideal} and \eqref{eq:Vzeta2}.
\end{proof}
\end{lemma}

\begin{lemma}\label{lem:J-or-Vodd}
Let $a_j \in V(\bar{\gamma}_j, \bar{\mu}_j)_{\zeta_j}$ for $j =1 , \ldots, s$ such that $\sum_{j=1}^s\bar{\gamma}_j = \sum_{j=1}^s\bar{\mu}_j$. Suppose that $\bar{\gamma}_j \neq \bar{\mu}_j$ or $a_j \in (\gH{J})_{\zeta_j}$ for some $j$. Then $a_1 * \cdots * a_s \in (\gH{J})_{\zeta_1 + \cdots + \zeta_s}$.
\begin{proof}
First, consider the case that $\bar{\gamma}_j \neq \bar{\mu}_j$ for some $j$. Let $j_0 = \min\{j \mid a_j \notin \gH{V} \}$. Then $a_j \in \gH{V}$ for all $j = 1, \ldots, j_0 -1$ and $a_{j_0+1} * \cdots * a_s \in V(\bar{\gamma}, \bar{\mu})$ such that $\bar{\gamma}_{j_0}+\bar{\gamma} = \bar{\mu}_{j_0}+\bar{\mu}$ and $\bar{\gamma} \neq \bar{\mu}$ by \eqref{eq:Vgamma-mu}. Thus, by Lemma \ref{lem:J-ideal} and \eqref{eq:Vzeta2},
\begin{align*}
a_1 * \cdots * a_s \in a_1 * \cdots * a_{j_0 -1} * (\gH{J})_{\zeta_{j_0} + \cdots + \zeta_s}
\subset (\gH{J})_{\zeta_1 + \cdots + \zeta_s}.
\end{align*}
Next, we suppose that $a_i \in \gH{V}$ for all $i$ and $a_j \in (\gH{J})_{\zeta_j}$ for some $j$. Then again by Lemma \ref{lem:J-ideal}, $a_1 * \cdots * a_s \in (\gH{J})_{\zeta_1 + \cdots + \zeta_s}$. This completes the proof.
\end{proof}
\end{lemma}

\begin{lemma}\label{lem:W=J}
$W_\zeta = (\gH{J})_\zeta$.
\begin{proof}
Using \eqref{eq:der-*-eq1}, Lemma \ref{lem:Vzeta} (7) and the proof of \cite[lemma 3.21]{DK}, it follows that $W_{1, \zeta} \subset (\gH{V} \circ \gH{V})_{\zeta}$. By Lemma \ref{lem:J-or-Vodd}, we also have $W_{2,\zeta} \subset (\gH{J})_\zeta$. Thus $W_\zeta \subset (\gH{J})_\zeta$. Now it suffices to show that $(\gH{J})_\zeta \subset W_\zeta$ by induction on $\zeta$.

First of all, we will show that $(\gH{V})_{\zeta_1} \circ (\gH{V})_{\zeta_2} \subset W_{\zeta}$ for $\zeta_1, \zeta_2 \in \monoid_{\geq0}$ with $\zeta_1+\zeta_2 = \zeta$. By \eqref{eq:VZzeta}, it suffices to show that
\begin{align}\label{eq:W=J}
(e_{i_1} * \cdots * e_{i_s}) \circ (e_{j_1} * \cdots * e_{j_t}) \in W_\zeta
\end{align}
for $i_\ell, j_\ell \in \mathbb{I}$ such that $\zeta(\overrightarrow{i}) = \zeta_1$, $\zeta(\overrightarrow{j}) = \zeta_2$, $\epsilon(\overrightarrow{i}) = \epsilon(\overrightarrow{j})=0$ and $\#(\overrightarrow{i}) = \#(\overrightarrow{j}) = 0$ for $\overrightarrow{i} = (i_1, \ldots, i_s)$, $\overrightarrow{j} = (j_1, \ldots, j_t)$ with $\zeta_1 + \zeta_2 = \zeta$ under the inductive assumption that $(\gH{J})_{\zeta'} = W_{\zeta'}$ for all $\zeta' \in \monoid_{\geq0}$ such that $\zeta' < \zeta$. For $s=1$, $e_{i_1} \in (\gH{V})_{\zeta_1}$ and $e_{j_1} * \cdots * e_{j_t} \in (\gH{V})_{\zeta_2}$. Hence \eqref{eq:W=J} holds by Lemma \ref{lem:W-ideal} (2). For $s \geq2$, we set $a=e_{i_1}$, $b=e_{i_2}* \cdots * e_{i_s}$ and $c=e_{j_1}* \cdots * e_{j_t}$. Then
\begin{align*}
(e_{i_1} * \cdots * e_{i_s}) \circ (e_{j_1} * \cdots * e_{j_t})
= (a*b)\circ c.
\end{align*}
We show that $(a*b)\circ c \in W_\zeta$. If $a = e_{i_1} \in \gH{V}$, we have $b \in \gH{V}$ and using Proposition \ref{prop:Borcherds-id} with $n=-1$ and $k=-2$,
\begin{align*}
&(e_{i_1} * \cdots * e_{i_s}) \circ (e_{j_1} * \cdots * e_{j_t})
= (a*b)\circ c\\
&=a*(b\circ c) + a \circ(b*c) + \sum_{j=0}^\infty\left(
a*_{-j-3}(b*_j c) +(-1)^{p(a)p(b)} b*_{-j-3}(a*_j c)
\right).
\end{align*}
First of all, since $b\circ c \in (\gH{J})_{\zeta-\zeta(i_1)} = W_{\zeta-\zeta(i_1)}$ by the inductive assumption, $a * (b \circ c) \in W_{\zeta}$ by Lemma \ref{lem:W-ideal} (1). Next, since $b* c \in (\gH{V})_{\zeta-\zeta(i_1)}$, we also have $a \circ (b*c) \in W_\zeta$ by Lemma \ref{lem:W-ideal} (2). Finally, since $b*_j c \in V_{(\zeta-\zeta(i_1))_{-}}$ and $a*_j c \in V_{(\zeta(i_1)+\zeta_2)_{-}}$ for all $j \geq 0$ by \eqref{eq:Vzeta2}, we have $a*_{-j-3}(b*_jc), b*_{-j-3}(a*_jc) \in (\gH{J})_{\zeta_{-}} = W_{\zeta_{-}} \subset W_\zeta$ by Lemma \ref{lem:J-trivial} and the inductive assumption. Thus $(a*b)\circ c \in W_\zeta$ for $a \in \gH{V}$. If $a=e_{i_1} \in V(\bar{\gamma}, \bar{\mu})$ with $\bar{\gamma} \neq \bar{\mu}$, we have $\epsilon_a\neq 0 $ and $\epsilon_b \neq 0$ so that $\chi_{a,b} =1$. Then using Proposition \ref{prop:Borcherds-id} with $n=-1$ and $k=-2$ again,
\begin{align*}
&(e_{i_1} * \cdots * e_{i_s}) \circ (e_{j_1} * \cdots * e_{j_t})
=(a*b)\circ c\\
&=a*(b\circ c) + a \circ (b*c) + a*(b*c) \\
&+\sum_{j=0}^\infty\left(
a*_{-j-3}(b*_jc)+(-1)^{p(a)p(b)}b*_{-j-3}(a*_jc)
\right)\\
&+\sum_{j=0}^\infty\left(
a*_{-j-2}(b*_jc) + +(-1)^{p(a)p(b)}b*_{-j-2}(a*_jc)
\right).
\end{align*}
Then all $a*(b\circ c)$, $a\circ(b*c)$, $a*(b*c)$ and $a*_{-j-2}(b*_jc)$ for $j\geq0$ belong to $W_{\zeta}$ by Lemma \ref{lem:W-ideal} (3). Moreover, since $a*_jc \in V(\bar{\gamma}, \bar{\mu})_{(\zeta(i_1)+\zeta_2)_{-}}$ by \eqref{eq:Vgamma-mu}\eqref{eq:Vzeta2}, we have $b*_{-j-2}(a*_jc) \in (\gH{J})_{\zeta_{-}} = W_{\zeta_{-}} \subset W_\zeta$ by Lemma \ref{lem:J-trivial} and the inductive assumption on $\zeta$. Therefore, we obtain \eqref{eq:W=J}.

Next, we will show that $V(\bar{\gamma}_1, \bar{\mu}_1)_{\zeta_1} * V(\bar{\gamma}_2, \bar{\mu}_2)_{\zeta_2} \subset W_{\zeta}$ for $\zeta_1, \zeta_2 \in \monoid_{\geq0}$ with $\bar{\gamma}_1 \neq \bar{\mu}_1$, $\bar{\gamma}_2 \neq \bar{\mu}_2$, $\bar{\gamma}_1 + \bar{\gamma}_2 = \bar{\mu}_1 + \bar{\mu}_2$ and $\zeta_1+\zeta_2 = \zeta$. By \eqref{eq:Vgammazeta}, it suffices to show that
\begin{align}\label{eq:W=J_2}
(e_{i_1} * \cdots * e_{i_s}) * (e_{j_1} * \cdots * e_{j_t}) \in W_\zeta
\end{align}
for $i_\ell, j_\ell \in \mathbb{I}$ such that $\zeta(\overrightarrow{i}) = \zeta_1$, $\zeta(\overrightarrow{j}) = \zeta_2$, $\bar{\gamma}(\overrightarrow{i}) = \bar{\gamma}_1$, $\bar{\gamma}(\overrightarrow{j}) = \bar{\gamma}_2$, $\bar{\mu}(\overrightarrow{i}) = \bar{\mu}_1$, $\bar{\mu}(\overrightarrow{j}) = \bar{\mu}_2$ and $\#(\overrightarrow{i}) = \#(\overrightarrow{j}) = 0$ for $\overrightarrow{i} = (i_1, \ldots, i_s)$, $\overrightarrow{j} = (j_1, \ldots, j_t)$ with $\zeta_1 + \zeta_2 = \zeta$ under the inductive assumption that $(\gH{J})_{\zeta'} = W_{\zeta'}$ for all $\zeta' \in \monoid_{\geq0}$ such that $\zeta' < \zeta$. For $s = 1$, $e_{i_1} \in V(\bar{\gamma}_1, \bar{\mu}_1)_{\zeta_1}$ so nothing to prove. For $s \geq2$, using Corollary \ref{cor:quasi-asso} for $a = e_{i_1}$, $b=e_{i_2}* \cdots * e_{i_s}$ and $c=e_{j_1}* \cdots * e_{j_t}$, we have
\begin{align*}
&(e_{i_1} * \cdots * e_{i_s}) * (e_{j_1} * \cdots * e_{j_t})
=(a*b)*c\\
&=a*(b*c)+\sum_{j=0}^\infty\left(
a*_{-j-2}(b*_jc) + (-1)^{p(a)p(b)}b*_{-j-2}(a*_jc)
\right)\\
&+\chi_{a,b}\sum_{j=0}^\infty\left(
a*_{-j-1}(b*_jc) + (-1)^{p(a)p(b)}b*_{-j-1}(a*_jc)
\right).
\end{align*}
First of all, $a*_{-j-2}(b*_j c), b*_{-j-2}(a*_jc) \in (\gH{J})_{\zeta_{-}} = W_{\zeta_{-}} \subset W_\zeta$ for all $ j\geq 0$ by \eqref{eq:Vzeta2}, Lemma \ref{lem:J-trivial} and the inductive assumption on $\zeta$. Next, if $a \in \gH{V}$, we have $\epsilon_a = 0$ so that $\chi_{a, b} = 0$. Then $b * c \in (\gH{J})_{\zeta-\zeta(i_1)} = W_{\zeta-\zeta(i_1)}$ by the inductive assumption so $a*(b*c) \in W_\zeta$ by Lemma \ref{lem:W-ideal} (1). Furthermore, if $b \in \gH{V}$, we also have $\chi_{a, b} = 0$ and then $a*(b*c) \in W_\zeta$ by Lemma \ref{lem:W-ideal} (3). Finally, if $a \notin \gH{V}$ and $b \notin \gH{V}$, $a*(b*c) \in W_\zeta$ again by Lemma \ref{lem:W-ideal} (3), and both of $a*_{-j-1}(b*_jc)$ and $b*_{-j-1}(a*_jc)$ for $j \geq0$ belong to $(\gH{J})_{\zeta_{-}}=W_{\zeta_{-}}\subset W_\zeta$ by \eqref{eq:Vzeta2}, Lemma \ref{lem:J-trivial} and the inductive assumption on $\zeta$. Therefore we have $(a*b)*c \in W_\zeta$. This proves \eqref{eq:W=J_2}.

As a consequence, it follows from \eqref{eq:W=J} and \eqref{eq:W=J_2} that
\begin{align*}
(\gH{J})_\zeta 
&= \sum_{
\begin{subarray}{c}
\zeta_1, \zeta_2 \in \monoid_{\geq 0}\\
\zeta_1 + \zeta_2 \leq \zeta
\end{subarray}}
\Bigl(
(\gH{V})_{\zeta_1} \circ (\gH{V})_{\zeta_2}
+ \sum_{\begin{subarray}{c} \bar{\gamma}_1 + \bar{\gamma}_2 = \bar{\mu}_1 + \bar{\mu}_2 \\ \bar{\gamma}_1 \neq \bar{\mu}_1, \bar{\gamma}_2 \neq \bar{\mu}_2 \end{subarray}} V(\bar{\gamma}_1, \bar{\mu}_1)_{\zeta_1} * V(\bar{\gamma}_2, \bar{\mu}_2)_{\zeta_2}
\Bigr)\\
&\subset \sum_{
\begin{subarray}{c}
\zeta_1, \zeta_2 \in \monoid_{\geq 0}\\
\zeta_1 + \zeta_2\leq \zeta
\end{subarray}}
W_{\zeta_1+\zeta_2}
=W_\zeta.
\end{align*}
This completes the proof.
\end{proof}
\end{lemma}

\begin{lemma}\label{lem:W1subsetSigma}
$W_{1, \zeta} \subset \Span_\C \Sigma_\zeta$
\begin{proof}
The assertion follows thanks to the proof of \cite[Proposition 3.21]{DK} with \eqref{eq:Vzeta2}\eqref{eq:Vzeta3}, Lemma \ref{lem:Vzeta} and Lemma \ref{lem:W=J}.
\end{proof}
\end{lemma}

\begin{lemma}\label{lem:Sigma=W}
$\Span_\C \Sigma_\zeta = W_\zeta$.
\begin{proof}
It is clear that $\Span_\C \Sigma_\zeta \subset W_\zeta$. We will show that $W_\zeta \subset \Span_\C \Sigma_\zeta$. We have $W_{1, \zeta} \subset \Span_\C \Sigma_\zeta$ by Lemma \ref{lem:W1subsetSigma}. Thus, it suffices to show that $W_{2, \zeta} \subset \Span_\C \Sigma_\zeta$. We will prove that
\begin{align}\label{eq:Sigma=W}
e_{i_1} * \cdots * e_{i_s} \in W_{2, \zeta}
\end{align}
for $s>1$, $i_\ell \in \mathbb{I}$ with $\zeta(\overrightarrow{i}) = \zeta$ for $\overrightarrow{i} = (i_1, \ldots, i_s)$ such that $e_{i_j} \notin \gH{V}$ for some $j$. If $\#(\overrightarrow{i}) = 0$, $e_{i_1} * \cdots * e_{i_s} \in \Sigma_\zeta$ so nothing to prove. Suppose that $\#(\overrightarrow{i}) \geq 1$. If $i_p > i_q$ for some $1 \leq p <  q \leq s$, we may choose $p = \min\{ \ell \mid i_\ell > i_{\ell+1}\}$. Then $i_p > i_{p+1}$. Otherwise we have $i_\ell \leq i_{\ell+1}$ for all $\ell$. Then there exists $i_p = i_{p+1}$ such that $p(e_{i_p}) = \bar{1}$.
\smallskip

\hspace{-3.5mm}(1) Consider the case that $i_p > i_{p+1}$. Let $\bar{\gamma}_j, \bar{\mu}_j \in \R/\Z$ for $j=1,2,3$ such that $e_{i_p} \in V(\bar{\gamma}_1, \bar{\mu}_1)$, $e_{i_{p+1}} \in V(\bar{\gamma}_2, \bar{\mu}_{2})$ and $a_0 = e_{p+2} * \dots * e_{i_s} \in V(\bar{\gamma}_3, \bar{\mu}_3)$. Set $\zeta_1 = \zeta(i_p)$, $\zeta_2 = \zeta(i_{p+1})$ and $\zeta_3 = \zeta(i_{p+2}, \ldots, i_s)$. Denote by $P(a) = e_{i_1} * \cdots * e_{i_{p-1}} * a$ for $a \in V$.
\smallskip

\hspace{-2mm}(i) Suppose that $\chi_{e_{i_p}, e_{i_{p+1}}} = 0$. We have
\begin{align*}
&e_{i_1} * \cdots * e_{i_s} = P(a_1) + P(a_2 * a_0) + P(a_3 * a_0) + (-1)^{p(e_{i_p})p(e_{i_{p+1}})}P(a_4),\\
&a_1=e_{i_p} * e_{i_{p+1}} * a_0 - (e_{i_p} * e_{i_{p+1}}) * a_0,\\
&a_2 = e_{i_p} * e_{i_{p+1}}  - (-1)^{p(e_{i_p})p(e_{i_{p+1}})}e_{i_{p+1}} * e_{i_{p}} - [{e_{i_p}}_*e_{i_{p+1}}],\\
&a_3=[{e_{i_p}}_*e_{i_{p+1}}],\quad
a_4 = {e_{i_{p+1}}} * e_{i_{p}}*a_0.
\end{align*}
Since $\chi_{e_{i_p}, e_{i_{p+1}}} = 0$, it follows from Lemma \ref{lem:Vzeta} (3) that
\begin{align*}
a_1 \in \left(V(\bar{\gamma}_1, \bar{\mu}_1) \circ V(\bar{\gamma}_2+\bar{\gamma}_3, \bar{\mu}_2+\bar{\mu}_3) + V(\bar{\gamma}_2, \bar{\mu}_2) \circ V(\bar{\gamma}_1+\bar{\gamma}_3, \bar{\mu}_1+\bar{\mu}_3)\right)_{(\zeta_1+\zeta_2+\zeta_3)_{-}}.
\end{align*}
Then $a_1$ belongs to $(\gH{J})_{(\zeta_1+\zeta_2+\zeta_3)_{-}}$ or $V(\bar{\gamma}, \bar{\mu})_{(\zeta_1+\zeta_2+\zeta_3)_{-}}$ with some $\bar{\gamma} \neq \bar{\mu}$. Thus, $P(a_1) \in (\gH{J})_{\zeta_{-}}$ by Lemma \ref{lem:J-or-Vodd}. Since $(\gH{J})_{\zeta_{-}} = W_{\zeta_{-}} = \Span_\C\Sigma_{\zeta_{-}} \subset \Span_\C\Sigma_{\zeta}$ by Lemma \ref{lem:W=J} and the inductive assumption on $\zeta$, we have $P(a_1) \in \Span_\C\Sigma_{\zeta}$. By Lemma \ref{lem:Vzeta} (4),
\begin{align*}
a_2 \in \left(V(\bar{\gamma}_1+\bar{\gamma}_2, \bar{\mu}_1 + \bar{\mu}_2) \circ \vac \right)_{(\zeta_1+\zeta_2)_{-}}.
\end{align*}
Then $a_2$ belongs to $(\gH{J})_{(\zeta_1+\zeta_2)_{-}}$ or $V(\bar{\gamma}, \bar{\mu})_{(\zeta_1+\zeta_2)_{-}}$ with some $\bar{\gamma} \neq \bar{\mu}$. Similarly, we have $P(a_2 *a_0) \in \Span_\C\Sigma_{\zeta}$. By \eqref{eq:Vgamma-mu}\eqref{eq:Vzeta3}, $a_3 \in V(\bar{\gamma}_1+\bar{\gamma}_2, \bar{\mu}_1+\bar{\mu}_2)_{(\zeta_1+\zeta_2)_{-}}$. If $\bar{\gamma}_1+\bar{\gamma}_2 \neq \bar{\mu}_1+\bar{\mu}_2$, $P(a_3*a_0) \in (\gH{J})_{\zeta_{-}}$ by Lemma \ref{lem:J-or-Vodd}. Hence $P(a_3*a_0) \in \Span_\C\Sigma_{\zeta}$. Otherwise $\bar{\gamma}_1+\bar{\gamma}_2 = \bar{\mu}_1+\bar{\mu}_2$. Then if $\bar{\gamma}_1 \neq \bar{\mu}_1$, we have $\chi_{e_{i_p}, e_{i_{p+1}}} = 1$, contrary to our assumption. Thus $\bar{\gamma}_1=\bar{\mu}_1$ and $\bar{\gamma}_2 = \bar{\mu}_2$. Since $e_{i_1} * \cdots * e_{i_s} \in W_{2, \zeta}$, there exists some $j \neq p, p+1$ such that $e_{i_j} \notin \gH{V}$. Then
\begin{align*}
P(a_3 * a_0) = e_{i_1} * \cdots * [{e_{i_p}}_*e_{i_{p+1}}] * \cdots * e_{i_s} \in (\gH{J})_{\zeta_{-}}
\end{align*}
by Lemma \ref{lem:J-or-Vodd}. Hence $P(a_3 * a_0) \in \Span_\C\Sigma_{\zeta}$. Finally, $P(a_4) \in \Span_\C\Sigma_{\zeta}$ by the inductive assumption on $\#(\overrightarrow{i})$. Therefore $e_{i_1} * \cdots * e_{i_s} \in \Span_\C\Sigma_{\zeta}$.
\smallskip

\hspace{-2mm}(ii) Suppose that $\chi_{e_{i_p}, e_{i_{p+1}}} = 1$. We have
\begin{align*}
&e_{i_1} * \cdots * e_{i_s} = P(a_1) + P(a_2' * a_0) + (-1)^{p(e_{i_p})p(e_{i_{p+1}})}P(a_4),\\
&a_1=e_{i_p} * e_{i_{p+1}} * a_0 - (e_{i_p} * e_{i_{p+1}}) * a_0,\quad
a_4 = {e_{i_{p+1}}} * e_{i_{p}}*a_0,\\
&a_2' = e_{i_p} * e_{i_{p+1}}  - (-1)^{p(e_{i_p})p(e_{i_{p+1}})}e_{i_{p+1}} * e_{i_{p}}.
\end{align*}
By Lemma \ref{lem:Vzeta} (3), we have
\begin{align*}
a_1 \in \left(V(\bar{\gamma}_1, \bar{\mu}_1) * V(\bar{\gamma}_2+\bar{\gamma}_3, \bar{\mu}_2+\bar{\mu}_3) + V(\bar{\gamma}_2, \bar{\mu}_2) * V(\bar{\gamma}_1+\bar{\gamma}_3, \bar{\mu}_1+\bar{\mu}_3)\right)_{(\zeta_1+\zeta_2+\zeta_3)_{-}}.
\end{align*}
Notice that $\bar{\gamma}_1 \neq \bar{\mu}_1$ and $\bar{\gamma}_2 \neq \bar{\mu}_2$ due to our assumption that $\chi_{e_{i_p}, e_{i_{p+1}}} = 1$. Thus, $a_1$ belongs to $(\gH{J})_{(\zeta_1+\zeta_2+\zeta_3)_{-}}$ or $V(\bar{\gamma}, \bar{\mu})_{(\zeta_1+\zeta_2+\zeta_3)_{-}}$ with some $\bar{\gamma} \neq \bar{\mu}$ again. Hence $P(a_1) \in \Span_\C\Sigma_{\zeta}$. By Lemma \ref{lem:Vzeta} (4),
\begin{align*}
a_2' \in \left(V(\bar{\gamma}_1+\bar{\gamma}_2, \bar{\mu}_1 + \bar{\mu}_2) \circ \vac \right)_{(\zeta_1+\zeta_2)_{-}}.
\end{align*}
As in case of $a_2$, we also have $P(a_2'*a_0) \in \Span_\C\Sigma_{\zeta}$. In the same way, $P(a_4) \in \Span_\C\Sigma_{\zeta}$ by the inductive assumption on $\#(\overrightarrow{i})$. Therefore $e_{i_1} * \cdots * e_{i_s} \in \Span_\C\Sigma_{\zeta}$. This complete the proof of \eqref{eq:Sigma=W} in case that $i_p > i_{p+1}$.
\smallskip

\hspace{-3.5mm}(2) Consider the case that $i_p = i_{p+1}$ such that $p(e_{i_p}) = \bar{1}$. We have
\begin{align*}
&e_{i_1} * \cdots * e_{i_s} =P(b_1) + \frac{1}{2}P(b_2*a_0) + \frac{1}{2}(1-\chi_{e_{i_p}, e_{i_{p}}})P(b_3*a_0),\\
&b_1=e_{i_p} * e_{i_{p}} * a_0 - (e_{i_p} * e_{i_{p}})*a_0,\\
&b_2=e_{i_p} * e_{i_{p}} + e_{i_p} * e_{i_{p}} -(1-\chi_{e_{i_p}, e_{i_{p}}})[{e_{i_p}}_*e_{i_{p}}],\\
&b_3=[{e_{i_p}}_*e_{i_{p}}].
\end{align*}
Then we can prove that $e_{i_1} * \cdots * e_{i_s} \in \Span_\C\Sigma_{\zeta}$ similarly to the case $i_p > i_{p+1}$. This completes the proof.
\end{proof}
\end{lemma}

\begin{theorem}\label{thm:J-PBW}
$\Span_\C \Sigma = W = \gH{J}$. In particular, $\Sigma$ is a basis of $\gH{J}$.
\begin{proof}
The first assertion is immediate from Lemma \ref{lem:W=J} and Lemma \ref{lem:Sigma=W}. The second assertion follows from the fact that elements of $\Sigma$ are linearly independent.
\end{proof}
\end{theorem}

Recall that $\conf = \C[\partial]\otimes\base$ for $\base = \Span_\C\{e_\alpha\}_{\alpha \in I}$ is a non-linear Lie conformal superalgebra such that $V \simeq V(\conf)$.

\begin{corollary}\label{cor:ZhuR}
Let $\gH{\conf} = \conf \cap \gH{V}$. Then $\gH{\conf} \cap \gH{J} = (\partial + H)\gH{\conf}$.
\begin{proof}
$\gH{\conf} \cap \gH{J} \supset (\partial + H)\gH{\conf}$ is clear. On the other hand, $\{\partial^{(k)}e_\alpha - \binom{-\Delta_{\alpha}}{k}e_\alpha \mid k \geq1, \alpha \in I, \epsilon_{\alpha}=0\}$ is a basis of $\gH{\conf} \cap \gH{J}$ by Theorem \ref{thm:J-PBW}. Since
\begin{align*}
\partial^{(k)}e_\alpha - \binom{-\Delta_{\alpha}}{k}e_\alpha = \frac{1}{k}(\partial + H)\partial^{(k-1)}e_\alpha - \frac{\Delta_{\alpha}+k-1}{k}\left(\partial^{(k-1)}e_\alpha-\binom{-\Delta_{\alpha}}{k-1}e_\alpha\right),
\end{align*}
we have $\partial^{(k)}e_\alpha - \binom{-\Delta_{\alpha}}{k}e_\alpha \in (\partial + H)\gH{\conf}$ by induction on $k$. Therefore $\gH{\conf} \cap \gH{J} \subset (\partial + H)\gH{\conf}$. This completes the proof.
\end{proof}
\end{corollary}

Let
\begin{align*}
\projZhu \colon \gH{V} \twoheadrightarrow \gH{\Zhu} V
\end{align*}
be the canonical projection.

\begin{corollary}\label{cor:ZhuV-PBW}
$\gH{\Zhu} V$ is PBW generated by $\{\projZhu(e_\alpha) \mid \alpha \in I, \epsilon_{\alpha}=0\}$ and pregraded by $\monoid_{\geq0}$.
\begin{proof}
The first statement is clear from Theorem \ref{thm:J-PBW}. By Proposition \ref{prop:J-ideal} (5) and \eqref{eq:Vzeta3},
\begin{align}\label{eq:ZhuR}
a * b - (-1)^{p(a)p(b)}b*a \equiv [a_*b] \in V_{(\zeta_1+\zeta_2)_{-}}
\mod. \gH{J}
\end{align}
for $a \in \base_{\zeta_1} \cap \gH{V}$ and $b \in \base_{\zeta_2} \cap \gH{V}$ with $\zeta_1, \zeta_2 \in \monoid_{>0}$. This implies the second statement.
\end{proof}
\end{corollary}

Recall that $\twist = \mathrm{e}^{2\pi i H}$, and that $\epsilon_\alpha = 0 \iff (g \circ \theta_H^{-1})(e_\alpha) = e_\alpha$. Set
\begin{align*}
\gH{\Zhu} \conf := \projZhu(\gH{\conf}),\quad
\gH{\base} := \base\cap \gH{V} = \Span_\C\{e_\alpha \mid \alpha \in I,\ (g \circ \theta_H^{-1})(e_\alpha) = e_\alpha\}.
\end{align*}
By Corollary \ref{cor:ZhuR},
\begin{align*}
\gH{\Zhu} \conf = \gH{\conf}/(\partial + H)\gH{\conf} = \projZhu(\gH{\base}) \simeq \gH{\base}.
\end{align*}
Let $\rho \colon \gH{\Zhu} V \hookrightarrow T(\gH{\Zhu} \conf)$ be a linear embedding by $e_{\alpha_1} \cdots e_{\alpha_s} \mapsto e_{\alpha_1} \otimes \cdots \otimes e_{\alpha_s}$ for $\alpha_i \in I$ such that $\alpha_i \leq \alpha_{i+1}$ and $\alpha_i < \alpha_{i+1}$ if $p(e_{\alpha_i}) = \bar{1}$. Since $\gH{\Zhu} \conf$ has a $\monoid_{>0}$-grading induced from that of $\gH{\base}$,
\begin{align*}
\gH{\Zhu} \conf \otimes \gH{\Zhu} \conf \rightarrow T(\gH{\Zhu} \conf),\quad
\projZhu(a) \otimes \projZhu(b) \mapsto (\rho\circ\projZhu)([a_*b]),\quad
a, b \in \gH{\base}
\end{align*}
defines a structure of non-linear Lie superalgebra on $\gH{\Zhu} \conf$ such that
\begin{align*}
\gH{\Zhu} V \simeq U(\gH{\Zhu} \conf)
\end{align*}
by Proposition \ref{prop:DK-U(r)} and Corollary \ref{cor:ZhuV-PBW}. Therefore, we obtain the following theorem:

\begin{theorem}\label{thm:main PBW thm}
Let $V$ be a vertex superalgebra with a Hamiltonian operator $H$ whose eigenvalues belong to $\R$ and $g \in \operatorname{Aut}_H V$ be a diagonalizable automorphism with modulus $1$ eigenvalues. Suppose that $V$ is freely generated by free generators $\{e_\alpha\}_{\alpha \in I}$ homogeneous for $g$, $H$ and the parity, with an ordered index set $I$ and is pregraded by $\monoid_{\geq0}$ in the sense of Section \ref{sec:def-va}. Let $\projZhu \colon \gH{V} \twoheadrightarrow \gH{\Zhu} V$ be the canonical projection, $\gH{\conf} = \C[\partial]\otimes\gH{\base} \subset \gH{V}$ for $\gH{\base} = \Span_\C\{e_\alpha \mid \alpha \in I,\ (g \circ \theta_H^{-1})(e_\alpha) = e_\alpha\}$. Then $\gH{\Zhu} \conf := \projZhu(\gH{\conf})$ has a structure of non-linear Lie superalgebra such that $\gH{\Zhu} V \simeq U(\gH{\Zhu} \conf)$. In particular, $\gH{\Zhu} V$ is an associative superalgebra generated by $\projZhu(\gH{\base}) \simeq \gH{\base}$ with the defining relations
\begin{align}\label{eq:defininig_rel_Zhu}
[\projZhu(a), \projZhu(b)] = \sum_{j=0}^{\infty}\binom{\Delta_a-1}{j}\projZhu(a_{(j)}b),\quad
a, b \in \gH{\base}.
\end{align}
\end{theorem}

\section{Associated varieties}\label{sec:asso var}
Let $V$ be a vertex superalgebra with a Hamiltonian operator $H$ whose eigenvalues belong to $\R$, and $g \in \operatorname{Aut}_H V$ be a diagonalizable automorphism with modulus $1$ eigenvalues.

\subsection{Poisson superalgebras from Zhu}
Recall that $V$ has a increasing filtration $F^pV = \Span_\C\{a \in V \mid \Delta_a \leq p\}$ defined by \cite{Li}. Set
\begin{align*}
F^p \gH{V} := F^pV \cap \gH{V} = \bigoplus_{\Delta \leq p}V(\bar{\Delta}, \Delta),\quad
p \in \R.
\end{align*}
Since, for all $a \in F^p\gH{V}$ and $b\in F^q \gH{V}$,
\begin{align*}
a * b \in F^{p+q}\gH{V},\quad
a*b - (-1)^{p(a)p(b)}b*a \in F^{p+q-1}\gH{V},
\end{align*}
$F^p \gH{V}$ induces an increasing filtration on $\gH{\Zhu} V$ by
\begin{align*}
F^p\gH{\Zhu} V := F^p \gH{V}/(F^p\gH{V} \cap \gH{J}) \subset \gH{\Zhu}V,
\end{align*}
and makes $\gH{\Zhu}V$ an almost supercommutative superalgebra. That is, as proved in \cite[Proposition 2.7(b)]{DK}, the associated graded superspace
\begin{align*}
\mathrm{gr}\gH{\Zhu} V := \bigoplus_{p \in \R}\mathrm{gr}^p\gH{\Zhu}V,\quad
\mathrm{gr}^p\gH{\Zhu}V = \frac{F^p\gH{\Zhu} V}{F^{p_-}\gH{\Zhu} V}
\end{align*}
becomes a Poisson superalgebras with the associative supercommutative product
\begin{align*}
\mathrm{gr}^p\gH{\Zhu}V \times \mathrm{gr}^q\gH{\Zhu}V \rightarrow \mathrm{gr}^{p+q}\gH{\Zhu}V,\quad
(\overline{a}, \overline{b}) \mapsto \overline{a*b},
\end{align*}
and the Lie superbracket
\begin{align*}
\mathrm{gr}^p\gH{\Zhu}V \times \mathrm{gr}^q\gH{\Zhu}V \rightarrow \mathrm{gr}^{p+q-1}\gH{\Zhu}V,\quad
(\overline{a}, \overline{b}) \mapsto \overline{[a_*b]},
\end{align*}
where
\begin{align*}
F^{p_-}\gH{\Zhu} V := F^{p_-} \gH{V}/(F^{p_-}\gH{V} \cap \gH{J}),\quad
F^{p_-} \gH{V} := \bigoplus_{\Delta < p}V(\bar{\Delta}, \Delta).
\end{align*}

On the other hand, let
\begin{align*}
\gH{C}(V) := \Span_\C\{a_{(-2+\chi_{a, b})}b \in V \mid a, b \in V\} \cap \gH{V},
\end{align*}
where $\chi_{a, b}$ is defined by the equation \eqref{eq:def-chi}. Then, by \cite[Proposition 2.7(a)]{DK},
\begin{align*}
\gH{R}(V) := \gH{V}/\gH{C}(V)
\end{align*}
is also a Poisson superalgebra with the associative commutative product
\begin{align*}
\overline{a} \cdot \overline{b} := \overline{a_{(-1)}b}
\end{align*}
and the Lie superbracket
\begin{align*}
\{\overline{a}, \overline{b}\} = \overline{a_{(0)}b}
\end{align*}
for $\overline{a}, \overline{b} \in \gH{R}(V)$. Moreover, $H$ defines an algebra grading on $\gH{R}(V)$:
\begin{align*}
\gH{R}(V) = \bigoplus_{p \in \R}\gH{R}^p(V).
\end{align*}
Notice that
\begin{align*}
\mathrm{gr}^p\gH{\Zhu}V = \frac{F^p \gH{V}/(F^p\gH{V} \cap \gH{J})}{F^{p_-} \gH{V}/(F^{p_-}\gH{V} \cap \gH{J})}
\simeq \frac{F^p\gH{V}}{F^{p_-}\gH{V}+F^p\gH{V}\cap\gH{J}},
\end{align*}
and that
\begin{align*}
\gH{R}^p(V) \simeq \frac{F^p\gH{V}}{F^{p_-}\gH{V}+F^p\gH{V} \cap \gH{C}(V)}.
\end{align*}
Since
\begin{align*}
&F^p\gH{V} \cap \gH{C}(V)\\
&=
\sum_{\Delta_1+\Delta_2\leq p-1}
V(\bar{\Delta}_1, \Delta_1)_{(-2)}V(\bar{\Delta}_2, \Delta_2)
+
\sum_{
\begin{subarray}{c}
\Delta_1+\Delta_2 \leq p \\
\Delta_1+\Delta_2 \in \bar{\gamma}_1+\bar{\gamma}_2 \\
\Delta_1 \notin \bar{\gamma}_1, \Delta_2 \notin \bar{\gamma}_2
\end{subarray}
}
V(\bar{\gamma}_1, \Delta_1)_{(-1)}V(\bar{\gamma}_2, \Delta_2)\\
&\subset F^{p_-}\gH{V}+F^p\gH{V}\cap\gH{J}
\end{align*}
we have \cite[Proposition 2.17(c)]{DK} a canonical surjective morphism of graded Poisson superalgebras
\begin{align}\label{eq:R_surj}
\gH{R}(V) \twoheadrightarrow \mathrm{gr}\gH{\Zhu} V.
\end{align}
\begin{theorem}\label{thm:R_surj}
$\gH{R}(V) \simeq \mathrm{gr}\gH{\Zhu} V$ under the assumptions in Theorem \ref{thm:main PBW thm}.
\begin{proof}
By Theorem \ref{thm:J-PBW}, we have
\begin{align*}
&F^p\gH{V} \cap \gH{J}\\
&\subset \sum_{\Delta_1 + \Delta_2  \leq p-1}
V(\bar{\Delta}_1, \Delta_1) \circ V(\bar{\Delta}_2, \Delta_2)
+ \sum_{
\begin{subarray}{c}
\Delta_1+\Delta_2 \leq p \\
\Delta_1+\Delta_2 \in \bar{\gamma}_1+\bar{\gamma}_2 \\
\Delta_1 \notin \bar{\gamma}_1, \Delta_2 \notin \bar{\gamma}_2
\end{subarray}
}
V(\bar{\gamma}_1, \Delta_1) * V(\bar{\gamma}_2, \Delta_2)\\
&\subset F^{p_-}\gH{V} + F^p\gH{V} \cap \gH{C}(V).
\end{align*}
Thus, the assertion follows from \eqref{eq:R_surj}.
\end{proof}
\end{theorem}

\subsection{Lisse conditions}
Let $\gH{X}(V)$ be a Poisson variety defined by
\begin{align*}
\gH{X}(V) := \operatorname{Specm}\gH{R}(V).
\end{align*}
In case $g=\twist$,
\begin{align*}
X_V := X_{\twist, H}(V),
\end{align*}
called the associated variety of $V$ by Arakawa in \cite{Arakawa12}.
\smallskip

Suppose that $V$ satisfies the same assumptions as in Theorem \ref{thm:main PBW thm}. Then $X_V$ is an affine variety isomorphic to $\Span_\C\{e_\alpha \mid \alpha \in I,\,p(e_\alpha) = \bar{0}\}$. Thus, $X_V$ has a $\C^*$-action by $g \circ \twist^{-1}$ so that
\begin{align*}
X_V = \bigsqcup_{\bar{\gamma} \in \R/\Z}(X_V)_{\bar{\gamma}},\quad
(X_V)_{\bar{\gamma}} = \{ x \in X_V \mid (g \circ \twist)(x) = \mathrm{e}^{2\pi i \gamma}x\}.
\end{align*}
Since $\gH{\Zhu}V$ is PBW generated by $\{e_\alpha\mid\alpha \in I,\,(g \circ \twist^{-1})(e_\alpha)=e_\alpha\}$, $\gH{X}(V)$ is an affine variety isomorphic to $(X_V)_{\bar{0}}$ by Theorem \ref{thm:R_surj}. Therefore,
\begin{align*}
X_V = \gH{X}(V) \sqcup (X_V)_{\neq\bar{0}},\quad
(X_V)_{\neq\bar{0}} = \bigsqcup_{\bar{\gamma} \neq \bar{0}}(X_V)_{\bar{\gamma}}.
\end{align*}
Let $U$ be a quotient vertex superalgebra of $V$ by a $(g, H)$-invariant ideal. Then the projection $V \twoheadrightarrow U$ induces a projection $\gH{R}(V) \twoheadrightarrow \gH{R}(U)$. In particular,
\begin{align*}
X_U \hookrightarrow X_V.
\end{align*}
Thus, $X_U$ is a subvariety of $V$. Following \cite{BFM}, we call $X_U$ \textit{lisse along a subspace $Y \subset X_V$} if
\begin{align*}
\dim(X_U \cap Y^\perp) = 0,\quad
Y^\perp = \{ x \in X_V \mid x \notin Y\}.
\end{align*}
\begin{proposition}\label{prop:lisse along}
$\gH{X}(U) \subset X_U \cap \gH{X}(V)$. Therefore, $\dim\gH{X}(U) = 0$ if $X_U$ is lisse along $(X_V)_{\neq \bar{0}}$.
\begin{proof}
Notice that $R_{\twist, H}(T) = T/C_2(T)$ for $T = U, V$, where
\begin{align*}
C_2(T) = \Span_\C \{ a_{(-2)}b \mid a, b \in T\}.
\end{align*}
Recall that the projection $V \twoheadrightarrow U$ induces a projection $\gH{R}(V) \twoheadrightarrow \gH{R}(U)$. Using the fact that $C_2(T)\cap\gH{T} \subset \gH{C}(T)$, it follows that there exist the projections $R_{\twist, H}(T) \twoheadrightarrow \gH{R}(T)$ for $T=U, V$.  Thus, we have the following commutative diagram:
\begin{equation*}
\xymatrix@W10pt@H12pt@R20pt@C20pt{
R_{\twist, H}(V) \ar@{->>}[r] \ar@{->>}[d]& R_{\twist, H}(U) \ar@{->>}[d] \\
\gH{R}(V) \ar@{->>}[r] & \gH{R}(U) \makebox[0pt][l]{\,.} }
\end{equation*}
Finally, the functor $\operatorname{Specm}(?)$ applies to the diagram above.
\end{proof}
\end{proposition}

\section{Cohomology on non-linear superalgebras}\label{sec:coh of non-linear}
Let $\monoid$ be a discrete additive subgroup of $\R$. Denote by $\monoid_{\geq0} = \monoid \cap \R_{\geq0}$, by $\monoid_{>0} = \monoid \cap \R_{>0}$ and by $\zeta_-$ the largest number in $\monoid_{\geq0}$ that is strictly smaller than $\zeta$ for each $\zeta \in \monoid_{>0}$ as in Section \ref{sec:non-linear salg}. Set
\begin{align*}
\Charge = \left\{(p, q) \in \frac{1}{N}\Z \times \frac{1}{N}\Z \,\middle|\, p+q \in \Z_{\geq0}\right\},\quad
N \in \N.
\end{align*}

\subsection{Cohomology on non-linear Lie superalgebras}

Let $\base$ be a $\monoid_{>0} \times \Charge$-graded finite-dimensional superspace
\begin{align}\label{eq:r-gr}
\base = \bigoplus_{\zeta \in \monoid_{>0}, (p, q)\in \Charge}\base^{p, q}_\zeta.
\end{align}
Denote by $\base^{p, q} = \bigoplus_{\zeta \in \monoid_{>0}}\base^{p, q}_\zeta$ and by $\base_\zeta = \bigoplus_{(p, q) \in \Charge}\base^{p, q}_\zeta$. Then the tensor superalgebra $T(\base)$ is also $\Charge$-graded and has an increasing $K_{\geq0}$-filtration by
\begin{align*}
T(\base)_\zeta = \bigoplus_{(p, q)\in \Charge}T(\base)^{p, q}_\zeta,\quad
T(\base)^{p, q}_\zeta = \sum_{\begin{subarray}{c} p_1+\cdots+p_s = p\\
q_1+\cdots+q_s = q\\ \zeta_1+\cdots+\zeta_s \leq \zeta\end{subarray}}\base^{p_1, q_1}_{\zeta_1}\otimes \cdots \otimes\base^{p_s, q_s}_{\zeta_s},
\end{align*}
where $T(\base)^{p, q}_0 = \delta_{p,0}\delta_{q,0}\C1$. Define a $\Z_{\geq0}$-grading
\begin{align*}
\base = \bigoplus_{n =0}^\infty \base^n,\quad
\base^n = \bigoplus_{p+q=n}\base^{p,q}.
\end{align*}
Suppose that $\base$ is a non-linear Lie superalgebra with respect to the $\monoid_{>0}$-grading such that
\begin{align}\label{eq:r-cond1}
[\base^{p_1, q_1}_{\zeta_1}, \base^{p_2, q_2}_{\zeta_2}] \subset \bigoplus_{l=0}^\infty T(\base)^{p_1+p_2+l, q_1+q_2-l}_{(\zeta_1+\zeta_2)_-}
\end{align}
Let $U(\base)$ be the universal enveloping algebra of $\base$. Set
\begin{align*}
F^pU(\base)^n_\zeta
= \Span_\C\left\{a_1 \cdots a_s \in U(\base)\ \middle|\ 
\begin{lgathered}
a_k \in \base^{p_k, q_k}_{\zeta_k}, \sum_k(p_k+q_k) = n,\\
\sum_k p_k \geq p,\ \sum_k\zeta_k \leq \zeta
\end{lgathered}
\ \right\}
\end{align*}
for $p \in \frac{1}{N}\Z, n \in \Z_{\geq0}, \zeta \in \monoid_{\geq0}$. Denote by $F^pU(\base)^n = \sum_{\zeta}F^pU(\base)^n_\zeta$ and by $F^pU(\base) = \bigoplus_{n=0}^\infty F^pU(\base)^n$. Then $F^pU(\base)$ defines a decreasing $\frac{1}{N}\Z$-filtration on $U(\base)$. Let
\begin{align*}
&\gr U(\base) = \bigoplus_{(p, q) \in \Charge}\gr U(\base)^{p, q},\quad
\gr U(\base)^{p, q} = F^pU(\base)^{p+q}/F^{p+\frac{1}{N}}U(\base)^{p+q},\\
&\gr U(\base)_\zeta = \bigoplus_{(p, q) \in \Charge}\gr U(\base)^{p, q}_\zeta,\quad
\gr U(\base)^{p, q}_\zeta = F^pU(\base)^{p+q}_\zeta/F^{p+\frac{1}{N}}U(\base)^{p+q}_\zeta.
\end{align*}
Then $\gr U(\base)_\zeta$ defines an increasing $\monoid_{\geq0}$-filtration on the associated graded superalgebra $\gr U(\base)$. Let $F^p\base^n_\zeta = \bigoplus_{p' \geq p}\base^{p', n-p'}_\zeta$ and $\gr\base = \bigoplus_{(p,q) \in \Charge, \zeta \in \monoid_{>0}}\gr^{p, q}\base_\zeta$ with $\gr^{p, q}\base_\zeta = F^p\base^{p+q}_\zeta/F^{p+\frac{1}{N}}\base^{p+q}_\zeta$. Then $\gr\base$ is isomorphic to $\base$ as a superspace and is naturally equipped with a structure of a non-linear Lie superalgebra induced from $\base$. Let $U(\gr \base)$ be the universal enveloping algebra of $\gr\base$. Set
\begin{align*}
U(\gr\base)^{p,q}_\zeta
= \Span_\C\left\{\bar{a}_1 \cdots \bar{a}_s \in U(\gr\base)\ \middle|\ 
\begin{lgathered}
\bar{a}_k \in \gr^{p_k, q_k}\base_{\zeta_k}, \sum_k p_k = p,\\
\sum_k q_k = q,\ \sum_k\zeta_k \leq \zeta
\end{lgathered}
\ \right\}.
\end{align*}
By \cite[Lemma 4.8]{DK}, there exists a canonical isomorphism of superalgebras
\begin{align*}
U(\gr\base) \simeq \gr U(\base),
\end{align*}
which induces $U(\gr\base)^{p, q}_\zeta \simeq \gr U(\base)^{p, q}_\zeta$.

An odd linear map $d \colon U(\base) \rightarrow U(\base)$ is called an odd differential on $U(\base)$ if $d^2=0$ and $d(a\cdot b) = d(a)\cdot b + (-1)^{p(a)}a\cdot d(b)$ for $a, b \in U(\base)$. Denote by $H(U(\base), d) = \Ker d/ \Img d$. Then $H(U(\base), d)$ has a structure of an associative $\C$-superalgebra induced from $U(\base)$. Suppose that
\begin{align}\label{eq:r-cond2}
d\cdot F^pU(\base)^n_\zeta \subset F^pU(\base)^{n+1}_\zeta.
\end{align}
Then $d$ maps $U(\base)^n$ to $U(\base)^{n+1}$ and thus $(U(\base), d)$ forms a cochain complex. Let $H^n(U(\base), d)$ be the $n$-th cohomology of the cochain complex $(U(\base), d)$. Then
\begin{align*}
H(U(\base), d) = \bigoplus_{n \in \Z}H^n(U(\base), d).
\end{align*}
An odd differential $d$ on $U(\base)$ is called \textit{almost linear} if
\begin{align}\label{eq:r-cond3}
d(\base^{p,q}_{\zeta}) \subset \base^{p,q+1}_{\zeta} \oplus F^{p+\frac{1}{N}}U(\base)^{p+q+1}_\zeta.
\end{align}
Denote by $[a]$ the element in $\gr U(\base)$ corresponding to $a \in U(\base)$. Define $\gr d \colon \gr U(\base) \rightarrow \gr U(\base)$ by $\gr d \cdot [a] = [d a]$. Then the condition \eqref{eq:r-cond3} is equivalent to say that
\begin{align*}
\gr d \cdot \gr^{p,q}\base_\zeta \subset \gr^{p, q+1}\base_\zeta.
\end{align*}
Thus $(\gr^{p, \bullet} \base, \gr d)$ froms a complex. Let
\begin{align*}
H^{p, q}(\gr \base, \gr d) = \frac{\Ker(\gr d \colon \gr^{p, q}\base \rightarrow \gr^{p, q+1}\base)}{\Img(\gr d \colon \gr^{p, q-1}\base \rightarrow \gr^{p, q}\base)}
\end{align*}
be the cohomology of the complex $(\gr^{p, \bullet} \base, \gr d)$. Then $d$ is called good if
\begin{align}\label{eq:r-cond4}
\base^{p,q}_{\zeta}\cap d^{-1}(F^{p+\frac{1}{N}}U(\base)^{p+q+1}_\zeta) \subset d(\base^{p, q-1}_{\zeta}) + F^{p+\frac{1}{N}}U(\base)^{p+q}_\zeta,\quad
p+q \neq0.
\end{align}
The condition \eqref{eq:r-cond4} is equivalent to say that
\begin{align*}
H^{p, q}(\gr\base, \gr d) = 0,\quad
p+q \neq0.
\end{align*}

\begin{theorem}[{\cite[Lemma 4.13, Theorem 4.14]{DK}}]\label{thm:r-coh}
Let $\base$ be a $\monoid_{>0} \times \Charge$-graded finite-dimensional non-linear Lie superalgebra satisfying \eqref{eq:r-cond1} and $d$ be an odd differential on $U(\base)$ satisfying \eqref{eq:r-cond2}\eqref{eq:r-cond3}\eqref{eq:r-cond4}. Then
\begin{enumerate}[label=(\arabic*)]
\item $H^n(U(\base), d) = 0$ for $n \neq 0$.
\item Let $\{e_\alpha\}_{\alpha \in \Cohindex(p, \zeta)}$ be a parity-homogeneous basis of $\base^{p, -p}_{\zeta} \cap d^{-1}(F^{p+\frac{1}{N}}U(\base)^1_\zeta)$ for $p \in \frac{1}{N}\Z, \zeta \in \monoid_{>0}$ with some index set $\Cohindex(p, \zeta)$. Then there exist elements $\{E_\alpha\}_{\alpha \in \Cohindex(p, \zeta)}$ in $F^pU(\base)^0_\zeta\cap\Ker d$ such that
\begin{align*}
E_\alpha -e_\alpha \in F^{p+\frac{1}{N}}U(\base)^0_\zeta,\quad
p(E_\alpha) = p(e_\alpha).
\end{align*}
\item Let $\pi_0 \colon \Ker d \twoheadrightarrow H(U(\base), d)$ be the canonical projection and $\Cohindex = \bigcup_{p, \zeta}\Cohindex(p, \zeta)$. Then $H^0(U(\base), d)$ is PBW generated by $\{\pi_0(E_\alpha)\}_{\alpha \in \Cohindex}$.
\item Let $H(\base) = \Span_\C\{\pi_0(E_\alpha)\mid \alpha \in \Cohindex\}$. Then $H(\base)$ has a structure of a non-linear Lie superalgebra such that $H^0(U(\base), d) \simeq U(H(\base))$.
\end{enumerate}
\end{theorem}

\subsection{Cohomology on non-linear Lie conformal superalgebras}\label{sec:cohomology on conf}

Let $\conf$ be a non-linear Lie conformal superalgebra freely generated by a finite-dimensional subspace $\base$ as a $\C[\partial]$-module. Suppose that $\base$ is $\monoid_{>0} \times \Charge \times \R$-graded as follows:
\begin{align*}
\base = \bigoplus_{
\begin{subarray}{c}
\zeta \in \monoid_{>0},
\Delta \in \R\\
(p, q)\in \Charge
\end{subarray}
}
\base^{p, q}_\zeta(\Delta).
\end{align*}
Define a $\monoid_{>0} \times \Charge \times \R$-grading on $\conf$ by
\begin{align}\label{eq:R-gr}
\conf = \bigoplus_{
\begin{subarray}{c}
\zeta \in \monoid_{>0},
\Delta \in \R,\\
(p, q) \in \Charge
\end{subarray}
}
\conf^{p, q}_{\zeta}(\Delta),\quad
\conf^{p, q}_{\zeta}(\Delta) = \bigoplus_{n=0}^\infty\partial^n\,\base^{p, q}_\zeta(\Delta-n).
\end{align}
Then tensor superalgebra $T(\conf)$ is $\Charge \times \R$-graded and has an increasing $\monoid_{\geq0}$-filtration by
\begin{align*}
&T(\conf)_\zeta = \bigoplus_{(p, q)\in \Charge, \Delta \in \R}T(\conf)^{p, q}_\zeta(\Delta),\\
&T(\conf)^{p, q}_\zeta(\Delta)
= \sum_{\begin{subarray}{c}
p_1 + \cdots + p_s = p, \\
q_1 + \cdots + q_s = q
\end{subarray}
}
\sum_{
\begin{subarray}{c}
\Delta_1+\cdots+\Delta_s = \Delta, \\
\zeta_1 + \cdots + \zeta_s \leq \zeta
\end{subarray}
}
\conf^{p_1, q_1}_{\zeta_1}(\Delta_1) \otimes \cdots \otimes \conf^{p_s, q_s}_{\zeta_s}(\Delta_s).
\end{align*}
Suppose that for all $n \in \Z_{\geq0}$,
\begin{align}\label{eq:R-cond1}
{\conf^{p_1, q_1}_{\zeta_1}(\Delta_1)}_{(n)} \conf^{p_2, q_2}_{\zeta_2}(\Delta_2) \subset \bigoplus_{l=0}^\infty T(\conf)^{p_1+p_2+l, q_1+q_2-l}_{\zeta_{-}}(\Delta_1+\Delta_2-n-1).
\end{align}
Define a decreasing $\frac{1}{N}\Z$-filtration on $\conf$ by
\begin{align*}
F^p\conf = \bigoplus_{n=0}^\infty\bigoplus_{\zeta \in \monoid_{>0}}F^p\conf^n_\zeta,\quad
F^p\conf^n_\zeta = \bigoplus_{p' \geq p, \Delta \in \R}\conf^{p', n-p'}_\zeta(\Delta).
\end{align*}
Then, the associated graded superspace
\begin{align*}
\gr\conf = \bigoplus_{(p,q) \in \Charge, \zeta \in \monoid_{>0}}\gr\conf^{p,q}_\zeta,\quad
\gr\conf^{p,q}_\zeta = F^p\conf^{p+q}_\zeta/F^{p+\frac{1}{N}}\conf^{p+q}_\zeta
\end{align*}
has a structure of a non-linear Lie conformal superalgebra induced from $\conf$. Let $V(\conf)$ be the universal enveloping vertex algebra of $\conf$. Set
\begin{align*}
F^pV(\conf)^n_\zeta(\Delta)
= \Span_\C\left\{\NO{a_1 \cdots a_s}\ \in V(\conf)\ \middle|\ 
\begin{gathered}
 a_k \in \conf^{p_k, q_k}_{\zeta_k}(\Delta_k), \sum_k(p_k+q_k) = n,\\
\sum_k p_k \geq p, \sum_k \Delta_k = \Delta, \sum_k\zeta_k \leq \zeta
\end{gathered}
\right\}.
\end{align*}
Then for $n \in \Z_{\geq0}$,
\begin{align*}
F^{p_1}V(\conf)^{n_1}_{\zeta_1}(\Delta_1)_{(n)}F^{p_2}V(\conf)^{n_2}_{\zeta_2}(\Delta_2) \subset F^{p_1+p_2}V(\conf)^{n_1+n_2}_{(\zeta_1+\zeta_2)_-}(\Delta_1+\Delta_2-n-1).
\end{align*}
Denote the associated graded superspace by
\begin{align*}
&\gr V(\conf) = \sum_{\zeta \in \monoid_{\geq0}}\bigoplus_{(p, q)\in \Charge, \Delta\in\R}\gr V(\conf)^{p,q}_\zeta(\Delta),\\
&\gr V(\conf)^{p,q}_\zeta(\Delta) = F^p V(\conf)^{p+q}_\zeta(\Delta)/F^{p+\frac{1}{N}}V(\conf)^{p+q}_\zeta(\Delta).
\end{align*}
Then it is easy to see that $\gr V(\conf)$ is $\Charge \times \R$-graded and has an increasing $\monoid_{\geq0}$-filtration. Let $V(\gr \conf)$ be the universal enveloping vertex algebra of $\gr\conf$. Similarly, $V(\gr \conf)$ is also $\Charge \times \R$-graded and has an increasing $\monoid_{\geq0}$-filtration. Using the same proof as \cite[Lemma 4.8]{DK}, it follows that there exists a canonical isomorphism $\gr V(\conf) \simeq V(\gr \conf)$, which induces $\gr V(\conf)^{p, q}_\zeta(\Delta) \simeq V(\gr \conf)^{p, q}_\zeta(\Delta)$.

An odd linear map $d \colon V(\conf) \rightarrow V(\conf)$ is called an odd differential on $V(\conf)$ if $d^2=0$, $d(\partial a) = \partial(da)$, $d(a_{(n)}b) = (da)_{(n)}b + (-1)^{p(a)}a_{(n)}(db)$ for $a, b \in V(\conf)$. In addition, we suppose that
\begin{align}\label{eq:R-cond2}
d \cdot F^pV(\conf)^n_\zeta(\Delta) \subset F^pV(\conf)^{n+1}_\zeta(\Delta).
\end{align}
Then $d \colon V(\conf)^n \rightarrow V(\conf)^{n+1}$ defines a cochain complex. Denote by $H^n(V(\conf), d)$ the $n$-th cohomology of $(V(\conf)^\bullet, d)$. Then $H(V(\conf), d) = \bigoplus_{n \in \Z}H^n(V(\conf), d)$ is a vertex superalgebra inherited from $V(\conf)$. Denote by $[a]$ the element in $\gr V(\conf)$ corresponding to $a \in V(\conf)$. Then $\gr d \colon \gr V(\conf)^{p,q} \rightarrow \gr V(\conf)^{p, q+1}$ is defined by $\gr [a] = [d(a)]$. Set
\begin{align*}
H^{p,q}(\gr V(\conf), \gr d) = \frac{\gr V(\conf)^{p,q} \cap \Ker \gr d}{\gr V(\conf)^{p,q} \cap \Img \gr d}.
\end{align*}
An odd differential $d$ on $V(\conf)$ is called \textit{almost linear} if
\begin{align}\label{eq:R-cond3}
d \cdot \base^{p,q}_\zeta(\Delta) \subset \base^{p,q+1}_\zeta(\Delta)\oplus F^{p+\frac{1}{N}}V(\conf)^{p+q+1}_\zeta(\Delta),
\end{align}
or equivalently,
\begin{align*}
\gr d \cdot \gr \base^{p,q}_\zeta(\Delta) \subset \gr \base^{p,q+1}_\zeta(\Delta).
\end{align*}
Then $\gr d$ acts on $\gr\conf = \C[\partial]\otimes(\gr\base)$ and the $q$-th cohomologies $H^{p, q}(\gr \conf, \gr d)$ of $(\gr\conf^{p, \bullet}, \gr d)$ and $H^{p, q}(\gr \base, \gr d)$ of $(\gr\base^{p, \bullet}, \gr d)$ satisfy that
\begin{align*}
&H^{p, q}(\gr \conf, \gr d) = \C[\partial]\otimes H^{p, q}(\gr \base, \gr d),\\
&H^{p, q}(\gr \base, \gr d) = \sum_{\zeta \in \monoid_{\geq0}}\bigoplus_{\Delta \in \R}H^{p, q}(\gr \base, \gr d)_\zeta(\Delta),\\
&H^{p, q}(\gr \base, \gr d)_\zeta(\Delta) = \frac{\Ker(\gr d \colon \gr \base^{p,q}_\zeta(\Delta) \rightarrow \gr \base^{p,q+1}_\zeta(\Delta))}{\Img(\gr d \colon \gr \base^{p,q-1}_\zeta(\Delta) \rightarrow \gr \base^{p,q}_\zeta(\Delta))}.
\end{align*}
An almost linear odd differential $d$ is called \textit{good} if
\begin{align}\label{eq:R-cond4}
\base^{p,q}_\zeta(\Delta)\cap d^{-1}(F^{p+\frac{1}{N}}V^{p+q+1}_\zeta(\Delta)) \subset d(\base^{p, q-1}_\zeta(\Delta)) + F^{p+\frac{1}{N}}V_\zeta(\Delta),\quad
p+q \neq0,
\end{align}
or equivalently,
\begin{align*}
H^{p, q}(\gr \conf, \gr d) = 0,\quad
p+q \neq 0.
\end{align*}
Since $\gr\conf^{p,-p}_\zeta(\Delta)\cap\Img\gr d = 0$, we have
\begin{align*}
H^{p,-p}(\gr\conf,\gr d) _\zeta(\Delta) = \conf^{p,-p}_\zeta(\Delta)\cap\Ker\gr d = \conf^{p, -p}_\zeta(\Delta)\cap d^{-1}\left(F^{p+\frac{1}{N}}V(\conf)^1_\zeta(\Delta)\right).
\end{align*}
Then $H^{p,-p}(\gr\conf,\gr d)$ is naturally embedded into $\gr\conf$ and thus is a non-linear Lie conformal superalgebra. Therefore by \cite[Theorem 4.7]{DK},
\begin{align*}
H(\gr V(\conf), \gr d) \simeq V(H(\gr \conf, \gr d)).
\end{align*}
Now, using the same proof as that of \cite[Lemma 4.11]{DK}, it follows that
\begin{align*}
\gr d(\gr V(\conf)^{p,q}) \cap \gr V(\conf)^{p,q+1}_\zeta(\Delta) = \gr d(\gr V(\conf)^{p,q}_\zeta(\Delta)).
\end{align*}
Hence there exist natural embeddings
\begin{align*}
H^{p,q}(\gr V(\conf)_\zeta(\Delta), \gr d) \hookrightarrow H^{p,q}(\gr V(\conf), \gr d).
\end{align*}
Since $H^{p,q}(\gr V(\conf), \gr d) = H^{p,q}(V(\gr\conf), \gr d) = 0$ for $p+q \neq 0$, we have
\begin{align*}
H^{p,q}(\gr V(\conf)_\zeta(\Delta), \gr d) = 0,\quad
p+q \neq 0.
\end{align*}
Since $\dim\base<\infty$, we also have $\dim V(\conf)_\zeta(\Delta)<\infty$ so that for each $n \in \Z_{\geq0}$, $F^pV(\conf)^n_\zeta(\Delta) = 0$ with $p\gg0$. By \cite[Lemma 4.2]{DK},
\begin{align}
\label{eq:locallyfinite1}&F^pV(\conf)^n_\zeta(\Delta)\cap\Ker d = F^pV(\conf)^n_\zeta(\Delta)\cap\Img d,\quad
n\neq 0,\\
\label{eq:locallyfinite2}&F^pV(\conf)^0_\zeta(\Delta)\cap d^{-1}(F^{p+\frac{1}{N}}V(\conf)^1_\zeta(\Delta)) = F^{p+\frac{1}{N}}V(\conf)^0_\zeta(\Delta)+F^pV(\conf)^0_\zeta(\Delta)\cap\Ker d.
\end{align}
Let $\{e_\alpha\}_{\alpha\in \Cohindex(p,\zeta,\Delta)}$ be a parity-homogeneous basis of $\base^{p,-p}_\zeta(\Delta)\cap d^{-1}(F^{p+\frac{1}{N}}V(\conf)^1_\zeta(\Delta))$. Then
\begin{align*}
\left\{\ 
\NO{(\partial^{(k_1)}e_{\alpha_1})\cdots(\partial^{(k_s)}e_{\alpha_s})}
\ \middle|\ 
\begin{lgathered}
k_j \in \Z_{\geq0}, i_j \in \Cohindex(p_j, \zeta_j, \Delta_j),\sum_j p_j = p, \sum_j \zeta_j \leq \zeta,\\
\sum_j \Delta_j = \Delta,\ i_j \leq i_{j+1},\ 
i_j < i_{j+1}\ \mathrm{if}\ p(e_{\alpha_j}) = \bar{1}
\end{lgathered}
\ \right\}
\end{align*}
gives a basis of $H^{p,-p}(\gr V(\conf), \gr d)_\zeta(\Delta)$ since
\begin{align*}
H^{p,-p}(\gr V(\conf), \gr d)_\zeta(\Delta) \simeq V(H(\gr\base, \gr d))^{p,-p}_\zeta(\Delta).
\end{align*}
By \eqref{eq:locallyfinite2},
\begin{align*}
\base^{p,-p}_\zeta(\Delta)\cap d^{-1}(F^{p+\frac{1}{N}}V(\conf)^1_\zeta(\Delta)) \subset
F^{p+\frac{1}{N}}V(\conf)^0_\zeta(\Delta)+F^pV(\conf)^0_\zeta(\Delta)\cap\Ker d,
\end{align*}
there exist $E_\alpha \in F^pV(\conf)^0_\zeta(\Delta)\cap\Ker d$ and $\widetilde{e}_\alpha \in F^{p+\frac{1}{N}}V(\conf)^0_\zeta(\Delta)$ for $\alpha \in \Cohindex(p,\zeta,\Delta)$ such that $E_\alpha = e_\alpha + \widetilde{e}_\alpha$ and $p(E_\alpha) = p(e_\alpha) = p(\widetilde{e}_\alpha)$. Recall that $H^0(V(\conf), d) = V(\conf) \cap \Ker d$ so that $E_\alpha \in H^0(V(\conf), d)$ for all $\alpha \in \Cohindex = \bigcup_{p, \zeta, \Delta}\Cohindex(p,\zeta,\Delta)$. Using the same proof as that of \cite[Lemma 4.13]{DK}, we have that $H^n(V(\conf), d) = 0$ for all $n \neq 0$ and $H^0(V(\conf), d)$ is freely generated by $\{\partial^{(k)}E_\alpha\}_{k\in\Z_{\geq0}, \alpha \in \Cohindex}$. Let $H(\conf)=\C[\partial]\otimes\Span_\C\{E_\alpha\}_{\alpha \in \Cohindex}$. Then $H(\conf)$ has a structure of non-linear Lie conformal superalgebra inherited from $\conf$ so that $H(\conf)$ is pregraded by $\monoid_{\geq0}$. As a consequence, we obtain the following theorem as a generalization of \cite[Lemma 4.18, Theorem 4.19]{DK}:

\begin{theorem}\label{thm:R-coh}
Let $\conf$ be a non-linear Lie conformal superalgebra freely generated by a $\monoid_{>0} \times \Charge \times \R$-graded finite-dimensional subspace $\base$ as a $\C[\partial]$-module satisfying \eqref{eq:R-cond1} and $d$ be an odd differential on $V(\conf)$ satisfying \eqref{eq:R-cond2}\eqref{eq:R-cond3}\eqref{eq:R-cond4}. Then
\begin{enumerate}[label=(\arabic*)]
\item $H^n(V(\conf), d) = 0$ for $n \neq 0$.
\item Let $\{e_\alpha\}_{\alpha \in \Cohindex(p, \zeta, \Delta)}$ be a parity-homogeneous basis of $\base^{p, -p}_\zeta(\Delta) \cap d^{-1}(F^{p+\frac{1}{N}}V(\conf)^1_\zeta(\Delta))$ for $p \in \frac{1}{N}\Z, \zeta \in \monoid_{>0}, \Delta \in \R$ with some index set $\Cohindex(p, \zeta, \Delta)$. Then there exist elements $\{E_\alpha\}_{\alpha \in \Cohindex(p, \zeta, \Delta)}$ in $F^pV(\conf)^0_\zeta(\Delta)\cap\Ker d$ such that
\begin{align*}
E_\alpha -e_\alpha \in F^{p+\frac{1}{N}}V(\conf)^0_\zeta(\Delta),\quad
p(E_\alpha) = p(e_\alpha).
\end{align*}
\item Let $\pi \colon \Ker d \twoheadrightarrow H(V(\conf), d)$ be the canonical projection and $\Cohindex =\bigcup_{p, \zeta, \Delta}\Cohindex(p, \zeta, \Delta)$. Then $H^0(V(\conf), d)$ is freely generated by $\{\pi(\partial^{(k)}E_\alpha)\}_{k \in \Z_{\geq0}, \alpha \in \Cohindex}$.
\item Let $H(\conf) = \C[\partial]\otimes\Span_\C\{\pi(E_\alpha)\mid \alpha \in \Cohindex\}$. Then $H(\conf)$ has a structure of a non-linear Lie conformal superalgebra such that $H^0(V(\conf), d) \simeq V(H(\conf))$.
\end{enumerate}
\end{theorem}

\section{Commutativity of Zhu and Cohomology functor}\label{sec:commutativity}

Let $\monoid$ be a discrete additive subgroup of $\R$ and $\Charge = \{(p, q) \in \frac{1}{N}\Z \times \frac{1}{N}\Z \mid p+q \in \Z_{\geq0}\}$ with some $N \in \N$. Denote by $\monoid_{\geq0} = \monoid \cap \R_{\geq0}$, by $\monoid_{>0} = \monoid \cap \R_{>0}$ and by $\zeta_-$ the largest number in $\monoid_{\geq0}$ that is strictly smaller than $\zeta$ for each $\zeta \in \monoid_{>0}$ as in Section \ref{sec:non-linear salg}. Let $V$ be a vertex superalgebra with a Hamiltonian operator $H$ whose eigenvalues belong to $\R$ and $g \in \operatorname{Aut}_H V$ be a diagonalizable automorphism with modulus $1$ eigenvalues. Suppose that $V$ is freely generated by free generators $\{e_\alpha\}_{\alpha \in I}$ homogeneous for $g$, $H$ and the parity, with an ordered index set $I$ and is pregraded by $\monoid_{\geq0}$ in the sense of Section \ref{sec:def-va}.
\smallskip

Moreover, suppose that
\begin{enumerate}[label=(A.\arabic*)]
\item $I$ is finite,
\item the set $\{e_\alpha\}_{\alpha \in I}$ is $\Charge$-graded, denoted by $\deg_\Charge e_\alpha \in \Charge$, such that
\begin{align*}
[{\base_{\zeta_1}^{p_1, q_1}}_\lambda\base_{\zeta_2}^{p_2, q_2}] \subset \bigoplus_{l=0}^\infty\C[\lambda]\otimes V^{p_1+p_2-l, q_1+q_2-l}_{(\zeta_1+\zeta_2)_-}
\end{align*}
for all $\zeta_j \in \monoid_{>0}$ and $(p_j, q_j) \in \Charge$ for $j =1, 2$, where
\begin{align*}
&\base_\zeta^{p, q}=\Span_\C\{e_\alpha \mid \alpha \in I, \zeta_\alpha=\zeta, \deg_\Charge e_\alpha = (p, q)\},\\
&V^{p, q}_\zeta = \sum_{\begin{subarray}{c}k_1, \ldots, k_s \geq 0, \sum_i \zeta_i \leq \zeta \\ \sum_ip_i = p, \sum_i q_i =q \end{subarray}}\NO{(\partial^{k_1}\base^{p_1, q_1}_{\zeta_1})\cdots(\partial^{k_s}\base^{p_s, q_s}_{\zeta_s})},
\end{align*}
\item $V$ is $\Charge$-graded by
\begin{align*}
V = \bigoplus_{(p, q) \in \Charge}V^{p, q},\quad
V^{p, q} = \sum_{\zeta \in \monoid_{\geq0}}V^{p, q}_\zeta.
\end{align*}
\end{enumerate}
Then subspaces $V^{p, q}$ are invariant under the actions of $\partial$, $g$, and $H$ and are homogeneous for the parity. Denote by $A(\Delta) = \{ a \in A \mid \Delta_a = \Delta\}$ for any subspace $A$ in $V$. By Proposition \ref{prop:DK-V(R)}, $\conf = \C[\partial]\otimes\base$ with $\base = \Span_\C\{e_\alpha \mid \alpha \in I\}$ is a non-linear Lie conformal superalgebra such that $V \simeq V(\conf)$. Moreover,
\begin{align*}
\conf^{p_1, q_1}_{\zeta_1}(\Delta_1)_{(n)}\conf^{p_2, q_2}_{\zeta_2}(\Delta_2) \subset \bigoplus_{l=0}^\infty V^{p_1+p_2+l, q_1+q_2-l}_{(\zeta_1+\zeta_2)_{-}}(\Delta_1+\Delta_2-n-1)
\end{align*}
for all $\Delta_j \in \R$, $\zeta_j \in \monoid_{>0}$, $(p_j, q_j) \in \Charge$ for $j=1,2$ and $n \in \Z_{\geq0}$, where
\begin{align*}
\conf^{p, q}_\zeta(\Delta) = \bigoplus_{n=0}^\infty\partial^n\base^{p, q}_\zeta(\Delta-n) = \{ a \in \conf \mid \Delta_a = \Delta, \zeta_a = \zeta, \deg_\Charge a = (p, q) \}.
\end{align*}
Let $V^n_\zeta = \sum_{p+q = n}V^{p,q}_\zeta$ for $\zeta \in \monoid_{\geq0}$ and $n \in \Z$ and $V^n = \sum_{\zeta \in \monoid_{\geq0}}V^n_\zeta$. Define a descending filtration $F^pV_\zeta^n(\Delta)$ on $V_\zeta^n(\Delta)$ for $\zeta \in \monoid_{\geq0}$, $n \in \Z$, $\Delta \in \R$, $p \in \frac{1}{N}\Z$ by
\begin{align*}
F^pV_\zeta^n(\Delta)
= \Span_\C\left\{\ 
\NO{a_1 \cdots a_s}\ \in V\ \middle|\ 
\begin{lgathered}
a_k \in \conf^{p_k, q_k}_{\zeta_k}(\Delta_k),\ \sum_k p_k \geq p,\\
\sum_k(p_k+q_k)=n, \sum_k\zeta_k \leq \zeta, \sum_k\Delta_k = \Delta
\end{lgathered}
\ \right\}.
\end{align*}
Let $d$ be an odd differential on $V$ in the sense of Section \ref{sec:cohomology on conf} such that
\begin{enumerate}[start=4, label=(A.\arabic*)]
\item \label{it: d-comm} $d$ commutes with $g$ and $H$.
\item $d \cdot F^pV^n_\zeta(\Delta) \subset F^{p}V^{n+1}_\zeta(\Delta)$.
\item $d \cdot \base^{p, q}_\zeta(\Delta) \subset \base^{p, q+1}_\zeta(\Delta)\oplus F^{p+\frac{1}{N}}V^{p+q+1}_\zeta(\Delta)$.
\item $\base^{p, q}_\zeta(\Delta) \cap d^{-1}(F^{p+\frac{1}{N}}V^{p+q+1}(\Delta)) \subset d(\base^{p, q-1}_\zeta(\Delta)) + F^{p+\frac{1}{N}}V_\zeta^{p+q}(\Delta)$.
\end{enumerate}
Then $d \cdot V^n \subset V^{n+1}$ so that $(V^\bullet, d)$ forms a cochain complex. Denote by $H(V, d)$ the cohomology of the complex $(V^\bullet, d)$ as in Section \ref{sec:cohomology on conf}. By \ref{it: d-comm}, $H$ induces a Hamiltonian operator on $H(V, d)$, which we denote by $H$ by abuse of notation, and $g$ also defines an automorphism on $H(V, d)$ commuting with $H$, which we also denote by $g$. Thus we can define the $(g, H)$-twisted Zhu algebra of $H(V, d)$
\begin{align*}
\gH{\Zhu}H(V, d) = \gH{H(V, d)}/\gH{J}(H(V, d)).
\end{align*}
Let
\begin{align*}
\projZhu \colon \gH{V} \twoheadrightarrow \gH{\Zhu}V,\quad
\pi\colon\Ker d \twoheadrightarrow H(V, d)
\end{align*}
be the canonical projections and $\overline{d} \colon \gH{\Zhu}V \rightarrow \gH{\Zhu}V$ be an odd differential on $\gH{\Zhu}V$ defined by $\overline{d}\cdot\projZhu(a) = \projZhu(d a)$ for $a \in \gH{V}$. Since $d(a_{(n)}b) = (da)_{(n)}b + (-1)^{p(a)}a_{(n)}(db)$ for $a, b\in V$, $d$ is a derivation for the $*_n$-products. Then $d\cdot \gH{J} \subset \gH{J}$ so that $\overline{d}$ is well-defined. Let
\begin{align*}
\Coh{\projZhu} \colon \gH{H(V, d)} \twoheadrightarrow \gH{\Zhu}H(V, d),\quad
\Coh{\pi} \colon \Ker \overline{d} \twoheadrightarrow H(\gH{\Zhu}V, \overline{d})
\end{align*}
the canonical projections. The main purpose of this section is to prove the following theorem. The proof is essentially the same as that of \cite[Theorem 4.20]{DK}.

\begin{theorem}\label{thm:main_cohomology_thm}
Let $V$ be a vertex superalgebra with a Hamiltonian operator $H$ whose eigenvalues belong to $\R$ and $g \in \operatorname{Aut}_H V$ be a diagonalizable automorphism with modulus $1$ eigenvalues. Suppose that $V$ is freely generated by free generators $\{e_\alpha\}_{\alpha \in I}$ homogeneous with respect to $g$, $H$ and the parity, with an ordered index set $I$, that $V$ is pregraded by $\monoid_{\geq0}$ in the sense of Section \ref{sec:def-va}, and the conditions (A.1)--(A.3). Let $d \colon V \rightarrow V$ be an odd differential on $V$ satisfying the conditions (A.4)--(A.7). Then, the canonical map
\begin{align*}
\psi \colon \gH{\Zhu}H(V, d) \rightarrow H(\gH{\Zhu}V, \overline{d}),\quad
(\Coh{\projZhu} \circ \pi)(a) \mapsto (\Coh{\pi}\circ\projZhu)(a)
\end{align*}
is a well-defined isomorphism of associative superalgebras.
\end{theorem}

First of all, by the assumptions (A.1)--(A.7), the non-linear Lie conformal superalgebra $\conf$ satisfies all assumptions in Theorem \ref{thm:R-coh}. Thus $H^n(V,d) = 0$ for $n \neq 0$. Moreover, for a homogeneous basis $\{\widetilde{e}_\alpha\}_{\alpha \in \Cohindex(p, \zeta, \Delta)}$ of $\base^{p, -p}_\zeta(\Delta) \cap d^{-1}(F^{p+\frac{1}{N}}V^1_\zeta(\Delta))$ with respect to $g$ and the parity, there exist elements $\{E_\alpha\}_{\alpha \in \Cohindex(p, \zeta, \Delta)}$ in $F^pV^0_\zeta(\Delta)\cap\Ker d$ such that
\begin{align*}
E_\alpha -\widetilde{e}_\alpha \in F^{p+\frac{1}{N}}V^0_\zeta(\Delta).
\end{align*}
We may assume that $E_\alpha$ has the same eigenvalue of $g$ and parity as $\widetilde{e}_\alpha$. Let $H(\base) = \Span_\C\{\pi(E_\alpha)\mid\alpha\in \Cohindex\}$ with $\Cohindex=\bigcup_{p,\zeta,\Delta}\Cohindex(p,\zeta, \Delta)$ and $H(\conf) = \C[\partial]\otimes H(\base) \subset H(V, d)$. Then $H(\conf)$ has a structure of a non-linear Lie conformal superalgebra such that
\begin{align*}
H(V,d) = H^0(V,d) \simeq V(H(\conf))
\end{align*}
by Theorem \ref{thm:R-coh}. Let
\begin{align*}
&\gH{\Zhu} H(\conf) = \Coh{\projZhu}(\gH{H(\conf)}) = \gH{H(\conf)}/(\partial + H)\gH{H(\conf)} = \Coh{\projZhu}(H(\base)) \simeq \gH{H(\base)},\\
&\gH{H(\conf)} = H(\conf)\cap\pi(\gH{V}) = \C[\partial]\otimes \gH{H(\base)},\\
&\gH{H(\base)} = \Span_\C\{\pi(E_\alpha)\mid\alpha \in \gH{{\Cohindex}}\},\quad
\gH{{\Cohindex}} = \{\alpha \in \Cohindex \mid E_\alpha \in \gH{V}\}.
\end{align*}
By Theorem \ref{thm:main PBW thm}, $\gH{\Zhu} H(\conf)$ has a structure of a non-linear Lie superalgebra such that
\begin{align*}
\gH{\Zhu} H(V, d) \simeq U(\gH{\Zhu} H(\conf)).
\end{align*}
Thus, we get the following result:

\begin{proposition}\label{prop:ZH-PBW}
$\gH{\Zhu} H(V, d)$ is PBW generated by $\{(\Coh{\projZhu}\circ\pi)(E_\alpha) \mid \alpha \in \gH{\Cohindex}\}$.
\end{proposition}

Next, let $\gH{\Zhu}\conf = \projZhu(\gH{\conf}) = \projZhu(\gH{\base}) \simeq \gH{\base}$ for $\gH{\conf} = \C[\partial]\otimes\gH{\base}$ with $\gH{\base} = \base \cap \gH{V}$. Then, again by Theorem \ref{thm:main PBW thm}, $\gH{\Zhu}\conf$ is a non-linear Lie superalgebra by
\begin{align}\label{eq:tau_bracket}
[\projZhu(a), \projZhu(b)] = \sum_{j=0}^\infty\binom{\Delta_a-1}{j}(\rho\circ\projZhu)(a_{(j)}b),\quad
a, b \in \gH{\base}
\end{align}
and satisfies that
\begin{align*}
\gH{\Zhu} V \simeq U(\gH{\Zhu}\conf),
\end{align*}
where $\rho \colon \gH{\Zhu} V \hookrightarrow T(\projZhu(\gH{\base}))$. Let $(\gH{\base})^{p,q}_\zeta = \gH{\base} \cap \base^{p,q}_\zeta$. Then
\begin{align}\label{eq:tau-decomp}
\projZhu(\gH{\base}) = \bigoplus_{(p,q) \in \Charge, \zeta\in\monoid_{>0}}\projZhu((\gH{\base})^{p, q}_\zeta).
\end{align}
Let $(\gH{\conf})^{p, q}_\zeta = \gH{\conf}\cap\conf^{p,q}_\zeta = \C[\partial]\otimes(\gH{\base})^{p,q}_\zeta$. Using the equation $\projZhu(\partial^{(k)}a) = \binom{-\Delta_a}{k}\projZhu(a)$
 for $a\in \gH{\base}$, we have $\projZhu((\gH{\conf})^{p, q}_\zeta) = \projZhu((\gH{\base})^{p, q}_\zeta)$. Moreover, by \eqref{eq:tau_bracket},
\begin{align*}
&[\projZhu((\gH{\base})^{p_1, q_1}_{\zeta_1}), \projZhu((\gH{\base})^{p_2, q_2}_{\zeta_2})]
= (\rho\circ\projZhu)\left([{(\gH{\base})^{p_1, q_1}_{\zeta_1}}_*(\gH{\base})^{p_2, q_2}_{\zeta_2}]\right)\\
&\subset (\rho\circ\projZhu)\left(\bigoplus_{l=0}^\infty V^{p_1+p_2+l, q_1+q_2-l}_{(\zeta_1+\zeta_2)_-}\cap \gH{V}\right)
= \bigoplus_{l=0}^\infty T(\projZhu(\gH{\base}))^{p_1+p_2+l, q_1+q_2-l}_{(\zeta_1+\zeta_2)_-}
\end{align*}
Let
\begin{align*}
F^pU(\projZhu(\gH{\base}))^n_\zeta
=\Span_\C\left\{\ \projZhu(a_1)\cdots\projZhu(a_s)
\ \middle|\ 
\begin{lgathered}
a_k \in (\gH{\base})^{p_k, q_k}_{\zeta_k}, \sum_k(p_k+q_k)=n,\\
\sum_k p_k \geq p,\ \sum_k \zeta_k \leq \zeta
\end{lgathered}
\ \right\}
\end{align*}
for $p \in \frac{1}{N} \Z$, $n \in \Z_{\geq0}$, $\zeta \in \monoid_{\geq0}$.

\begin{lemma}\label{lem:FV_to_U}
$\projZhu(F^pV^n_\zeta\cap \gH{V}) = F^pU(\projZhu(\gH{\base}))^n_\zeta$.
\begin{proof}
Recall that
\begin{align*}
&F^pV^n_\zeta\cap \gH{V}\\
&= \Span_\C\left\{\ 
\NO{(\partial^{(k_1)}e_{\alpha_1}) \cdots (\partial^{(k_s)}e_{\alpha_s})}
\ \middle|\ 
\begin{lgathered}
k_j \in \Z_{\geq0}, e_{\alpha_j} \in \base^{p_j, q_j}_{\zeta_j},
\epsilon(\overrightarrow{i}) =0,\\
\mathrm{for}\ 
\overrightarrow{i} = (i_1, \ldots, i_s),
\sum_j(p_j+q_j) = n,\\
\sum_jp_j \geq p,\ 
\sum_j\zeta_j \leq \zeta
\end{lgathered}
\ \right\}\\
&= \Span_\C\left\{\ 
(\partial^{(k_1)}e_{\alpha_1})* \cdots *(\partial^{(k_s)}e_{\alpha_s})
\ \middle|\ 
\begin{lgathered}
k_j \in \Z_{\geq0}, e_{\alpha_j} \in \base^{p_j, q_j}_{\zeta_j},
\epsilon(\overrightarrow{i}) =0,\\
\mathrm{for}\ 
\overrightarrow{i} = (i_1, \ldots, i_s),
\sum_j(p_j+q_j) = n,\\
\sum_jp_j \geq p,\ 
\sum_j\zeta_j \leq \zeta
\end{lgathered}
\ \right\}.
\end{align*}
Therefore by Theorem \ref{thm:J-PBW},
\begin{align*}
&\projZhu(F^pV^n_\zeta\cap \gH{V})
= \Span_\C\left\{\ 
\projZhu(e_{\alpha_1})* \cdots *\projZhu(e_{\alpha_s})
\ \middle|\ 
\begin{lgathered}
e_{\alpha_j} \in (\gH{\base})^{p_j, q_j}_{\zeta_j},
\sum_j(p_j+q_j) = n,\\
\sum_jp_j \geq p,\ 
\sum_j\zeta_j \leq \zeta
\end{lgathered}
\ \right\}\\
&=F^pU(\projZhu(\gH{\base}))^n_\zeta.
\end{align*}
\end{proof}
\end{lemma}

Recall that $\overline{d} \colon \gH{\Zhu} V \rightarrow \gH{\Zhu} V$ is a derivation on $\gH{\Zhu} V \simeq U(\gH{\Zhu}\conf) = U(\projZhu(\gH{\base}))$ defined by $\projZhu(a) \mapsto \projZhu(d \cdot a)$. Then, by the assumption (A.5) and Lemma \ref{lem:FV_to_U},
\begin{align}\label{eq:dbar_FU}
\begin{split}
\overline{d} \cdot F^pU(\projZhu(\gH{\base}))^n_\zeta
= \overline{d} \cdot \projZhu(F^pV^n_\zeta\cap \gH{V})
= \projZhu(d \cdot F^pV^n_\zeta\cap \gH{V})\\
\subset \projZhu(F^pV^{n+1}_\zeta\cap \gH{V})
= F^pU(\projZhu(\gH{\base}))^{n+1}_\zeta.
\end{split}
\end{align}

\begin{lemma}\label{lem:FV_to_U2}
\begin{align*}
\projZhu((\gH{\base})^{p,q}_\zeta) \cap \overline{d}^{-1}(F^{p+\frac{1}{N}}U(\projZhu(\gH{\base}))^{p+q+1}_\zeta)
= \projZhu\left((\gH{\base})^{p,q}_\zeta \cap d^{-1}(F^{p+\frac{1}{N}}V^{p+q+1}_\zeta\cap \gH{V})\right).
\end{align*}
\begin{proof}
It suffices to show $\subset$ since $\supset$ is clear from Lemma \ref{lem:FV_to_U} and \eqref{eq:dbar_FU}. Let $a \in (\gH{\base})^{p,q}_\zeta$ such that $\overline{d}\cdot\projZhu(a) = \projZhu(da) \in F^{p+\frac{1}{N}}U(\projZhu(\gH{\base}))^{p+q+1}_\zeta$. By the assumption (A.6), there exist $a' \in  (\gH{\base})^{p,q+1}_\zeta$ and $A \in F^{p+\frac{1}{N}}V^{p+q+1}_\zeta\cap \gH{V}$ such that $da = a' + A$. Thus $\projZhu(a') = \projZhu(da)-\projZhu(A) \in F^{p+\frac{1}{N}}U(\projZhu(\gH{\base}))^{p+q+1}_\zeta$ by Lemma \ref{lem:FV_to_U}. But $\projZhu(a') \in \projZhu((\gH{\base})^{p,q+1}_\zeta)\cap F^{p+\frac{1}{N}}U(\projZhu(\gH{\base}))^{p+q+1}_\zeta = 0$ so that $\projZhu(a') = 0$. Hence $a' = 0$. Therefore $da = A$ and thus $a \in (\gH{\base})^{p,q}_\zeta \cap d^{-1}(F^{p+\frac{1}{N}}V^{p+q+1}_\zeta\cap \gH{V})$. This completes the proof.
\end{proof}
\end{lemma}

Now, by the assumption (A.6) and Lemma \ref{lem:FV_to_U},
\begin{align}\label{eq:dbar_almostlinear}
\begin{split}
\overline{d} \cdot \projZhu((\gH{\base})^{p, q}_\zeta)
&= \projZhu(d\cdot (\gH{\base})^{p, q}_\zeta)
\subset \projZhu((\gH{\base})^{p, q+1}_\zeta \oplus (F^{p+\frac{1}{N}}V^{p+q+1}_\zeta\cap \gH{V}))\\
&=\projZhu((\gH{\base})^{p, q+1}_\zeta)\oplus F^{p+\frac{1}{N}}U(\projZhu(\gH{\base}))^{p+q+1}_\zeta
\end{split}
\end{align}
so that $\overline{d}$ is almost linear. Moreover, by the assumption (A.7), Lemma \ref{lem:FV_to_U} and Lemma \ref{lem:FV_to_U2},
\begin{align}\label{eq:dbar_good}
\begin{split}
&\projZhu((\gH{\base})^{p, q}_\zeta) \cap \overline{d}^{-1}(F^{p+\frac{1}{N}}U(\projZhu(\gH{\base}))^{p+q+1}_\zeta)
= \projZhu\left((\gH{\base})^{p,q}_\zeta \cap d^{-1}(F^{p+\frac{1}{N}}V^{p+q+1}_\zeta\cap \gH{V})\right)\\
&\subset \projZhu(d((\gH{\base})^{p, q-1}_\zeta) + (F^{p+\frac{1}{N}}V_\zeta^{p+q})\cap \gH{V})\\
&= \overline{d}\cdot \projZhu((\gH{\base})^{p, q-1}_\zeta) + F^{p+\frac{1}{N}}U(\projZhu(\gH{\base}))^{p+q}_\zeta,\quad
\mathrm{for}\ 
p+q \neq0.
\end{split}
\end{align}
Thus, $\overline{d}$ is good. Let $\Cohindex(p, \zeta) = \bigcup_{\Delta}\Cohindex(p, \zeta, \Delta)$ and $\gH{\Cohindex}(p, \zeta) = \{ \alpha \in \Cohindex(p, \zeta) \mid \widetilde{e}_\alpha \in \gH{V}\}$. Then $\{\widetilde{e}_\alpha\}_{\alpha \in D_g(p, \zeta)}$ is a parity-homogeneous basis of $(\gH{\base})^{p, -p}_\zeta \cap d^{-1}(F^{p+\frac{1}{N}}V^1_\zeta\cap \gH{V})$. Then $E_\alpha \in F^pV^0_\zeta\cap \gH{V} \cap \Ker d$ for $\alpha \in \gH{\Cohindex}(p, \zeta)$ has the same parity as $\widetilde{e}_\alpha$ and satisfies that $E_\alpha - \widetilde{e}_\alpha \in F^{p+\frac{1}{N}}V^0_\zeta \cap \gH{V}$. By Lemma \ref{lem:FV_to_U2}, $\{\projZhu(\widetilde{e}_\alpha)\}_{\alpha \in \gH{\Cohindex}(p, \zeta)}$ is a parity-homogeneous basis of $\projZhu((\gH{\base})^{p,-p}_\zeta) \cap \overline{d}^{-1}(F^{p+\frac{1}{N}}U(\projZhu(\gH{\base}))^{1}_\zeta)$. By Lemma \ref{lem:FV_to_U}, $\projZhu(E_\alpha) \in F^pU(\projZhu(\gH{\base}))^0_\zeta\cap \Ker \overline{d}$ for $\alpha \in \gH{\Cohindex}(p, \zeta)$ has the same parity as $\projZhu(\widetilde{e}_\alpha)$ and satisfies that $\projZhu(E_\alpha) - \projZhu(\widetilde{e}_\alpha) \in F^{p+\frac{1}{N}}U(\projZhu(\gH{\base}))^0_\zeta$. Therefore by Theorem \ref{thm:r-coh} together with \eqref{eq:tau-decomp}\eqref{eq:dbar_FU}\eqref{eq:dbar_almostlinear}\eqref{eq:dbar_good}, we have $H^n(\gH{\Zhu} V, \overline{d}) = 0$ for $n \neq 0$ and  $H^0(\gH{\Zhu} V, \overline{d})$ is PBW generated by $\{(\Coh{\pi}\circ\projZhu)(E_\alpha)\}_{\alpha \in \gH{\Cohindex}}$. Therefore, we obtain the following result:

\begin{proposition}\label{prop:HZ-PBW}
$H(\gH{\Zhu} V, \overline{d})$ is PBW generated by $\{(\Coh{\pi}\circ\projZhu)(E_\alpha)\mid \alpha \in \gH{\Cohindex}\}$.
\end{proposition}

\hspace{-3.5mm}\textit{Proof of Theorem \ref{thm:main_cohomology_thm}}. Denote by $\gH{H(V, d)} = \bigoplus_{\bar{\gamma}}H(V, d)(\bar{\gamma}, \bar{\gamma})$ and by
\begin{align*}
\gH{J}(H(V, d)) = \gH{H(V, d)} \circ \gH{H(V, d)}
+ \sum_{
\begin{subarray}{c}
\bar{\gamma}_1 \neq \bar{\mu}_1, \bar{\gamma}_2 \neq \bar{\mu}_2 \\
\bar{\gamma}_1 + \bar{\gamma}_2 = \bar{\mu}_1 + \bar{\mu}_2
\end{subarray}
}
H(V, d)(\bar{\gamma}_1, \bar{\mu}_1)*H(V, d)(\bar{\gamma}_2, \bar{\mu}_2).
\end{align*}
Then $\gH{\Zhu} H(V, d) = \gH{H(V, d)}/\gH{J}(H(V, d))$. Recall that $\projZhu \colon \gH{V} \twoheadrightarrow \gH{\Zhu} V$ satisfies that $\projZhu(\Ker d \cap \gH{V}) \subset \Ker \overline{d}$ and $\projZhu(d \cdot \gH{V}) = \overline{d}(\gH{\Zhu} V)$. Thus, $\projZhu$ induces a canonical map
\begin{align*}
&\overline{\projZhu} \colon \gH{H(V,d)} = \Ker d\cap \gH{V} / \Img d \cap \gH{V} \rightarrow \Ker \overline{d}/\Img\overline{d} = H(\gH{\Zhu} V, \overline{d}),\\
&\pi(a) \mapsto (\Coh{\pi} \circ \projZhu)(a),\quad
a \in \Ker d \cap \gH{V}.
\end{align*}
Since $\pi$ is a vertex superalgebra morphism commuting with $g$ and $H$, $\pi(a)*_n\pi(b) = \pi(a*_n b)$ for $a, b \in \Ker d \cap \gH{V}$. Clearly, $(\Coh{\pi}\circ\projZhu)(a\circ b) = 0$ for $a, b \in \Ker d \cap \gH{V}$ and $(\Coh{\pi}\circ\projZhu)(a*b) = 0$ for $a \in \Ker d \cap V(\bar{\gamma}_1, \bar{\mu}_1)$ and $b \in \Ker d \cap V(\bar{\gamma}_2, \bar{\mu}_2)$ such that $\bar{\gamma}_1 \neq \bar{\mu}_1$, $\bar{\gamma}_2 \neq \bar{\mu}_2$ with $\bar{\gamma}_1 + \bar{\gamma}_2 = \bar{\mu}_1 + \bar{\mu}_2$. Thus $\gH{J}(H(V, d))$ is in the kernel of $\overline{\projZhu}$. Therefore
\begin{align*}
\psi \colon \gH{\Zhu} H(V, d) \rightarrow H(\gH{\Zhu} V, \overline{d}),\quad
(\Coh{\projZhu} \circ \pi)(a) \mapsto (\Coh{\pi}\circ\projZhu)(a)
\end{align*}
is a well-defined map. Moreover, $\psi((\Coh{\projZhu} \circ \pi)(a) \cdot (\Coh{\projZhu} \circ \pi)(b)) = \psi((\Coh{\projZhu} \circ \pi)(a*b)) = (\Coh{\pi}\circ\projZhu)(a*b) = (\Coh{\pi}\circ\projZhu)(a) \cdot (\Coh{\pi}\circ\projZhu)(b) = \psi((\Coh{\projZhu} \circ \pi)(a)) \cdot \psi((\Coh{\projZhu} \circ \pi)(b))$, so $\psi$ is a morphism of associative superalgebras. Finally, $\psi$ is an isomorphism since $\psi$ maps a PBW basis of $\gH{\Zhu} H(V, d)$ to one of $H(\gH{\Zhu} V, \overline{d})$ by Proposition \ref{prop:ZH-PBW} and Proposition \ref{prop:HZ-PBW}. This completes the proof. \qed

\section{Applications}\label{sec:applications}

\subsection{Automorphisms of non-linear Lie conformal superalgebras}

Let $\conf$ be a non-linear Lie conformal superalgebra freely generated by a $\monoid_{>0}$-graded subspace $\base$ as a $\C[\partial]$-module and $g$ be a parity-preserving linear automorphism on $\base$ such that
\begin{align*}
g(\base_\zeta)\subset \base_{\leq \zeta},\quad
\zeta \in \monoid_{>0}.
\end{align*}
Then we can extend $g$ to a $\C[\partial]$-module automorphism on $\conf$ by $g \circ \partial = \partial \circ g$. Hence $g$ is naturally extended to a superalgebra automorphism on $T(\conf)$ by $g(a\otimes b) = g(a)\otimes g(b)$ for $a, b \in T(\conf)$ and $g(1) = 1$. Then $g$ is called an automorphism of a non-linear Lie conformal superalgebra $\conf$ if
\begin{align*}
[g(a)_\lambda g(b)] = g([a_\lambda b]),\quad
a, b \in \conf,
\end{align*}
where $g$ is extended $\C[\lambda]$-linearly on $\C[\lambda]\otimes T(\conf)$.

\begin{lemma}\label{lem:auto_conf}
Let $g$ be an automorphism of a non-linear Lie conformal superalgebra $\conf$. Then $g$ induces an automorphism of the universal enveloping vertex superalgebra $V(\conf)$.
\begin{proof}
Recall that the $\lambda$-bracket on $\conf$ is uniquely extended to $T(\conf)$ and the normally ordered operator $T(\conf)\otimes T(\conf) \ni a\otimes b \rightarrow\ \NO{ab}\ \in T(\conf)$ is uniquely defined by equations \eqref{eq:lambda-NO-def}. Then $[g(a)_\lambda g(b)] = g([a_\lambda b])$ and $g(\NO{ab}) = \NO{g(a)g(b)}$ are well-defined for all $a, b \in T(\conf)$ since all equations in \eqref{eq:lambda-NO-def} are preserved by the action of $g$. Moreover, $g(M(\conf)_\zeta) \subset M(\conf)_\zeta$ by definition \eqref{eq:M(conf)_zeta} and thus $M(\conf) = \sum_\zeta M(\conf)_\zeta$ is $g$-invariant. Therefore $g$ induces an automorphism on $V(\conf) = T(\conf)/M(\conf)$, which satisfies the desired properties.
\end{proof}
\end{lemma}

\subsection{Affine vertex superalgebras}

Let $\mathfrak{g}$ be a simple basic classical Lie superalgebra. Then
\begin{align*}
\mathrm{Cur}^k(\mathfrak{g}):=\C[\partial]\otimes\mathfrak{g}
\end{align*}
has a structure of a non-linear Lie conformal superalgebra by
\begin{align}\label{eq:g-lambda}
[a_\lambda b] = [a, b] + k(a|b)\lambda,\quad
a, b \in \mathfrak{g}
\end{align}
for $k \in \C$, where $(\cdot|\cdot)$ is the normalized even supersymmetric invariant bilinear form on $\mathfrak{g}$ such that $(\theta|\theta) = 2$ for an even highest root $\theta$ of $\mathfrak{g}$. Here $\monoid=\Z$ and $\zeta_a = 1$ for all $a \in \mathfrak{g}$. Then the universal enveloping vertex algebra of $\mathrm{Cur}^k(\mathfrak{g})$ is isomorphic to the affine vertex superalgebra $V^k(\mathfrak{g})$ of $\mathfrak{g}$ at level $k$.

Let $L^\mathfrak{g}$ be a Sugawara conformal vector of $V^k(\mathfrak{g})$ for $k \neq - h^\vee$, where $h^\vee$ is the dual Coxeter number of $\mathfrak{g}$, and $x$ a semisimple element of $\mathfrak{g}$ such that $\ad x$ defines an $\R$-grading on $\mathfrak{g}$:
\begin{align*}
\mathfrak{g} = \bigoplus_{j \in \R}\mathfrak{g}_j,\quad
\mathfrak{g}_j = \{ a \in \mathfrak{g} \mid [x, a]=j a\}.
\end{align*}
Then $\omega_x := L^\mathfrak{g} + \partial x$ is a conformal vector of $V^k(\mathfrak{g})$ provided that $k \neq - h^\vee$. Thus
\begin{align*}
H := (\omega_x)_0 = L^\mathfrak{g}_0 - x_{(0)}
\end{align*}
defines a Hamiltonian operator on $V^k(\mathfrak{g})$ such that $\Delta_a = 1-j$ for $a \in \mathfrak{g}_j$, which is well-defined for all $k \in \C$. 

Let $g$ be a diagonalizable Lie superalgebra automorphism of $\mathfrak{g}$ with modulus $1$ eigenvalues such that $g(x) = x$ and $(g(a)|g(b)) = (a|b)$ for $a, b \in \mathfrak{g}$. Note that the first condition guarantees that $g \circ \ad x = \ad x \circ g$ and the second condition is automatically true if $h^\vee \neq 0$ since the normalized form is proportional to the Killing form. Then $g$ induces an automorphism of $V^k(\mathfrak{g})$ and commutes with $H$ since $g$ preserves the conformal weights. Denote by $\mathfrak{g}^h = \{ a \in \mathfrak{g} \mid h(a) = a\}$ for any automorphism $h$ of $\mathfrak{g}$ and $U(\mathfrak{g}^h)$ the universal enveloping algebra of $\mathfrak{g}^h$. Note that
\begin{align*}
\twist=\mathrm{e}^{2\pi iH}= \mathrm{e}^{-2\pi i x_{(0)}}\in \operatorname{Aut}_H V^k(\mathfrak{g}).
\end{align*}
The following gives a full generalization of results in \cite{FZ, SRW, Yang}:

\begin{theorem}\label{thm:twisted Zhu of affine}
Let $\mathfrak{g}$ be a simple basic classical Lie superalgebra, $L^\mathfrak{g}$ the Sugawara conformal vector of $V^k(\mathfrak{g})$, $x$ a semisimple element of $\mathfrak{g}$ such that $\ad x$ defines an $\R$-grading on $\mathfrak{g}$ and $g$ a Lie superalgebra automorphism of $\mathfrak{g}$ such that $g(x) = x$ and $(g(a)|g(b)) = (a|b)$ for $a, b \in \mathfrak{g}$, which induces an automorphism of $V^k(\mathfrak{g})$ and commutes with a Hamiltonian operator $H=L^\mathfrak{g}_0 - x_{(0)}$ of $V^k(\mathfrak{g})$. Then 
\begin{align*}
\Zhu_{g, H}V^k(\mathfrak{g}) \simeq U(\mathfrak{g}^{g\circ \twist^{-1}}).
\end{align*}
\begin{proof}
Let $\projZhu \colon V^k(\mathfrak{g})^{g\circ\twist^{-1}} \twoheadrightarrow \Zhu_{g, H}V^k(\mathfrak{g})$ be the canonical projection. By Theorem \ref{thm:main PBW thm}, $\Zhu_{g, H}V^k(\mathfrak{g})$ is PBW generated by $\projZhu(\mathfrak{g}^{g \circ \twist^{-1}}) \simeq  \mathfrak{g}^{g \circ \twist^{-1}} = \{a \in \mathfrak{g} \mid (g \circ \twist^{-1})(a) =a \}$ with defining relations
\begin{align*}
\bar{a} * \bar{b} - (-1)^{p(a)p(b)}\bar{b} * \bar{a} = \overline{[a, b]}- k(x|[a, b]),\quad
a, b \in \mathfrak{g}^{g \circ \twist^{-1}}.
\end{align*}
Here, we write $\bar{a}$ for $\projZhu(a)$. Since the defining relations above can be written as
\begin{align*}
\widetilde{a} * \widetilde{b} - (-1)^{p(a)p(b)}\widetilde{b} * \widetilde{a} = \widetilde{[a, b]},\quad
\widetilde{a} := a - k(x|a),
\end{align*}
we obtain the desired isomorphism.
\end{proof}
\end{theorem}

For example, choose a semisimple element $x$ such that
\begin{align*}
-1 < \alpha(x) < 0,\quad
\alpha(x) \in \Q,\quad
\alpha \in \Delta_+,
\end{align*}
where $\Delta_+$ is a set of all positive roots of $\mathfrak{g}$. Then the Hamiltonian operator $H$ defines a $\Q$-grading on $V^k(\mathfrak{g})$. By Theorem \ref{thm:twisted Zhu of affine},
\begin{align*}
\Zhu_{1, H}V^k(\mathfrak{g}) \simeq U(\mathfrak{h}),
\end{align*}
where $\mathfrak{h}$ is a Cartan subalgebra of $\mathfrak{g}$ containing $x$. Thus, any untwisted Zhu algebras of quotient vertex superalgebra of $V^k(\mathfrak{g})$ are commutative. This fact seems crucial to prove the $\omega_x$-rationality of simple affine vertex superalgebras at admissible levels in \cite{DLM97, Lin, LW}.

\subsection{Affine $W$-algebras}\label{sec:affine W}
We recall the definitions of the (affine) $W$-algebras associated with $\mathfrak{g}, f$ at level $k$, introduced by \cite{FF90a, KRW}.

\subsubsection{Good pair and automorphism $g$}\label{sec:good}
Let $\mathfrak{g}$ be a simple basic classical Lie superalgebra with the normalized invariant bilinear form $(\cdot|\cdot)$ on $\mathfrak{g}$, $f$ a nilpotent element in the even part $\mathfrak{g}_{\bar{0}}$ of $\mathfrak{g}$ and $x$ a semisimple element of $\mathfrak{g}$ such that (1) $\ad x$ defines a $\frac{1}{2}\Z$-grading
\begin{align*}
\mathfrak{g} = \bigoplus_{j \in \frac{1}{2}\Z}\mathfrak{g}_j,\quad
\mathfrak{g}_j = \{ a \in \mathfrak{g} \mid [x, a]=j a\},
\end{align*}
(2) $f \in \mathfrak{g}_{-1}$, and (3) $\ad f \colon \mathfrak{g}_j \rightarrow \mathfrak{g}_{j-1}$ is injective for $j \geq \frac{1}{2}$ and surjective for $j \leq \frac{1}{2}$. Then $(f, x)$ is called a good pair.

Let $g$ be a diagonalizable Lie superalgebra automorphism of $\mathfrak{g}$ with modulus $1$ eigenvalues such that
\begin{align*}
g(x) = x,\quad
g(f) = f,\quad
(g(a)|g(b)) = (a|b),\quad
a, b \in \mathfrak{g}.
\end{align*}
Then $g$ defines an automorphism of $V^k(\mathfrak{g})$. Moreover, $g$ commutes with a Hamiltonian operator $H=L^\mathfrak{g}_0 - x_{(0)}$ of $V^k(\mathfrak{g})$. Recall that $\theta_H = \mathrm{e}^{2\pi iH}$ defines an automorphism of $\mathfrak{g}$. Since $g$ commutes with $H$, we have the simeltaneous eigenspace decomposition
\begin{align*}
\mathfrak{g} = \bigoplus_{\bar{\gamma} \in \R/\Z, \Delta \in \frac{1}{2}\Z}\mathfrak{g}(\bar{\gamma}, \Delta),\quad
\mathfrak{g}(\bar{\gamma}, \Delta) = \{ a \in \mathfrak{g} \mid g(a) = \mathrm{e}^{2\pi i\gamma}a, H(a) = \Delta a\}.
\end{align*}
Then $\mathfrak{g}_j$ is invariant under the $g$-action, and $\mathfrak{g}_j = \bigoplus_{\bar{\gamma} \in \R/\Z}\mathfrak{g}(\bar{\gamma}, 1-j)$. Let
\begin{align*}
\mathfrak{g}(\bar{\gamma}, \bar{\mu}) := \bigoplus_{\Delta \in \bar{\mu}}\mathfrak{g}(\bar{\gamma}, \Delta),\quad
\bar{\mu} \in \frac{1}{2}\Z/\Z \simeq \Z_2.
\end{align*}
Then
\begin{align*}
\mathfrak{g}[\bar{\gamma}] := \bigoplus_{\bar{\mu}\in\Z_2}\mathfrak{g}(\bar{\gamma}+\bar{\mu}, \bar{\mu})
\end{align*}
is the eigenspace of $\mathfrak{g}$ for the action of $g \circ \theta_H^{-1}$ with the eigenvalue $\mathrm{e}^{2\pi i \gamma}$ for $\bar{\gamma} \in \R/\Z$. Note that
\begin{align*}
\mathfrak{g}[0] = \mathfrak{g}^{g \circ \theta_H^{-1}}.
\end{align*}

Let $\mathfrak{g}_j[\bar{\gamma}] = \mathfrak{g}_j \cap \mathfrak{g}[\bar{\gamma}]$. Fix a basis $\{ e_\alpha\}_{\alpha \in S_j[\bar{\gamma}]}$ of $\mathfrak{g}_j[\bar{\gamma}]$ with an index set $S_j[\bar{\gamma}]$. Set $S_j = \bigsqcup_{\bar{\gamma} \in \R/\Z}S_j[\bar{\gamma}]$, $S_{>0} = \bigsqcup_{j>0}S_j$, $S=\bigsqcup_{j \in \frac{1}{2}\Z}S_j$, $S_{>0}[\bar{\gamma}] = \bigsqcup_{j >0}S_{j}[\bar{\gamma}]$, $S[\bar{\gamma}] = \bigsqcup_{j  \in \frac{1}{2}\Z}S_{j}[\bar{\gamma}]$. By Theorem \ref{thm:twisted Zhu of affine},
\begin{align*}
\gH{\Zhu}V^k(\mathfrak{g}) \simeq U(\mathfrak{g}^{g \circ \theta_H^{-1}}),\quad
\tau(a) \mapsto \widetilde{a} = a -k(x|a),\quad
a \in \mathfrak{g}^{g \circ \theta_H}.
\end{align*}

\subsubsection{Vertex superalgebra $\Phi(\mathfrak{g}_{\frac{1}{2}})$}

Since $\ad f \colon \mathfrak{g}_{\frac{1}{2}} \xrightarrow{\sim} \mathfrak{g}_{-\frac{1}{2}}$ is an isomorphism,
\begin{align*}
\langle a, b \rangle_\mathfrak{n} := (f|[a, b]),\quad
a, b \in \mathfrak{g}_{\frac{1}{2}}
\end{align*}
defines a non-degenerate skew-supersymmetric bilinear form on $\mathfrak{g}_{\frac{1}{2}}$.

Let $\Phi(\mathfrak{g}_{\frac{1}{2}})$ be a vertex superalgebra defined as the universal enveloping vertex algebra of a non-linear Lie conformal superalgebra
\begin{align*}
\conf_\mathfrak{n} = \C[\partial]\otimes\mathfrak{n},\quad
\mathfrak{n} = \Span_\C\{\Phi_\alpha \mid \alpha \in S_{\frac{1}{2}}\},
\end{align*}
where $\Phi_\alpha$ has the same parity as $e_\alpha$, and the $\lambda$-bracket on $\conf_\mathfrak{n}$ is defined by
\begin{align*}
[{\Phi_\alpha}_\lambda \Phi_\beta] = \langle e_\alpha, e_\beta \rangle_\mathfrak{n} = (f|[e_\alpha, e_\beta]),\quad
\alpha, \beta \in S_{\frac{1}{2}}.
\end{align*}
Then we have an isomorphism of superspaces $\mathfrak{g}_{\frac{1}{2}} \ni e_\alpha \mapsto \Phi_\alpha \in \mathfrak{n}$. Denote by $\Phi_a$ an element in $\mathfrak{n}$ corresponding to $a \in \mathfrak{g}_{\frac{1}{2}}$. Since $\langle \cdot, \cdot \rangle_\mathfrak{n}$ is non-degenerate, we have a dual basis $\{e^\alpha\}_{\alpha \in S_{\frac{1}{2}}}$ of $\mathfrak{g}_{\frac{1}{2}}$ such that $\langle e_\alpha, e^\beta \rangle_\mathfrak{n} = \delta_{\alpha, \beta}$. Set $\Phi^\alpha=\Phi_{e^\alpha}$. Then
\begin{align*}
L^\mathfrak{n} = \frac{1}{2}\sum_{\alpha \in S_{\frac{1}{2}}}\NO{(\partial\Phi^\alpha)\Phi_\alpha}
\end{align*}
defines a conformal vector of $\Phi(\mathfrak{g}_{\frac{1}{2}})$ such that $L^\mathfrak{n}_0 (\Phi_\alpha) = \frac{1}{2}\Phi_\alpha$ for all $\alpha \in S_{\frac{1}{2}}$. Thus, $H=L^\mathfrak{n}_0$ defines a Hamiltonian operator of $\Phi(\mathfrak{g}_{\frac{1}{2}})$.

Define a $g$-action on $\conf_\mathfrak{n}$ by $g(\Phi_a) = \Phi_{g(a)}$. This is well-defined because $\mathfrak{g}_{\frac{1}{2}}$ is $g$-invariant, and
\begin{align*}
[{g(\Phi_a)}_\lambda g(\Phi_b)]
&= (f|[g(a), g(b)]) = (f|g([a,b])) = (g^{-1}(f)|[a,b])\\
&= (f|[a,b]) = g([{\Phi_a}_\lambda \Phi_b]).
\end{align*}
Thus, $g$ defines an automorphism of $\Phi(\mathfrak{g}_{\frac{1}{2}})$ commuting with $H$. By Theorem \ref{thm:main PBW thm},
\begin{align*}
\gH{\overline{\Phi}} := \gH{\Zhu}\Phi(\mathfrak{g}_{\frac{1}{2}})
\end{align*}
is PBW generated by $\overline{\Phi}_\alpha:=\tau(\Phi_\alpha)$ for $\alpha \in S_{\frac{1}{2}}[0]$ with the following defining relations:
\begin{align*}
[\overline{\Phi}_\alpha, \overline{\Phi}_\beta] = \langle e_\alpha, e_\beta\rangle_\mathfrak{n},\quad
\alpha, \beta \in S_{\frac{1}{2}}[0].
\end{align*}

\subsubsection{Vertex superalgebra $F(\mathfrak{g}_{>0})$}

Let $\langle \cdot, \cdot \rangle_\mathfrak{f} \colon \mathfrak{g}_{>0} \times \mathfrak{g}_{>0}^* \rightarrow \C$ be the natural pairing, that is, $\langle a, s \rangle_\mathfrak{f} = s(a)$ for $a \in \mathfrak{g}_{>0}$ and $s \in \mathfrak{g}_{>0}^*$. Extend the pairing to a bilinear form on $\mathfrak{g}_{>0}\oplus\mathfrak{g}_{>0}^*$ by
\begin{align*}
\langle a, b \rangle_\mathfrak{f} = \langle s, t \rangle_\mathfrak{f}=0,\quad
\langle s, a \rangle_\mathfrak{f} = (-1)^{p(a)}\langle a, s \rangle_\mathfrak{f} = (-1)^{p(a)}s(a)
\end{align*}
for $a, b \in \mathfrak{g}_{>0}$ and $s, t \in \mathfrak{g}_{>0}^*$. Define $e_\alpha^* \in \mathfrak{g}_{>0}^*$ by $\langle e_\alpha, e_\beta^* \rangle_\mathfrak{f} = \delta_{\alpha, \beta}$ for $\alpha, \beta \in S_{>0}$. Then $\{e_\alpha^*\}_{\alpha \in S_{>0}}$ forms a basis of $\mathfrak{g}_{>0}^*$. Let
\begin{align*}
\mathfrak{f} = \Span_\C\{\varphi_\alpha, \varphi^\alpha \mid \alpha \in S_{>0}\}
\end{align*}
be a superspace spanned by $\varphi_\alpha$ and $\varphi^\alpha$ ($\alpha \in S_{>0}$) which have the opposite parity to $e_\alpha$. Then, we have an isomorphism of superspaces
\begin{align*}
\Pi(\mathfrak{g}_{>0} \oplus \mathfrak{g}_{>0}^*) \xrightarrow{\sim} \mathfrak{f},\quad
e_\alpha \mapsto \varphi_\alpha,\quad
e_\alpha^* \mapsto \varphi^\alpha,\quad
\alpha \in S_{>0},
\end{align*}
where $\Pi$ denotes the parity change of a superspace. Denote by $\varphi_a$, $\varphi^s$ elements in $\mathfrak{f}$ corresponding to $a \in \mathfrak{g}_{>0}$, $s \in \mathfrak{g}_{>0}^*$, respectively. Let $F(\mathfrak{g}_{>0})$ be the universal enveloping vertex algebra of a non-linear Lie conformal superalgebra
\begin{align*}
\conf_\mathfrak{f} = \C[\partial]\otimes\mathfrak{f},
\end{align*}
with the $\lambda$-bracket
\begin{align*}
[{\varphi_a + \varphi^s}_\lambda \varphi_b + \varphi^t] 
= \langle a+s, b + t\rangle_\mathfrak{f} 
= t(a)+(-1)^{p(b)}s(b),\quad
a, b \in \mathfrak{g}_{>0},\quad
s, t \in \mathfrak{g}_{>0}^*.
\end{align*}
For $\alpha \in S$, denote by $j_\alpha \in \frac{1}{2}\Z$ the eigenvalue of $e_\alpha$ for $\ad x$. Then
\begin{align*}
L^\mathfrak{f} = -\sum_{\alpha \in S_{>0}}j_\alpha\NO{\varphi^\alpha\partial\varphi_\alpha} + \sum_{\alpha \in S_{>0}}(1-j_\alpha)\NO{(\partial\varphi^\alpha)\varphi_\alpha}
\end{align*}
defines a conformal vector of $F(\mathfrak{g}_{>0})$ such that $L^\mathfrak{f}_0\varphi_\alpha = (1-j_\alpha)\varphi_\alpha$ and $L^\mathfrak{f}_0\varphi^\alpha = j_\alpha\varphi^\alpha$ for all $\alpha \in S_{>0}$. Thus, $H = L^\mathfrak{f}_0$ defines a Hamitonian operator of  $F(\mathfrak{g}_{>0})$.

Define a $g$-action on $\mathfrak{g}_{>0}^*$ by
\begin{align*}
\langle a, g(s) \rangle_\mathfrak{f} := \langle g^{-1}(a), s\rangle_\mathfrak{f},\quad
a \in \mathfrak{g}_{>0},\quad
s \in \mathfrak{g}_{>0}^*.
\end{align*}
Then a $g$-action on $\conf_\mathfrak{f}$ is defined by $g(\varphi_a) = \varphi_{g(a)}$ and $g(\varphi^s) = \varphi^{g(s)}$. This is well-defined because
\begin{align*}
[{g(\varphi_a + \varphi^s)}_\lambda g(\varphi_b + \varphi^t)]
&= \langle g(a)+g(s), g(b) + g(t)\rangle_\mathfrak{f} 
= \langle a+s, b + t\rangle_\mathfrak{f}\\
&= g([{\varphi_a + \varphi^s}_\lambda \varphi_b + \varphi^t])
\end{align*}
for $a, b \in \mathfrak{g}_{>0}$, $s, t \in \mathfrak{g}_{>0}^*$. Thus, $g$ defines an automorphism of $F(\mathfrak{g}_{>0})$ commuting with $H$. By Theorem \ref{thm:main PBW thm},
\begin{align*}
\gH{\overline{F}} := \gH{\Zhu}F(\mathfrak{g}_{>0})
\end{align*}
is PBW generated by $\overline{\varphi}_\alpha:=\tau(\varphi_\alpha)$ and $\overline{\varphi}^\alpha:=\tau(\varphi^\alpha)$ for $\alpha \in S_{>0}[0]$ with the following defining relations:
\begin{align*}
[\overline{\varphi}_\alpha, \overline{\varphi}_\beta] = 
[\overline{\varphi}^\alpha, \overline{\varphi}^\beta] = 0,\quad
[\overline{\varphi}_\alpha, \overline{\varphi}^\beta] = \delta_{\alpha, \beta},\quad
\alpha, \beta \in S_{>0}[0].
\end{align*}

\subsubsection{Definitions of affine $W$-algebras}

Let $C^k(\mathfrak{g}, f, x)$ be a vertex superalgebra defined by
\begin{align*}
C^k(\mathfrak{g}, f, x) = V^k(\mathfrak{g}) \otimes \Phi(\mathfrak{g}_{\frac{1}{2}}) \otimes F(\mathfrak{g}_{>0})
\end{align*}
and an odd element $Q$ in $C^k(\mathfrak{g}, f, x)$ defined by
\begin{align*}
Q =
&\sum_{\alpha \in S_{>0}}(-1)^{p(\alpha)}e_\alpha \varphi^\alpha - \frac{1}{2}\sum_{\alpha, \beta, \gamma \in S_{>0}}(-1)^{p(\alpha)p(\gamma)}c_{\alpha, \beta}^\gamma \NO{\varphi_\gamma\varphi^\alpha\varphi^\beta}\\
&+ \sum_{\alpha \in S_{\frac{1}{2}}}\Phi_\alpha\varphi^\alpha + \sum_{\alpha \in S_{>0}}(f|e_\alpha)\varphi^\alpha. 
\end{align*}
Since $[Q_\lambda Q] = 0$, we have $(Q_{(0)})^2 = 0$ so that $d = Q_{(0)}$ defines an odd derivation on $C^k(\mathfrak{g}, f)$. Define a charge degree on $C^k(\mathfrak{g}, f)$ by $\deg a = \deg \Phi_\alpha = 0$ for $a \in \mathfrak{g}$ and $\alpha \in S_{\frac{1}{2}}$, $\deg \varphi^\alpha = - \deg \varphi_\alpha = 1$ for $\alpha \in S_{>0}$, $\deg\partial a = \deg a$ and $\deg\NO{a b} = \deg a + \deg b$ for $a, b \in C^k(\mathfrak{g}, f, x)$. Then $C^k(\mathfrak{g}, f, x) = \bigoplus_{n \in \Z}C^k(\mathfrak{g}, f, x)^n$ with $C^k(\mathfrak{g}, f, x)^n = \{a \in C^k(\mathfrak{g}, f, x) \mid \deg a = n\}$ and $d \cdot C^k(\mathfrak{g}, f, x)^n \subset C^k(\mathfrak{g}, f, x)^{n+1}$. Thus $(C^k(\mathfrak{g}, f, x)^\bullet, d)$ forms a cochain complex. The affine $W$-algebra of $\mathfrak{g}, f$ at level $k \in \C$ is defined as the cohomology of the complex
\begin{align*}
W^k(\mathfrak{g}, f) := H^\bullet(C^k(\mathfrak{g}, f, x), d).
\end{align*}
It follows from \cite{KW04} that $H^i(C^k(\mathfrak{g}, f, x), d) = 0$ for $i \neq 0$.

Provided that $k \neq -h^\vee$,
\begin{align*}
L = L^\mathfrak{g} + \partial x + L^\mathfrak{n} + L^\mathfrak{f}
\end{align*}
defines a conformal vector of $C^k(\mathfrak{g}, f, x)$ such that  $\Delta_a = 1-j$ for $a \in \mathfrak{g}_j$, $\Delta_{\Phi_\alpha} = \frac{1}{2}$ for $\alpha \in S_{\frac{1}{2}}$ and $\Delta_{\varphi_\alpha} = 1 - j_\alpha$, $\Delta_{\varphi^\alpha} = j_\alpha$ for $\alpha \in S_{>0}$. Thus, we will extend a Hamiltonian operator $H$ on $V^k(\mathfrak{g})$ to $C^k(\mathfrak{g}, f, x)$ by 
\begin{align}\label{eq:H-cano}
H = L^\mathfrak{g}_0 - x_{(0)} + L^\mathfrak{n}_0 + L^\mathfrak{f}_0.
\end{align}
The Hamiltonian operator $H$ is well-defined for all $k \in \C$. Since $d \cdot L = 0$, $L$ also defines a conformal vector of $\mathcal{W}^k(\mathfrak{g}, f, x)$ if $k + h^\vee \neq 0$. Thus, $H$ gives a Hamiltonian operator of $\mathcal{W}^k(\mathfrak{g}, f)$ for all $k \in \C$. Let
\begin{align*}
J_a = a + \sum_{\beta, \gamma \in S_{>0}}(-1)^{p(\gamma)}c_{a, \beta}^\gamma \NO{\varphi_\gamma \varphi^\beta},\quad
a \in \mathfrak{g}
\end{align*}
be fields on $C^k(\mathfrak{g}, f, x)$, where $c_{a, \beta}^\gamma$ is the structure constant defined by $[a, e_\beta] = \sum_{\gamma \in S} c_{a, \beta}^\gamma e_\gamma$. Then
\begin{align*}
[J_a{}_\lambda J_b] = J_{[a, b]} 
+ \nu_k(a|b)\lambda,\quad
\nu_k(a|b) = 
\left(
k(a|b)+\frac{1}{2}\kappa_\mathfrak{g}(a|b)-\frac{1}{2}\kappa_{\mathfrak{g}_0}(a|b)
\right)
\end{align*}
for $a, b \in \mathfrak{g}_{\geq0}$ or $a, b \in \mathfrak{g}_{\leq0}$, where $\kappa_{\mathfrak{g}}$ denotes the Killing form on $\mathfrak{g}$. Note that $\Delta_{J_a} = 1 - j$ for $a \in \mathfrak{g}_j$. Let $C^-$ be a vertex subalgebra of $C^k(\mathfrak{g}, f, x)$ generated by $\varphi_\alpha$ and $d \cdot \varphi_\alpha = J_{e_\alpha} + \Phi_\alpha + (f|e_\alpha)$ for $\alpha \in \Delta_{>0}$, where $\Phi_\alpha = 0$ for $\alpha \in S_{\geq1}$. We have
\begin{align*}
[\varphi_\alpha{}_\lambda (d \cdot \varphi_\beta)] = \sum_{\gamma \in S_{>0}}c_{\alpha, \beta}^\gamma \varphi_\gamma,\quad
[(d \cdot \varphi_\alpha)_\lambda (d \cdot \varphi_\beta)] = \sum_{\gamma \in S_{>0}}c_{\alpha, \beta}^\gamma (d \cdot \varphi_\gamma)
\end{align*}
for $\alpha, \beta \in \Delta_{>0}$.  Let
\begin{align*}
\conf^- = \C[\partial]\otimes\base^-,\quad
\base^- = \Span_\C\{\varphi_\alpha,\ d \cdot \varphi_\alpha \mid \alpha \in S_{>0}\}.
\end{align*}
Using the $\lambda$-brackets above, $\conf^-$ forms a Lie conformal superalgebra whose universal enveloping vertex algebra is isomorphic to $C^-$, and $d$ defines an odd derivation of $\conf^-$. The cohomology $H(\conf^-, d)$ equals zero. By \cite[Lemma 3.4]{KW04}, $H(C^-, d) \simeq V(H(\conf^-, d)) = \C$. Let $C^+$ be a vertex subalgebra of $C^k(\mathfrak{g}, f, x)$ generated by $J_a$ for $a \in \mathfrak{g}_{\leq0}$, $\Phi_\alpha$ for $\alpha \in S_{\frac{1}{2}}$ and $\varphi^\alpha$ for $\alpha \in S_{>0}$. We have
\begin{align*}
[\varphi^\alpha{}_\lambda J^a] = \sum_{\beta \in S_{>0}}c_{a, \beta}^\alpha\varphi^\beta,\quad
[\Phi_\alpha{}_\lambda J^a] = 0,\quad
a \in \mathfrak{g}_{\leq0},
\end{align*}
\begin{align*}
[Q_\lambda J^a] =
&\Biggl\{ \sum_{\beta\in S_{>0}} \Bigl( k(a|e_\beta)+\operatorname{str}_\mathfrak{g}((\operatorname{ad}a)p_+(\operatorname{ad}e_\beta)) \Bigr)\partial\varphi^\beta\\
&+ \sum_{\beta \in S_{>0}} :\Bigl( -\sum_{\alpha \in S_{\leq 0}} (-1)^{p(\alpha)}c_{a, \beta}^\alpha J^{e_\alpha} + \sum_{\alpha \in S_\frac{1}{2}}c_{a, \beta}^\alpha\Phi_\alpha + (f|[a, e_\beta]) \Bigr)\varphi^\beta: \Biggr\}\\
&+\lambda\sum_{\beta\in S_{>0}}\Bigl( k(a|e_\beta)+\operatorname{str}_\mathfrak{g}(p_+(\operatorname{ad}a)(\operatorname{ad}e_\beta)) \Bigr)\varphi^\beta
\end{align*}
for $a \in \mathfrak{g}_{\leq 0}$, where $p_+ \colon \mathfrak{g} \twoheadrightarrow \mathfrak{g}_{>0}$ is a projection, and
\begin{align*}
&[Q_\lambda\Phi_\alpha] = \sum_{\beta \in S_\frac{1}{2}}(f|[e_\beta, e_\alpha])\varphi^\beta,\quad
\alpha \in S_\frac{1}{2},\\
&[Q_\lambda\varphi^\alpha] = -\frac{1}{2}\sum_{\beta, \gamma \in S_{>0}}(-1)^{p(\alpha)p(\beta)}c_{\beta, \gamma}^\alpha:\varphi^\beta\varphi^\gamma:,\quad
\alpha \in S_{>0}.
\end{align*}
Since $C^k(\mathfrak{g}, f, x) \simeq C^+ \otimes C^-$ as super spaces and $(C^\pm, d)$ are subcomplexes of $(C^k(\mathfrak{g}, f, x), d)$,
\begin{align*}
H(C^k(\mathfrak{g}, f, x), d) = H(C^+, d) \otimes H(C^-, d) = H(C^k_+, d).
\end{align*}
Let
\begin{align*}
\conf^+ = \C[\partial]\otimes\base^+,\quad
\base^+ = \Span_\C\{J_{e_\alpha},\ \Phi_\beta,\ \varphi^\gamma \mid \alpha \in S_{\leq0},\ \beta \in S_{\frac{1}{2}},\ \gamma \in S_{>0}\}.
\end{align*}
Then $\conf^+$ forms a non-linear Lie conformal superalgebra whose universal enveloping vertex algebra is isomorphic to $C^+$. We have \cite[Theorem 5.7]{DK}:
\begin{align*}
W^k(\mathfrak{g}, f) \simeq H(C^+, d).
\end{align*}

\subsubsection{Twisted Zhu algebra of $C^k(\mathfrak{g}, f, x)$}
Using the facts that $\langle g(e_\alpha), g(e_\beta^*)\rangle_\mathfrak{f} = \delta_{\alpha, \beta}$ and $[g(e_\alpha), g(e_\beta)] = g([e_\alpha, e_\beta])$, one can show that
\begin{align*}
g(Q) = Q.
\end{align*}
Then $g(d \cdot a) = g(Q_{(0)}a) = g(Q)_{(0)}g(a) = Q_{(0)}g(a)$ for $a \in C^k(\mathfrak{g}, f, x)$. Hence $g$ commutes with $d$. Moreover, $H$ commutes with $d$. Therefore, $g$ defines an automorphism of $W^k(\mathfrak{g}, f)$ commuting $H$. Now one has
\begin{align*}
\gH{\Zhu}W^k(\mathfrak{g}, f) = \gH{\Zhu}H^\bullet(C^k(\mathfrak{g}, f, x) , d).
\end{align*}

On the other hand, $d$ induces an odd differential on
\begin{align*}
\gH{\overline{C}} := \gH{\Zhu}C^k(\mathfrak{g}, f, x) \simeq U(\mathfrak{g}^{g \circ \theta_H^{-1}}) \otimes \gH{\overline{\Phi}} \otimes \gH{\overline{F}},
\end{align*}
which we denote by $\overline{d}$. By direct computations,
\begin{align*}
\overline{Q}:= \tau(Q) = 
&\sum_{\alpha \in S_{>0}[0]}(-1)^{p(\alpha)}\widetilde{e}_\alpha * \overline{\varphi}^\alpha - \frac{1}{2}\sum_{\alpha, \beta, \gamma \in S_{>0}[0]}(-1)^{p(\alpha)p(\gamma)}c_{\alpha, \beta}^\gamma \overline{\varphi}_\gamma * \overline{\varphi}^\alpha * \overline{\varphi}^\beta\\
&+ \sum_{\alpha \in S_{\frac{1}{2}}[0]}\overline{\Phi}_\alpha*\overline{\varphi}^\alpha + \sum_{\alpha \in S_{>0}[0]}(f|e_\alpha)\overline{\varphi}^\alpha.
\end{align*}
Since $[\overline{Q}, \overline{a}] = [\tau(Q)_*\tau(a)] = \tau(Q_{(0)}a) = \overline{d} \cdot \overline{a}$ for $\overline{a} \in \gH{\overline{C}}$, we obtain that $\overline{d} = \ad \overline{Q}$. Therefore,
\begin{lemma}\label{lem: finite-W}
\begin{align*}
U(\mathfrak{g}^{g \circ \theta_H^{-1}}, f) := H^\bullet(\gH{\overline{C}}, \overline{d})
\end{align*}
is isomorphic to the finite $W$-algebra of $\mathfrak{g}^{g \circ \theta_H^{-1}}$ and $f$.
\end{lemma}

\subsubsection{Twisted Zhu algebras of $W$-algebras}
For $a \in \mathfrak{g}_{\leq 0}[0]$, set
\begin{align*}
\widetilde{J}_a := \tau(J_a) - \nu_k(x|a) = \widetilde{a} + \sum_{\beta, \gamma \in S_{>0}[0]}(-1)^{p(\gamma)}c_{a, \beta}^\gamma \overline{\varphi}_\gamma * \overline{\varphi}^\beta \in \gH{\overline{C}}.
\end{align*}
Then $[\widetilde{J}_a, \widetilde{J}_b] = \widetilde{J}_{[a, b]}$. Let $\gH{\overline{C}}^-$ be a subalgebra of $\gH{\overline{C}}$ generated by $\overline{\varphi}_\alpha$ and $\overline{d} \cdot \overline{\varphi}_\alpha$ for $\alpha \in S_{>0}[0]$, and $\gH{\overline{C}}^+$ a subalgebra of $\gH{\overline{C}}$ generated by $\overline{J}_{e_\alpha}, \overline{\Phi}_\beta, \overline{\varphi}_\gamma$ for $\alpha \in S_{\leq0}[0], \beta \in S_{\frac{1}{2}}[0], \gamma \in S_{>0}[0]$. Using the same arguments as in \cite[Section 5.4]{DK}, $(\gH{\overline{C}}^\pm, \overline{d})$ are subcomplexes of $\gH{\overline{C}}$ such that
\begin{align*}
\gH{\overline{C}} \simeq \gH{\overline{C}}^+ \otimes \gH{\overline{C}}^-
\end{align*}
as superspaces, and
\begin{align*}
H^i(\gH{\overline{C}}, \overline{d}) \simeq H^i(\gH{\overline{C}}^+, \overline{d}),\quad
\gH{\Zhu}C^+ \simeq \gH{\overline{C}}^+.
\end{align*}

Set $\monoid = \frac{1}{2}\Z$ and $\Charge = \{ (p, q) \in \frac{1}{2}\Z \times \frac{1}{2}\Z \mid p + q \in \Z_{\geq0}\}$. Finally, as in \cite[Section 5.5]{DK}, define a $\monoid_{>0} \times \Charge \times \R$-grading on $\base^+$ by
\begin{align*}
&\zeta_{J_a} = 1 - j,\quad
\deg_\Charge J_a = \left(j - \frac{1}{2}, -j+\frac{1}{2}\right),\quad
\Delta_{J_a} = 1 - j,\quad
a \in \mathfrak{g}_j,\\
&\zeta_{\Phi_\alpha} = \frac{1}{2},\quad
\deg_\Charge \Phi_\alpha = \left(0, 0\right),\quad
\Delta_{J^a} = \frac{1}{2},\quad
\alpha \in S_{\frac{1}{2}},\\
&\zeta_{\varphi^\alpha} = j_\alpha,\quad
\deg_\Charge \varphi^\alpha = \left(-j_\alpha + \frac{1}{2}, j_\alpha+\frac{1}{2}\right),\quad
\Delta_{\varphi^\alpha} = j_\alpha,\quad
\alpha \in S_{>0}.
\end{align*}

\begin{theorem}\label{thm:ZhuofW}
Let $g$ be a diagonalizable Lie superalgebra automorphism of $\mathfrak{g}$ with modulus $1$ eigenvalues such that $g(x) = x$, $g(f) = f$ and $(g(a)|g(b)) = (a|b)$ for $a, b \in \mathfrak{g}$, which induces an automorphism of $W^k(\mathfrak{g}, f)$ commuting with the natural Hamiltonian operator $H$. Then
\begin{align*}
\gH{\Zhu}W^k(\mathfrak{g}, f) \simeq U(\mathfrak{g}^{g \circ \theta_H^{-1}}, f).
\end{align*}
\begin{proof}
Using the same arguments as in \cite[Section 5.5]{DK}, one can check that $W^k(\mathfrak{g}, f) \simeq H^\bullet(C^+, d)$ satisfies all assumptions in Theorem \ref{thm:main_cohomology_thm}. Therefore,
\begin{align*}
\gH{\Zhu}W^k(\mathfrak{g}, f) \simeq H^0(\gH{\Zhu}C^+, \overline{d}) \simeq H^0(\gH{\overline{C}}^+, \overline{d}) \simeq U(\mathfrak{g}^{g \circ \theta_H^{-1}}, f)
\end{align*}
by Lemma \ref{lem: finite-W}. This completes the proof.
\end{proof}
\end{theorem}

\begin{corollary}\label{cor: g=1-finite-W}
\begin{align*}
\Zhu_{1, H}W^k(\mathfrak{g}, f) \simeq U(\mathfrak{g}_\Z, f),
\end{align*}
where $\mathfrak{g}_\Z = \bigoplus_{j \in \Z}\mathfrak{g}_j$.
\begin{proof}
Since $\Delta_a = 1 - j$ for $a \in \mathfrak{g}_j$ and
\begin{align*}
\theta_H(a) = \mathrm{e}^{2\pi i \Delta_a}a =
\begin{cases}
a & \mathrm{if}\ \Delta_a \in \Z,\\
-a & \mathrm{if}\ \Delta_a \in \Z + \frac{1}{2},
\end{cases}
\end{align*}
we have $\mathfrak{g}^{g \circ \theta_H^{-1}} = \mathfrak{g}^{\theta_H} = \mathfrak{g}_\Z$ if $g=1$. Apply Theorem \ref{thm:ZhuofW} for $g=1$.
\end{proof}
\end{corollary}

\bibliographystyle{halpha} 
\bibliography{refs} 

\end{document}